\documentclass[11pt]{article}
\usepackage{amsmath,amsthm,amsfonts,mathrsfs,bbm}
\usepackage{graphicx,psfrag}
\usepackage{url}

\usepackage{a4}
\usepackage{fullpage}
\newcommand{\N}{\mathbb{N}}

\newcommand{\E}{\mathbb{E}}
\newcommand{\LL}{\mathscr{L}}

\newcommand{\1}{\mathbbm{1}}

\newcommand{\da}{\downarrow}

\newcommand{\tr}{{\rm tr}}

\renewcommand{\Pr}{\mathbb{P}}

 \makeatletter
 \@addtoreset{equation}{section}
 \makeatother

 \makeatletter
 \@addtoreset{enunciato}{section}
 \makeatother

 \newcounter{enunciato}[section]

 \newtheorem{ittheorem}{Theorem}
 \newtheorem{itlemma}{Lemma}
 \newtheorem{itproposition}{Proposition}
 \newtheorem{itdefinition}{Definition}
 \newtheorem{itremark}{Remark}
 \newtheorem{itclaim}{Claim}
 \newtheorem{itfact}{Fact}
 \newtheorem{itconjecture}{Conjecture}
 \newtheorem{itobservation}{Observation} 
 \newtheorem{itcorollary}{Corollary}

 \newenvironment{theorem}{\addtocounter{enunciato}{1}
 \begin{ittheorem}}{\end{ittheorem}}

 \newenvironment{lemma}{\addtocounter{enunciato}{1}
 \begin{itlemma}}{\end{itlemma}}

 \newenvironment{proposition}{\addtocounter{enunciato}{1}
 \begin{itproposition}}{\end{itproposition}}

 \newenvironment{definition}{\addtocounter{enunciato}{1}
 \begin{itdefinition}}{\end{itdefinition}}

 \newenvironment{remark}{\addtocounter{enunciato}{1}
 \begin{itremark}}{\end{itremark}}

 \newenvironment{claim}{\addtocounter{enunciato}{1}
 \begin{itclaim}}{\end{itclaim}}

 \newenvironment{fact}{\addtocounter{enunciato}{1}
 \begin{itfact}}{\end{itfact}}

 \newenvironment{conjecture}{\addtocounter{enunciato}{1}
 \begin{itconjecture}}{\end{itconjecture}}

 \newenvironment{corollary}{\addtocounter{enunciato}{1}
 \begin{itcorollary}}{\end{itcorollary}}

 \newcommand{\be}[1]{\begin{equation}\label{#1}}
 \newcommand{\ee}{\end{equation}}

 \newcommand{\bl}[1]{\begin{lemma}\label{#1}}
 \newcommand{\el}{\end{lemma}}

 \newcommand{\br}[1]{\begin{remark}\label{#1}}
 \newcommand{\er}{\end{remark}}

 \newcommand{\bt}[1]{\begin{theorem}\label{#1}}
 \newcommand{\et}{\end{theorem}}

 \newcommand{\bd}[1]{\begin{definition}\label{#1}}
 \newcommand{\ed}{\end{definition}}

 \newcommand{\bcl}[1]{\begin{claim}\label{#1}}
 \newcommand{\ecl}{\end{claim}}

 \newcommand{\bfact}[1]{\begin{fact}\label{#1}}
 \newcommand{\efact}{\end{fact}}

 \newcommand{\bp}[1]{\begin{proposition}\label{#1}}
 \newcommand{\ep}{\end{proposition}}

 \newcommand{\bc}[1]{\begin{corollary}\label{#1}}
 \newcommand{\ec}{\end{corollary}}

 \newcommand{\bcj}[1]{\begin{conjecture}\label{#1}}
 \newcommand{\ecj}{\end{conjecture}}

 \newcommand{\bpr}{\begin{proof}}
 \newcommand{\epr}{\end{proof}}

 \newcommand{\bprl}[1]{\begin{proofof}{\it\ref{#1}}.\,\,}
 \newcommand{\eprl}{\end{proofof}}

 \newcommand{\bi}{\begin{itemize}}
 \newcommand{\ei}{\end{itemize}}

 \newcommand{\ben}{\begin{enumerate}}
 \newcommand{\een}{\end{enumerate}}

\title{Quenched large deviation principle for words in a letter sequence}

\date{13th May 2009}

\author{

Matthias Birkner
\footnotemark[1]%
\\

Andreas Greven
\footnotemark[2]
\\

Frank den Hollander
\footnotemark[3]\,\, \footnotemark[4]
}

\footnotetext[1]{%
Department Biologie II, Abteilung Evolutionsbiologie,
University of Munich (LMU),
Grosshaderner Str. 2,
82152 Planegg-Martinsried,
Germany}

\footnotetext[2]{
Mathematisches Institut, Universit\"at Erlangen-N\"urnberg,
Bismarckstrasse $1\frac12$, 91054 Erlangen, Germany}

\footnotetext[3]{
Mathematical Institute, Leiden University, P.O.\ Box 9512,
2300 RA Leiden, The Netherlands}

\footnotetext[4]{
EURANDOM, P.O.\ Box 513, 5600 MB Eindhoven, The Netherlands}

\begin{document} 

\maketitle

\begin{abstract} 
When we cut an i.i.d.\ sequence of letters into words according to
an independent renewal process, we obtain an i.i.d.\ sequence of
words. In the \emph{annealed} large deviation principle (LDP) for
the empirical process of words, the rate function is the specific
relative entropy of the observed law of words w.r.t.\ the reference
law of words. In the present paper we consider the \emph{quenched}
LDP, i.e., we condition on a typical letter sequence. We focus on
the case where the renewal process has an \emph{algebraic} tail. The
rate function turns out to be a sum of two terms, one being the
annealed rate function, the other being proportional to the specific
relative entropy of the observed law of letters w.r.t.\ the reference 
law of letters, with the former being obtained by concatenating the 
words and randomising the location of the origin.  The proportionality 
constant equals the tail exponent of the renewal process. Earlier work 
by Birkner considered the case where the renewal process has an exponential 
tail, in which case the rate function turns out to be the first term on 
the set where the second term vanishes and to be infinite elsewhere.
In a companion paper the annealed and the quenched LDP are applied to 
the collision local time of transient random walks, and the existence 
of an intermediate phase for a class of interacting stochastic systems 
is established. 

\vspace{0.5cm}
\emph{Acknowledgement:} This work was supported in part by DFG and NWO 
through the Dutch-German Bilateral Research Group ``Mathematics of Random 
Spatial Models from Physics and Biology''. MB and AG are grateful for
hospitality at EURANDOM. 
We also thank an anonymous referee for her/his careful reading and 
helpful comments. 
\end{abstract}

\newpage


\section{Introduction and main results}
\label{S1} 

\subsection{Problem setting} 
\label{S1.1}
Let $E$ be a finite set of \emph{letters}. Let $\widetilde{E} =
\cup_{n\in\N} E^n$ be the set of finite \emph{words} drawn from
$E$. Both $E$ and $\widetilde{E}$ are Polish spaces under the discrete
topology. Let $\mathcal{P}({E}^\N)$ and $\mathcal{P}(\widetilde{E}^\N)$ 
denote the set of probability measures on sequences drawn from $E$,
respectively, $\widetilde{E}$, equipped with the topology of weak convergence.  
Write $\theta$ and $\widetilde{\theta}$ for the left-shift acting on $E^\N$, 
respectively, $\widetilde{E}^\N$. Write $\mathcal{P}^{\mathrm{inv}}(E^\N),
\mathcal{P}^{\mathrm{erg}} (E^\N)$ and $\mathcal{P}^{\mathrm{inv}}
(\widetilde{E}^\N),\mathcal{P}^{\mathrm{erg}}(\widetilde{E}^\N)$ for 
the set of probability measures that are invariant and ergodic under 
$\theta$, respectively, $\widetilde{\theta}$.

For $\nu \in \mathcal{P}(E)$, let $X=(X_i)_{i\in\N}$ be i.i.d.\ with law 
$\nu$. Without loss of generality we will assume that ${\rm supp}(\nu)=E$
(otherwise we replace $E$ by ${\rm supp}(\nu)$). For $\rho \in \mathcal{P}(\N)$, 
let $\tau=(\tau_i)_{i\in\N}$ be i.i.d.\ with law $\rho$ having infinite 
support and satisfying the \emph{algebraic tail property}
\be{rhocond}
\lim_{ {n\to\infty} \atop {\rho(n)>0} } \frac{\log\rho(n)}{\log n}  
=: -\alpha, \quad \alpha \in (1,\infty). 
\ee 
(No regularity assumption will be necessary for ${\rm supp}(\rho)$.) Assume 
that $X$ and $\tau$ are independent and write $\Pr$ to denote their joint law. 
Cut words out of $X$ according to 
$\tau$, i.e., put (see Figure 1)
\be{Tdefs}
T_0:=0 \quad \mbox{ and } \quad T_i:=T_{i-1}+\tau_i,\quad i\in\N,
\ee 
and let
\be{eqndefYi} 
Y^{(i)} := \bigl( X_{T_{i-1}+1}, X_{T_{i-1}+2},\dots, X_{T_{i}}\bigr), 
\quad i \in \N. 
\ee 
Then, under the law $\Pr$, $Y = (Y^{(i)})_{i\in\N}$ is an i.i.d.\ sequence of 
words with marginal law $q_{\rho,\nu}$ on $\widetilde{E}$ given by 
\be{q0def}
\begin{aligned}
&q_{\rho,\nu}\big((x_1,\dots,x_n)\big) := \Pr\big(Y^{(1)}=(x_1,\dots,x_n)\big) 
= \rho(n) \,\nu(x_1) \cdots \nu(x_n),\\
&n\in\N,\,x_1,\dots,x_n\in E. 
\end{aligned}
\ee  

\begin{figure}[hctb] 
\label{cutwords}
\begin{center}
\psfrag{tau1}{$\tau_1$}
\psfrag{tau2}{$\tau_2$}
\psfrag{tau3}{$\tau_3$}
\psfrag{tau4}{$\tau_4$}
\psfrag{tau5}{$\tau_5$}
\psfrag{T1}{$T_1$}
\psfrag{T2}{$T_2$}
\psfrag{T3}{$T_3$}
\psfrag{T4}{$T_4$}
\psfrag{T5}{$T_5$}
\psfrag{Y1}{$Y^{(1)}$}
\psfrag{Y2}{$Y^{(2)}$}
\psfrag{Y3}{$Y^{(3)}$}
\psfrag{Y4}{$Y^{(4)}$}
\psfrag{Y5}{$Y^{(5)}$}
\psfrag{X}{$X$}
\includegraphics[width=120mm,height=22mm]{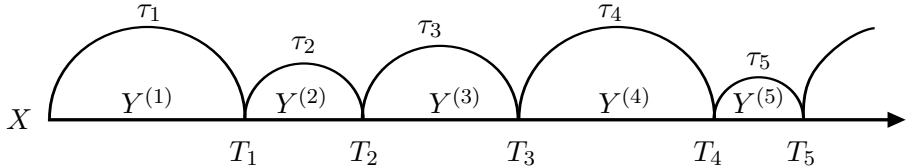}
\end{center}
\caption{Cutting words from a letter sequence according to a renewal process.}
\end{figure}

For $N\in\N$, let $(Y^{(1)},\dots,Y^{(N)})^\mathrm{per}$ stand for the periodic 
extension of $(Y^{(1)},\dots,Y^{(N)})$ to an element of $\widetilde{E}^\N$, and 
define
\be{eqndefRN} 
R_N := \frac{1}{N} \sum_{i=0}^{N-1} 
\delta_{\widetilde{\theta}^i (Y^{(1)},\dots,Y^{(N)})^\mathrm{per}} 
\; \in \mathcal{P}^{\mathrm{inv}}(\widetilde{E}^\N),
\ee 
the \emph{empirical process of $N$-tuples of words}. By the ergodic theorem, we have 
\be{ergthm}
{\rm w}\!-\!\!\lim_{N\to\infty} R_N = q_{\rho,\nu}^{\otimes\N} \quad\mbox{$\Pr$--a.s.},
\ee
with ${\rm w}-\lim$ denoting the weak limit. The following large deviation principle
(LDP) is standard (see e.g.\ Dembo and Zeitouni~\cite{DeZe98}, Corollaries~6.5.15 
and 6.5.17). For $Q \in \mathcal{P}^{\mathrm{inv}}(\widetilde{E}^\N)$ let
\be{spentrdef}
H(Q \mid q_{\rho,\nu}^{\otimes\N}) 
:= \lim_{N\to\infty} \frac{1}{N}\,h\left(Q_{|_{\mathscr{F}_N}}
\mid (q_{\rho,\nu}^{\otimes\N})_{|_{\mathscr{F}_N}}\right) \in [0,\infty]
\ee
be the \emph{specific relative entropy of $Q$ w.r.t.\ $q_{\rho,\nu}^{\otimes\N}$}, where 
$\mathscr{F}_N=\sigma(Y^{(1)},\dots,Y^{(N)})$ is the sigma-algebra generated by the first 
$N$ words, $Q_{|_{\mathscr{F}_N}}$ is the restriction of $Q$ to $\mathscr{F}_N$, and 
$h(\,\cdot\mid\cdot\,)$ denotes relative entropy. (For general properties of entropy, 
see Walters~\cite{Wa82}, Chapter~4.)

\begin{theorem}
\label{LDPann}
{\bf [Annealed LDP]}
The family of probability distributions $\Pr(R_N \in \cdot\,)$, $N\in\N$, satisfies 
the LDP on $\mathcal{P}^{\mathrm{inv}}(\widetilde{E}^\N)$ with rate $N$ and with rate 
function $I^\mathrm{ann}\colon\,\mathcal{P}^{\mathrm{inv}}(\widetilde{E}^\N) \to [0,\infty]$ 
given by
\be{Iann}
I^\mathrm{ann}(Q) = H(Q \mid q_{\rho,\nu}^{\otimes\N}).
\ee 
This rate function is lower semi-continuous, has compact level sets, has a unique
zero at $Q=q_{\rho,\nu}^{\otimes\N}$, and is affine.
\end{theorem}

\noindent
The LDP for $R_N$ arises from the LDP for $N$-tuples via a projective limit theorem. The 
ratio under the limit in (\ref{spentrdef}) is the rate function for $N$-tuples according 
to Sanov's theorem (see e.g.\ den Hollander~\cite{dHo00}, Section II.5), and is non-decreasing
in $N$.

\subsection{Main theorems} 
\label{S1.2}

Our aim in the present paper is to derive the LDP for $\Pr(R_N \in\cdot\mid X)$, $N\in\N$. 
To state our result, we need some more notation.

Let $\kappa\colon\widetilde{E}^\N\to E^\N$ denote the \emph{concatenation map} that glues 
a sequence of words into a sequence of letters. For $Q \in \mathcal{P}^{\mathrm{inv}}
(\widetilde{E}^\N)$ such that 
\be{eq:defmQ}
m_Q := \E_Q[\tau_1] < \infty, 
\ee 
define $\Psi_Q \in \mathcal{P}^{\mathrm{inv}}({E}^\N)$ as
\be{PsiQdef}
\Psi_Q(\cdot) := \frac{1}{m_Q} \E_Q\left[ \sum_{k=0}^{\tau_1-1} 
\delta_{\theta^k\kappa(Y)}(\cdot)\right].
\ee
Think of $\Psi_Q$ as \emph{the shift-invariant version of the concatenation of $Y$ under
the law $Q$ obtained after randomising the location of the origin}. 

For $\tr\in\N$, let $[\cdot]_\tr\colon\,\widetilde{E} \to [\widetilde{E}]_\tr := \cup_{n=1}^\tr 
E^n$ denote the \emph{word length truncation} map defined by
\be{trunword} 
y = (x_1,\dots,x_n) \mapsto [y]_\tr := (x_1,\dots,x_{n \wedge \tr}),
\qquad  n\in\N,\,x_1,\dots,x_n\in E.
\ee
Extend this to a map from $\widetilde E^\N$ to $[\widetilde E]_\tr^\N$ via
\be{trunsent} 
\big[(y^{(1)}, y^{(2)},\dots)\big]_\tr := \big( [y^{(1)}]_\tr, [y^{(2)}]_\tr,\dots\big)
\ee
and to a map from $\mathcal{P}^{\mathrm{inv}}(\widetilde{E}^\N)$ to $\mathcal{P}^{\mathrm{inv}}
([\widetilde{E}]_\tr^\N)$ via
\be{trunprob}
[Q]_\tr(A) := Q(\{z\in\widetilde{E}^\N\colon\,[z]_\tr\in A\}), 
\qquad A \subset [\widetilde{E}]_\tr^\N \mbox{ measurable}.
\ee 
Note that if $Q \in \mathcal{P}^{\mathrm{inv}}(\widetilde{E}^\N)$, then $[Q]_\tr$ is an 
element of the set
\be{Pmean}
\mathcal{P}^{\mathrm{inv,fin}}(\widetilde{E}^\N)
= \{Q\in\mathcal{P}^{\mathrm{inv}}(\widetilde{E}^\N)\colon\,m_Q<\infty\}. 
\ee

\begin{theorem} 
\label{mainthm}
{\bf [Quenched LDP]}
Assume {\rm (\ref{rhocond})}. Then, for $\nu^{\otimes\N}$--a.s.\ all $X$, the family of 
(regular) conditional probability distributions $\Pr(R_N \in \cdot \mid X)$, $N\in\N$, 
satisfies the LDP on $\mathcal{P}^{\mathrm{inv}}(\widetilde{E}^\N)$ with rate $N$ and 
with deterministic rate function $I^\mathrm{que}\colon\,\mathcal{P}^{\mathrm{inv}}
(\widetilde{E}^\N) \to [0,\infty]$ given by
\be{eqgndefinitionI}
I^\mathrm{que}(Q) := \left\{\begin{array}{ll}
I^{\mathrm{fin}}(Q),
&\mbox{if } Q \in \mathcal{P}^{\mathrm{inv,fin}}(\widetilde{E}^\N),\\[1ex]
\lim\limits_{\tr \to \infty} I^{\mathrm{fin}}\big([Q]_\tr\big),
&\mbox{otherwise},
\end{array}
\right.
\ee 
where
\be{eqnratefctexplicit}
I^{\mathrm{fin}}(Q) := H(Q \mid q_{\rho,\nu}^{\otimes\N}) 
+ (\alpha-1) \, m_{Q}\,H(\Psi_{Q} \mid \nu^{\otimes \N}). 
\ee 
\end{theorem}

\begin{theorem} 
\label{ratefct}
The rate function $I^\mathrm{que}$ is lower semi-continuous, has compact level sets, has 
a unique zero at $Q=q_{\rho,\nu}^{\otimes\N}$, and is affine. Moreover, it is equal to 
the lower semi-continuous extension of $I^{\mathrm{fin}}$ from $\mathcal{P}^{\mathrm{inv,fin}}
(\widetilde{E}^\N)$ to $\mathcal{P}^{\mathrm{inv}}(\widetilde{E}^\N)$.
\end{theorem}

\noindent
Theorem~\ref{mainthm} will be proved in Sections~\ref{S3}--\ref{S5}, Theorem~\ref{ratefct} 
in Section~\ref{S6}.

A remarkable aspect of (\ref{eqnratefctexplicit}) in relation to (\ref{Iann}) is that 
it \emph{quantifies} the difference between the quenched and the annealed rate function.
Note the appearance of the tail exponent $\alpha$. 
We have not been able to 
find a simple formula for $I^\mathrm{que}(Q)$ when $m_Q=\infty$. In Appendix~\ref{SA} 
we will show that 
the annealed and the quenched rate function are continuous under truncation 
of word lengths, i.e., 
\be{conlimtr}
I^\mathrm{ann}(Q)=\lim_{\tr\to\infty} I^\mathrm{ann}([Q]_\tr),
\qquad
I^\mathrm{que}(Q)=\lim_{\tr\to\infty} I^\mathrm{que}([Q]_\tr),
\qquad
Q \in \mathcal{P}^{\mathrm{inv}}(\widetilde{E}^\N).
\ee 

Theorem~\ref{mainthm} is an extension of Birkner~\cite{Bi08}, Theorem~1. In that
paper, the quenched LDP is derived under the assumption that the law $\rho$ 
satisfies the \emph{exponential tail property}
\be{exptail}
\exists\,C<\infty,\,\lambda>0\colon\, 
\quad \rho(n) \leq Ce^{-\lambda n} \quad \forall\,n\in\N
\ee
(which includes the case where ${\rm supp}(\rho)$ is finite). The rate function 
governing the LDP is given by
\be{exptailrate}
I^\mathrm{que}(Q) := \left\{\begin{array}{ll}
H(Q \mid q_{\rho,\nu}^{\otimes\N}), &\mbox{if } Q\in\mathscr{R}_\nu,\\
\infty, &\mbox{if } Q\notin\mathscr{R}_\nu,
\end{array}
\right.
\ee
where 
\be{Rdef}
\mathscr{R}_\nu := \bigg\{Q \in \mathcal{P}^{\mathrm{inv}}(\widetilde{E}^\N)\colon\,
{\rm w}\!-\!\!\!\lim_{L\to\infty} \frac{1}{L} \sum_{k=0}^{L-1} \delta_{\theta^k\kappa(Y)}
= \nu^{\otimes\N}\,\,Q-a.s.\bigg\}. 
\ee
Think of $\mathscr{R}_\nu$ as the set of those $Q$'s for which \emph{the concatenation of words 
has the same statistical properties as the letter sequence $X$}. This set is not closed in the 
weak topology: its closure is $\mathcal{P}^{\mathrm{inv}}(\widetilde{E}^\N)$. 

We can include the cases where $\rho$ satisfies (\ref{rhocond}) with $\alpha=1$ or 
$\alpha=\infty$.

\begin{theorem} 
\label{mainthmboundarycases} 
{\rm (a)} If $\alpha=1$, then the quenched LDP holds with $I^\mathrm{que}=I^\mathrm{ann}$
given by {\rm (\ref{Iann})}.\\
{\rm (b)} If $\alpha=\infty$, then the quenched LDP holds with 
rate function 
\be{eq:ratefctalphainfty} 
I^\mathrm{que}(Q) = 
\begin{cases} 
H(Q \mid q_{\rho,\nu}^{\otimes\N}) & \mbox{if} \;\; 
\lim\limits_{\tr\to\infty} m_{[Q]_\tr} H( \Psi_{[Q]_\tr} \mid \nu^{\otimes \N}) = 0, \\
\infty & \mbox{otherwise}. 
\end{cases}
\ee
\end{theorem}

\noindent
Theorem~\ref{mainthmboundarycases} will be proved in Section~\ref{S7}. Part (a) says 
that the quenched and the annealed rate function are identical when $\alpha=1$. Part 
(b) says that (\ref{exptailrate}) can be viewed as the limiting case of 
(\ref{eqnratefctexplicit}) as $\alpha\to\infty$. Indeed, it was shown in 
Birkner~\cite{Bi08}, Lemma~2, that on $\mathcal{P}^{\mathrm{inv,fin}}(\widetilde{E}^\N)$:
\be{psiQprop}
\Psi_Q=\nu^{\otimes\N} \mbox{ if and only if } Q\in\mathscr{R}_\nu.
\ee 
Hence, (\ref{eq:ratefctalphainfty}) and (\ref{exptailrate}) agree on 
$\mathcal{P}^{\mathrm{inv,fin}}(\widetilde{E}^\N)$, and the rate function  
(\ref{eq:ratefctalphainfty}) is the lower semicontinuous extension of (\ref{exptailrate}) 
to $\mathcal{P}^{\mathrm{inv}}(\widetilde{E}^\N)$. 
By Birkner~\cite{Bi08}, Lemma~7, 
the expressions in (\ref{eq:ratefctalphainfty}) and (\ref{exptailrate}) are identical 
if $\rho$ has exponentially decaying tails. 
In this sense, Part (b) generalises the result in Birkner~\cite{Bi08}, Theorem 1, to arbitrary 
$\rho$ with a tail that decays faster than algebraic. 
\smallskip

Let $\pi_1\colon\,\widetilde{E}^\N\to\widetilde{E}$ be the projection onto the first word, 
and let $\mathcal{P}(\widetilde{E})$ be the set of probability measures on $\widetilde{E}$. 
An application of the contraction principle to Theorem~\ref{mainthm} yields the following.

\bc{cor:marginal}
Under the assumptions of Theorem~{\rm \ref{mainthm}}, for $\nu^{\otimes\N}$--a.s.\ all 
$X$, the family of (regular) conditional probability distributions $\Pr(\pi_1 R_N \in\cdot
\mid X)$, $N\in\N$, satisfies the LDP on $\mathcal{P}(\widetilde{E})$ with rate $N$ and with 
deterministic rate function $I^\mathrm{que}_1\colon\,\mathcal{P}(\widetilde{E}) \to [0,\infty]$ 
given by
\be{eq:contractedratefct} 
I^\mathrm{que}_1(q) 
:= \inf\big\{ I^\mathrm{que}(Q)\colon\, 
Q \in \mathcal{P}^{\mathrm{inv}}(\widetilde{E}^\N),\,\pi_1 Q = q \big\}. 
\ee
This rate function is lower semi-continuous, has compact levels sets, has a unique zero
at $q=q_{\rho,\nu}$, and is convex.
\ec
\noindent 
Corollary~\ref{cor:marginal} shows that the rate function in Birkner~\cite{Bi03}, 
Theorem~6, must be replaced by (\ref{eq:contractedratefct}).
It does not appear possible to evaluate the 
infimum in (\ref{eq:contractedratefct}) explicitly in general. For a 
$q \in \mathcal{P}(\widetilde{E})$ with finite mean length and 
$\Psi_{q^{\otimes\N}}=\nu^{\otimes\N}$, we have $I^\mathrm{que}_1(q)=h(q \mid q_{\rho,\nu})$. 
\smallskip 

By taking projective limits, it is possible to extend Theorems~\ref{mainthm}--\ref{ratefct} 
to more general letter spaces. 
See, e.g., Deuschel and Stroock\ \cite{DeStr89},  Section~4.4, or 
Dembo and Zeitouni \cite{DeZe98}, Section~6.5, for background on 
(specific) relative entropy in general spaces. 
The following corollary will be proved in Section~\ref{Scount}.

\bc{LDPcount}
The quenched LDP also holds when $E$ is a Polish space, 
with the same rate function as in
{\rm (\ref{eqgndefinitionI}--\ref{eqnratefctexplicit})}.
\ec 

\medskip 

In the companion paper \cite{BGdH08} the annealed and quenched LDP are applied to 
the collision local time of transient random walks, and the existence 
of an intermediate phase for a class of interacting stochastic systems 
is established. 

\subsection{Heuristic explanation of main theorems}
\label{S1.3}

To explain the background of Theorem~\ref{mainthm}, we begin by recalling a few properties of 
entropy. Let $H(Q)$ denote the \emph{specific entropy of $Q \in \mathcal{P}^{\mathrm{inv}}
(\widetilde{E}^\N)$} defined by
\be{entrdef}
H(Q) := \lim_{N\to\infty} \frac{1}{N}\,h\big(Q_{|_{\mathscr{F}_N}}\big)
\in [0,\infty],
\ee
where $h(\cdot)$ denotes entropy. The sequence under the limit in (\ref{entrdef}) 
is non-increasing in $N$. Since $q_{\rho,\nu}^{\otimes\N}$ is a product measure, we 
have the identity (recall (\ref{Tdefs}--\ref{q0def}))
\be{HPsirelQ}
\begin{aligned}
H(Q \mid q_{\rho,\nu}^{\otimes\N}) &= -H(Q) - \E_Q[\log q_{\rho,\nu}(Y_1)]\\
&= -H(Q) - \E_Q[\log\rho(\tau_1)] - m_Q\,\E_{\Psi_Q}[\log\nu(X_1)]. 
\end{aligned}
\ee
Similarly,
\be{HPsirel}
H(\Psi_Q \mid\nu^{\otimes\N}) = -H(\Psi_Q) - \E_{\Psi_Q}[\log \nu(X_1)]. 
\ee

Below, for a discrete random variable $Z$ with a law $Q$ on a state space $\mathcal{Z}$ we 
will write $Q(Z)$ for the random variable $f(Z)$ with $f(z)=Q(Z=z)$, $z\in\mathcal{Z}$. 
Abbreviate 
\be{eq:defKN}
K^{(N)} := \kappa(Y^{(1)},\dots,Y^{(N)}) 
\quad \mbox{and} \quad K^{(\infty)} := \kappa(Y).
\ee 
In analogy with 
(\ref{Pmean}), define
\be{Pmeanerg}
\mathcal{P}^{\mathrm{erg,fin}}(\widetilde{E}^\N) := 
\Big\{ Q \in \mathcal{P}^{\mathrm{erg}}(\widetilde{E}^\N)\colon\, m_Q < \infty \Big\}.
\ee

\begin{lemma} 
\label{Birkner}
{\rm [Birkner~\cite{Bi08}, Lemmas~3 and~4]}\\ \hfill
Suppose that $Q\in\mathcal{P}^{\mathrm{erg,fin}}(\widetilde{E}^\N)$
and $H(Q)<\infty$. Then, $Q$-a.s.,
\be{eqnconcatse}
\begin{aligned}
\lim_{N\to\infty} \frac{1}{N} 
\log Q(K^{(N)}) &= - m_Q H(\Psi_Q),\\
\lim_{N\to\infty} \frac{1}{N} 
\log Q\big(\tau_1,\dots,\tau_N \mid K^{(N)} \big) &=: - H_{\tau|K}(Q),\\
\lim_{N\to\infty} \frac{1}{N} 
\log Q\big(Y^{(1)},\dots,Y^{(N)}\big) &= -H(Q), 
\end{aligned}
\ee
with
\be{Hrels}
m_Q H(\Psi_Q) + H_{\tau|K}(Q) = H(Q).
\ee
\end{lemma} 

\medskip\noindent
Equation (\ref{Hrels}), which follows from (\ref{eqnconcatse}) and the identity
\be{Qrels}
Q(K^{(N)})Q(\tau_1,\dots,\tau_N \mid K^{(N)}) 
= Q(Y^{(1)},\dots,Y^{(N)}),
\ee 
identifies $H_{\tau|K}(Q)$. Think of $H_{\tau|K}(Q)$ as the \emph{conditional specific 
entropy of word lengths under the law $Q$ given the concatenation}. Combining 
(\ref{HPsirelQ}--\ref{HPsirel}) and (\ref{Hrels}), we have
\be{eqnsre1}
H(Q \mid q_{\rho,\nu}^{\otimes\N}) = m_Q H(\Psi_Q \mid \nu^{\otimes\N}) 
- H_{\tau|K}(Q) - \E_Q[\log\rho(\tau_1)].
\ee 
The term $- H_{\tau|K}(Q)-\E_Q[\log\rho(\tau_1)]$ in (\ref{eqnsre1}) can be interpreted as the 
\emph{conditional specific relative entropy of word lengths under the law $Q$ w.r.t.\ 
$\rho^{\otimes\N}$ given the concatenation}. 

Note that $m_Q<\infty$ and $H(Q)<\infty$ imply that $H(\Psi_Q)<\infty$, as can be seen
from (\ref{Hrels}). Also note that $-\E_{\Psi_Q}[\log \nu(X_1)]<\infty$ because $E$ is
finite, and $-\E_Q[\log\rho(\tau_1)]<\infty$ because of (\ref{rhocond}) and 
$m_Q<\infty$, implying that (\ref{HPsirelQ}--\ref{HPsirel}) are proper.

We are now ready to give a heuristic explanation of Theorem~\ref{mainthm}. Let 
\be{RNdefind}
R^N_{j_1,\dots,j_N}(X), \qquad 0<j_1<\dots<j_N<\infty,
\ee 
denote the empirical process of $N$-tuples of words when $X$ is cut at the points $j_1,
\dots,j_N$ (i.e., when $T_i=j_i$ for $i=1,\dots,N$; see (\ref{xidef}--\ref{RNaltdef}) 
for a precise definition). Fix $Q \in \mathcal{P}^{\mathrm{erg,fin}}(\widetilde{E}^\N)$. 
The probability $\Pr(R_N \approx Q \mid X)$ is a sum over all $N$-tuples $j_1,\dots,j_N$ 
such that $R^N_{j_1,\dots,j_N}(X) \approx Q$, weighted by $\prod_{i=1}^N \rho(j_i-j_{i-1})$ 
(with $j_0=0$). The fact that $R^N_{j_1,\dots,j_N}(X) \approx Q$ has three consequences: 

\begin{enumerate}
\item[(1)] 
The $j_1,\dots,j_N$ must cut $\approx N$ substrings out of $X$ of total length $\approx N m_Q$ 
that look like the concatenation of words that are $Q$-typical, i.e., that look as if 
generated by $\Psi_Q$ (possibly with gaps in between). This means that most of the cut-points 
must hit atypical pieces of $X$. We expect to have to shift $X$ by $\approx\exp
[N m_Q H(\Psi_Q \mid \nu^{\otimes\N})]$ in order to find the first contiguous 
substring of length $N m_Q$ whose empirical shifts lie in a small neighbourhood 
of $\Psi_Q$. By (\ref{rhocond}), the probability for the single increment $j_1-j_0$ 
to have the size of this shift is $\approx \exp[-N\alpha\,m_Q H(\Psi_Q \mid \nu^{\otimes\N})]$. 
\item[(2)] 
The combinatorial factor $\exp[NH_{\tau|K}(Q)]$ counts how many ``local perturbations'' 
of $j_1,\dots,j_N$ preserve the property that $R^N_{j_1,\dots,j_N}(X)\approx Q$. 
\item[(3)] 
The statistics of the increments $j_1-j_0,\dots,j_N-j_{N-1}$ must be close to the 
distribution of word lengths under $Q$. Hence, the weight factor $\prod_{i=1}^N 
\rho(j_i-j_{i-1})$ must be $\approx \exp[N \E_Q[\log\rho(\tau_1)]]$ (at least, for 
$Q$-typical pieces).
\end{enumerate}

\noindent
The contributions from (1)--(3), together with the identity in (\ref{eqnsre1}), explain the 
formula in (\ref{eqnratefctexplicit}) on $\mathcal{P}^{\mathrm{erg,fin}}(\widetilde{E}^\N)$. 
Considerable work is needed to extend (1)--(3) from $\mathcal{P}^{\mathrm{erg,fin}}(\widetilde{E}^\N)$ 
to $\mathcal{P}^{\mathrm{inv}}(\widetilde{E}^\N)$. This is explained in Section~\ref{S3.5}.

In (1), instead of having a single large increment preceding a single contiguous 
substring of length $Nm_Q$, it is possible to have several large increments preceding 
several contiguous substrings, which together have length $Nm_Q$. The latter gives rise 
to the same contribution, and so there is some entropy associated with the choice of
the large increments. Lemma~\ref{lowtoylemma} in Section \ref{S2.1} is needed to control 
this entropy, and shows that it is negligible.

\subsection{Outline}
\label{S1.6}

Section~\ref{S2} collects some preparatory facts that are needed for the proofs of 
the main theorems, including a lemma that controls the entropy associated with
the locations of the large increments in the renewal process. In Section~\ref{S3} 
and \ref{S4} we prove the large deviation upper, respectively, lower bound. The 
proof of the former is long (taking up about half of the paper) and requires a 
somewhat lengthy construction with combinatorial, functional analytic and ergodic 
theoretic ingredients. In particular, extending the lower bound from ergodic 
to non-ergodic probability measures is technically involved. The proofs of 
Theorems~\ref{mainthm}--\ref{mainthmboundarycases} are in 
Sections~\ref{S5}--\ref{S7}, that of Corollary~\ref{LDPcount} is 
in Section~\ref{Scount}.
Appendix~\ref{SA} contains a proof that the annealed and the 
quenched rate function are continuous under the truncation of the word length 
approximation.


\section{Preparatory facts} 
\label{S2}

Section~\ref{S2.1} proves a core lemma that is needed to control the entropy
of large increments in the renewal process. Section~\ref{S2.2} shows that the 
tail property of $\rho$ is preserved under convolutions. 

\subsection{A core lemma}
\label{S2.1}

As announced at the end of Section~\ref{S1.3}, we need to account for the entropy 
that is associated with the locations of the large increments in the renewal process. 
This requires the following combinatorial lemma.

\begin{lemma}
\label{lowtoylemma} 
Let $\omega=(\omega_l)_{l\in\N}$ be i.i.d.\ with $\Pr(\omega_1=1)=1-\Pr(\omega_1=0)=p
\in (0,1)$, and let $\alpha \in (1,\infty)$. For $N\in\N$, let
\be{eqn1b} 
S_N(\omega) 
:= \sum_{ {0<j_1<\dots<j_N<\infty} \atop {\omega_{j_1} = \dots = \omega_{j_N}=1} }
\,\,\prod_{i=1}^N (j_i-j_{i-1})^{-\alpha} \qquad (j_0=0)
\ee
and put 
\be{LDlimsup} 
\limsup_{N\to\infty} \frac{1}{N} \log S_N(\omega) =: - \phi(\alpha,p) \qquad \omega-a.s. 
\ee
(the limit being $\omega$-a.s.\ constant by tail triviality). Then 
\be{phitildebound}
\lim_{p \da 0} \frac{\phi(\alpha,p)}{\alpha \log(1/p)} = 1. 
\ee
\end{lemma} 

\begin{proof}
Let $\tau_N := \min\{l \in \N\colon\, \omega_{l}=\omega_{l+1}=\cdots=\omega_{l+N-1}\}$. 
In (\ref{eqn1b}), choosing $j_1=\tau_N$ and $j_i=j_{i-1}+1$ for $i=2,\dots,N$, we 
see that $S_N(\omega) \geq \tau_N^{-\alpha}$. Since 
\be{tauLLN}
\lim_{N \to \infty} \frac{1}{N} \log \tau_N \to \log(1/p) \qquad \omega-a.s., 
\ee
we have
\be{phiub}
\phi(\alpha,p) \le \alpha \log(1/p) \qquad \forall\,p \in (0,1). 
\ee
To show that this bound is sharp in the limit as $p \da 0$, we estimate fractional moments 
of $S_N(\omega)$.

For any $\beta \in (1/\alpha,1]$, using that
$(u+v)^\beta \leq u^\beta+v^\beta$, $u,v \geq 0$, we get 
\be{philb1}
\begin{aligned}
\E\Big[ S_N(\omega)^\beta\Big]
&\leq \sum_{0<j_1<\dots<j_N<\infty} \E\Big[\1_{\{\omega_{j_1}=\dots=\omega_{j_N}=1\}} 
\prod_{i=1}^N (j_i-j_{i-1})^{-\alpha\beta}\Big]\\ 
&= \sum_{0<j_1<\dots<j_N<\infty} p^N \prod_{i=1}^N (j_i-j_{i-1})^{-\alpha\beta}\\ 
&= \big[ p \, \zeta(\alpha\beta) \big]^N,
\end{aligned}
\ee
where $\zeta(s) = \sum_{n\in\N} n^{-s}$, $s>1$, is Riemann's $\zeta$-function. 
Hence, for any $\varepsilon > 0$, Markov's inequality yields 
\be{philb2b}
\begin{aligned}
\Pr\Big( & \frac{1}{N} \log S_N(\omega) \geq \frac{1}{\beta} \big[ \log
  p + \log \zeta(\alpha\beta) + \varepsilon \big] \Big) \\ 
& = \Pr\Big(
S_N(\omega)^\beta \geq e^{\varepsilon N} \big[ p \, \zeta(\alpha\beta)
  \big]^N\Big) \leq e^{-\varepsilon N} 
  \big[ p \, \zeta(\alpha\beta) \big]^{-N} \E\Big[ S_N(\omega)^\beta\Big]
  \leq e^{-\varepsilon N}. 
\end{aligned}
\ee
Thus, by the first Borel-Cantelli Lemma, 
\be{} 
-\phi(\alpha,p) = 
\limsup_{N\to\infty} \frac1N \log S_N(\omega) \leq 
\frac{1}{\beta} \big[ \log
  p + \log \zeta(\alpha\beta) \big] \quad \mbox{a.s.}
\ee
Now let $p \da 0$, followed by $\beta \da 1/\alpha$ to 
obtain the claim.
\end{proof}

\begin{remark} \rm
Note that $\E[S_N(\omega)]=(p \zeta(\alpha))^N$, while typically 
$S_N(\omega) \approx p^{\alpha N}$. In the above computation, this is verified 
by bounding suitable non-integer moments of $S_N(\omega)/p^{\alpha N}$. 
Estimating non-integer moments in situations when the 
mean is inconclusive is a useful technique in a variety of different probabilistic 
contexts. See, e.g., Holley and Liggett\ \cite{HoLi81} and Toninelli\ \cite{To07}.  
The proof of Lemma~\ref{lowtoylemma} above is similar to that of 
Toninelli\ \cite{To07}, Theorem~2.1. 
\end{remark}

\subsection{Convolution preserves polynomial tail}
\label{S2.2}

The following lemma will be needed in Sections~\ref{S3.3} and \ref{S3.5}. 
For $m\in\N$, let $\rho^{*m}$ denote the $m$-fold convolution of $\rho$.

\begin{lemma} 
\label{lemmarhoconv} 
Suppose that $\rho$ satisfies $\rho(n) \leq C_\rho \, n^{-\alpha}$, $n \in \N$,
for some $C_\rho<\infty$. Then
\be{eqnlemmarhoconv}
\rho^{*m}(n) \le (C_\rho \vee 1) \, m^{\alpha+1} n^{-\alpha} \qquad 
\forall\,m,n\in\N. 
\ee
\end{lemma}

\begin{proof} 
If $n \leq m$, then the right-hand side of (\ref{eqnlemmarhoconv}) is $\geq 1$. So, 
let us assume that $n>m$. Then 
\be{count9}
\begin{aligned}
\rho^{*m}(n)  
&=  \hspace{-1em} 
\sum_{x_1,\dots,x_m \ge 1 \atop x_1+\cdots+x_m=n} \prod_{i=1}^m \rho(x_i) 
\leq \sum_{j=1}^m 
\sum_{\begin{array}{c} 
\scriptstyle x_1,\dots,x_m \geq 1 \\[-1ex] 
\scriptstyle x_1+\cdots+x_m=n \\[-1ex] 
\scriptstyle x_j = x_1 \vee\dots\vee x_m 
\end{array}} 
\rho(x_j)\, \prod_{i \neq j}^m \rho(x_i)\\ 
&\leq m\,C_\rho \, \lceil n/m \rceil^{-\alpha} 
\sum_{x_1,\dots,x_{m-1} \geq 1} \prod_{i=1}^{m-1} \rho(x_i)\\
&= m\,C_\rho \, \lceil n/m \rceil^{-\alpha} 
\leq C_\rho \,m^{\alpha+1}\,n^{-\alpha}.
\end{aligned}
\ee
\end{proof}


\section{Upper bound} 
\label{S3}

The following upper bound will be used in Section~\ref{S5} to derive
the upper bound in the definition of the LDP.

\begin{proposition}
\label{propub} 
For any $Q \in \mathcal{P}^{\mathrm{inv,fin}}(\widetilde{E}^\N)$ and any 
$\varepsilon>0$, there is an open neighbourhood $\mathcal{O}(Q)\subset
\mathcal{P}^{\mathrm{inv}}(\widetilde{E}^\N)$ of $Q$ such that  
\be{mainub} 
\limsup_{N\to\infty} \frac{1}{N} 
\log \Pr\big(R_N \in \mathcal{O}(Q) \mid X\big) 
\leq -I^{\mathrm{fin}}(Q) + \varepsilon \quad X-a.s.
\ee
\end{proposition}
\noindent
We remark that since $|E|<\infty$ we automatically have 
$I^{\mathrm{fin}}(Q) \in [0,\infty)$ for all 
$Q \in \mathcal{P}^{\mathrm{inv,fin}}(\widetilde{E}^\N)$, 
so the right-hand side of (\ref{mainub}) is finite.

\begin{proof} 
It suffices to consider the case $\Psi_Q \neq \nu^{\otimes \N}$. The case $\Psi_Q 
= \nu^{\otimes\N}$, for which $I^{\mathrm{fin}}(Q) = H(Q \mid q_{\rho,\nu}^{\otimes\N})$ 
as is seen from (\ref{eqnratefctexplicit}), is contained in the upper bound in Birkner 
\cite{Bi08}, Lemma 8. Alternatively, by lower semicontinuity of $Q' \mapsto H(Q'\mid 
q_{\rho,\nu}^{\otimes\N})$, there is a neighbourhood $\mathcal{O}(Q)$ such that 
\be{lscneight}
\inf_{Q' \in \overline{\mathcal{O}(Q)}} H(Q'\mid q_{\rho,\nu}^{\otimes\N}) 
\geq H(Q\mid q_{\rho,\nu}^{\otimes\N}) - \varepsilon = I^{\mathrm{fin}}(Q) 
- \varepsilon,
\ee 
where $\overline{\mathcal{O}(Q)}$ denotes the closure of $\mathcal{O}(Q)$ (in the weak
topology), and we can use the annealed bound.

In Sections~\ref{S3.1}--\ref{S3.5} we first prove Proposition~\ref{mainub} 
under the assumption that there exist $\alpha \in (1,\infty),\,C_\rho < \infty$ such that 
\be{eq:algtail}
\rho(n) \leq C_\rho \, n^{-\alpha}, 
\quad n \in \N,
\ee
which is needed in Lemma~\ref{lemmarhoconv}. In Section~\ref{S3.6} we show 
that this can be replaced by (\ref{rhocond}). In Sections~\ref{S3.1}--\ref{S3.4},
we first consider $Q \in \mathcal{P}^{\mathrm{erg,fin}}(\widetilde{E}^\N)$ (recall
(\ref{Pmeanerg})). Here, we turn the heuristics from Section~\ref{S1.3} into a rigorous 
proof. In Section~\ref{S3.5} we remove the ergodicity restriction. The proof is long
and technical (taking up more than half of the paper).

\subsection{Step 1: Consequences of ergodicity}
\label{S3.1}

We will use the ergodic theorem to construct specific neighborhoods of $Q\in
\mathcal{P}^{\mathrm{erg,fin}}(\widetilde{E}^\N)$ that are well adapted to formalize 
the strategy of proof outlined in our heuristic explanation of the main theorem
in Section~\ref{S1.3}. 

Fix $\varepsilon_1,\delta_1>0$. By the ergodicity of $Q$ and Lemma~\ref{Birkner}, 
the event (recall (\ref{eq:defmQ}) and (\ref{eq:defKN}))
\be{intevents}
\begin{aligned} 
&\left\{\frac{1}{M}|K^{(M)}| \in m_Q + [-\varepsilon_1, \varepsilon_1] \right\}\\[.2cm] 
&\cap 
\left\{- \frac{1}{M} \log Q(K^{(M)}) 
\in m_Q H(\Psi_Q) + [-\varepsilon_1, \varepsilon_1] \right\}\\[.2cm] 
&\cap 
\left\{- \frac{1}{M} \log Q(Y^{(1)},\dots,Y^{(M)}) 
\in H(Q) + [-\varepsilon_1, \varepsilon_1] \right\} \\[.2cm]
&\cap 
\left\{\frac{1}{M} \sum_{k=1}^{|K^{(M)}|} \log \nu((K^{(M)})_k) 
\in m_Q \E_{\Psi_Q}\big[\log \nu(X_1)\big] + [-\varepsilon_1, \varepsilon_1] \right\}\\[.2cm] 
&\cap 
\left\{\frac{1}{M} \sum_{i=1}^{M} \log \rho(\tau_i) 
\in \E_Q\big[\log \rho(\tau_1) \big] + [-\varepsilon_1, \varepsilon_1] \right\} 
\end{aligned}
\ee
has $Q$-probability at least $1-\delta_1/4$ for $M$ large enough (depending on $Q$), 
where $|K^{(M)}|$ is the length of the string of letters $K^{(M)}$. Hence, there is 
a finite number $A$ of sentences of length $M$, denoted by 
\be{zidefs}
(z_a)_{a=1,\dots,A} \mbox{ with } z_a:=(y^{(a,1)},\dots,y^{(a,M)}) \in \widetilde{E}^{M},
\ee 
such that for $a=1,\dots,A$,
\be{sentprop}
\begin{aligned} 
|\kappa(z_a)| 
&\in \Big[M(m_Q-\varepsilon_1),M(m_Q+\varepsilon_1)\Big],\\[.2cm]
Q(K^{(M)}=\kappa(z_a)) 
&\in \Big[\exp[-M(m_Q H(\Psi_Q) + \varepsilon_1)], 
\exp[-M(m_Q H(\Psi_Q) - \varepsilon_1)]\Big],\\[.2cm]
Q\big((Y^{(1)},\dots,Y^{(M)})=z_a\big) 
&\in \Big[\exp[-M(H(Q) + \varepsilon_1)],
\exp[-M(H(Q) - \varepsilon_1)]\Big],\\[.2cm]
\sum_{k=1}^{|\kappa(z_a)|} \log \nu((\kappa(z_a))_k) 
&\in \Big[M(m_Q \E_{\Psi_Q}[\log \nu(X_1)] - \varepsilon_1), 
M(m_Q \E_{\Psi_Q}[\log \nu(X_1)] + \varepsilon_1)\Big],\\[.2cm]
\sum_{i=1}^{M} \log \rho(|y^{(a,i)}|) 
&\in \Big[M(\E_Q[\log \rho(\tau_1)] - \varepsilon_1), 
M(\E_Q[\log \rho(\tau_1)] + \varepsilon_1)\Big], 
\end{aligned}
\ee
and
\be{probsent}
\sum_{a=1}^A Q\Big((Y^{(1)},\dots,Y^{(M)})=z_a\Big) 
\geq 1 - \frac{\delta_1}{2}. 
\ee
Note that (\ref{probsent}) and the third line of (\ref{sentprop}) imply that 
\be{Abdsint}
A \in \Big[\big(1-\frac{\delta_1}{2}\big) \exp\big[M(H(Q)-\varepsilon_1)\big],
\exp\big[M(H(Q)+\varepsilon_1)\big]\Big].
\ee 

Abbreviate 
\be{mathscrAdef}
\mathscr{A} := \{ z_a,\, a=1,\dots,A\}.
\ee 
Let 
\be{Bdef}
\mathscr{B} := \big\{ \zeta^{(b)}, b=1,\dots,B \big\} 
= \big\{ \kappa(z_a),\, a=1,\dots,A \big\} 
\ee
be the set of strings of letters arising from concatenations of the individual $z_a$'s, 
and let
\be{Ijdef}
I_b := \big\{ 1 \leq a \le A\colon\, \kappa(z_a) = \zeta^{(b)} \big\}, 
\quad b = 1, \dots, B,
\ee 
so that $|I_b|$ is the number of sentences in $\mathscr{A}$ giving a particular string 
in $\mathscr{B}$. By the second line of (\ref{sentprop}), we can bound $B$ as
\be{eqnsizeofB}
B \leq \exp\big[ M (m_Q H(\Psi_Q) + \varepsilon_1) \big], 
\ee
because $\sum_{b=1}^B Q(K^{(M)} = \zeta^{(b)}) \leq 1$ and each summand is at 
least $\exp[ -M (m_Q H(\Psi_Q) + \varepsilon_1)]$. Furthermore, we have 
\be{eqnsizeIj}
|I_b| \leq \exp\big[ M (H_{\tau|K}(Q) + 2 \varepsilon_1) \big] , \quad b=1,\dots,B,  
\ee
since
\be{sumQsand}
\begin{aligned}
\exp\big[-M(m_Q H(\Psi_Q) - \varepsilon_1)\big] 
&\geq 
Q\big(\kappa(Y^{(1)},\dots,Y^{(M)}) = \zeta^{(b)}\big) \\
&\geq 
\sum_{a \in I_b} Q\big((Y^{(1)},\dots,Y^{(M)}) =  z_a \big) 
\geq |I_b| \exp\big[-M(H(Q)+ \varepsilon_1)\big], 
\end{aligned}
\ee
and $H(Q) - m_Q H(\Psi_Q) = H_{\tau|K}(Q)$ by (\ref{Hrels}).%

\subsection{Step 2: Good sentences in open neighbourhoods}
\label{S3.2}

Define the following open neighbourhood of $Q$ (recall (\ref{mathscrAdef})) 
\be{Opnedef}
\mathcal{O}:= \Big\{Q' \in \mathcal{P}^{\mathrm{inv}}(\widetilde{E}^\N)\colon\, 
Q'_{|_{\mathscr{F}_M}}(\mathscr{A}) > 1 - \delta_1 \Big\}.
\ee
Here, $Q(z)$ is shorthand for $Q((Y^{(1)},\dots,Y^{(M)})=z)$. For $x \in E^\N$ 
and for a vector of cut-points $(j_1,\dots,j_N) \in \N^N$ with $0<j_1<\dots<j_N<\infty$ 
and $N > M$, let 
\be{xidef}
\xi_N := (\xi^{(i)})_{i=1,\dots,N} = 
\big( x|_{(0,j_1]}, x|_{(j_1,j_2]}, \dots, x|_{(j_{N-1},j_N]} \big) 
\in \widetilde{E}^N 
\ee
(with $(0,j_1]$ shorthand notation for $(0,j_1] \cap \N$, etc.) be the sequence of words 
obtained by cutting $x$ at the positions $j_i$, and let
\be{RNaltdef}
R^N_{j_1,\dots,j_N}(x) 
:= \frac{1}{N} \sum_{i=0}^{N-1} 
\delta_{ \textstyle \widetilde{\theta}^i (\xi_N)^\mathrm{per}} 
\ee
be the corresponding empirical process. By (\ref{Opnedef}), 
\be{eqnRNO}
\begin{aligned}
&R^N_{j_1,\dots,j_N}(x) \in \mathcal{O} \qquad \Longrightarrow\\[1ex] 
& \qquad \# \Big\{ 1\le i \leq N-M \colon\,  
\big(x|_{(j_{i-1},j_i]},\dots,x|_{(j_{i+M-1},j_{i+M}]}\big) 
\in \mathscr{A} \Big\} \geq N(1-\delta_1)-M.  
\end{aligned}
\ee 
Note that (\ref{eqnRNO}) implies that the sentence $\xi_N$ contains at least 
\be{defG}
C := \lfloor(1-\delta_1)N/M\rfloor -1
\ee
disjoint subsentences from the set $\mathscr{A}$, i.e., there are $1 \leq i_1,\dots,i_C 
\leq N-M$ with $i_c - i_{c-1} \geq M$ for $c=1,\dots,C$ such that 
\be{xiprop} 
\big( \xi^{(i_c)}, \xi^{(i_c+1)}, \dots, \xi^{(i_c+M-1)} \big) \in \mathscr{A} 
\ee 
(we implicitly assume that $N$ is large enough so that $C>1$). Indeed, we can e.g.\ 
construct the $i_c$'s iteratively as
\be{itis}
\begin{aligned}
i_0 &= -M,\\ 
i_c &= \min\Big\{k \geq i_{c-1}+M\colon\ 
\mbox{a sentence from } \mathscr{A} \mbox{ starts at position } k \mbox{ in } 
\xi_N \Big\},\\ 
&\qquad c=1,\dots,C,
\end{aligned}
\ee
and we can continue the iteration as long as $c M + \delta_1 N \leq N$.
But (\ref{xiprop}) in turn implies that the $j_{i_c}$'s cut out of $x$ at 
least $C$ disjoint subwords from $\mathscr{B}$, i.e.,
\be{eqnXsubwords} 
x|_{(j_{i_c},j_{i_c+M}]} \in \mathscr{B}, \quad c=1,\dots,C. 
\ee

\subsection{Step 3: Estimate of the large deviation probability}
\label{S3.3}

\begin{figure}[hctb]
\vspace{0.5cm}
\begin{center} 
\psfrag{filling_subsentences}{filling subsentences}
\psfrag{good_subsentences}{good subsentences}
\psfrag{atypical_medium}{medium $\approx \Psi_Q$}
\psfrag{X}{$X$}
\includegraphics{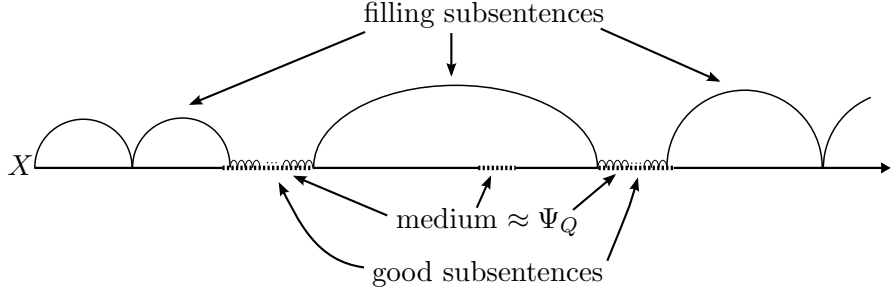}
\end{center}
\caption{Looking for good subsentences and filling subsentences
(see below (\ref{Xls})).}\label{fig:thequest}
\end{figure}

Using Steps 1 and 2, we estimate (recall (\ref{Opnedef})) 
\be{RNform} 
\Pr\big(R_N \in \mathcal{O} \mid X \big) 
= \sum_{0<j_1<\dots<j_N<\infty} 
\1_{\mathcal{O}}\left(R^N_{j_1,\dots,j_N}(X)\right) \,
\prod_{i=1}^N \rho(j_i-j_{i-1}) 
\ee
from above as follows. Fix a vector of cut-points $(j_1,\dots,j_N)$  giving rise 
to a non-zero contribution in the right-hand side of (\ref{RNform}). We think of 
this vector as describing a particular way of cutting $X$ into a sentence of $N$ 
words. By (\ref{eqnXsubwords}), at least $C$ (recall \ref{defG}) of the $j_c$'s 
must be cut-points where a word from $\mathscr{B}$ is written on $X$, and these $C$ 
subwords must be disjoint. As words in $\mathscr{B}$ arise from concatenations of 
sentences from $\mathscr{A}$, this means we can find 
\be{l1lG}
\ell_1 < \cdots < \ell_C, \quad \{\ell_1,\dots,\ell_C\} \subset \{0, j_1,\dots,j_N\} 
\quad \mbox{and} \quad \zeta_1,\dots,\zeta_C \in \mathscr{A}
\ee
such that 
\be{Xls}
X|_{(\ell_c, \ell_c+|\kappa(\zeta_c)|]} = \kappa(\zeta_c) =: \eta^{(c)} \in \mathscr{B}
\;\; \mbox{and} \;\; \ell_c \ge \ell_{c-1}+|\kappa(\zeta_{c-1})|, \quad c=1,\dots,C-1. 
\ee
We call $\zeta_1, \dots, \zeta_C$ the \emph{good} subsentences. 

Note that once we fix the $\ell_c$'s and the $\zeta_c$'s, this determines $C+1$ 
\emph{filling} subsentences (some of which may be empty) consisting of the words 
between the good subsentences. See Figure~\ref{fig:thequest} for an illustration. 
In particular, this determines numbers $m_1,\dots,m_{C+1}
\in \N$ such that $m_1+\cdots+m_{C+1}=N-CM$, where $m_c$ is the number of words we cut 
between the $(c-1)$-st and the $c$-th good subsentence (and $m_{C+1}$ is the number of 
words after the $C$-th good subsentence). 

Next, let us fix good $\ell_1 < \cdots < \ell_C$ and $\eta^{(1)}, \dots, \eta^{(C)} 
\in \mathscr{B}$, satisfying 
\be{eq:goodpoints}
X|_{(\ell_c,\ell_c+|\eta^{(c)}|]}=\eta^{(c)}, \quad 
\ell_{c} \geq \ell_{c-1}+|\eta^{(c-1)}|, \qquad c=1,\dots,C. 
\ee 
To estimate how many different choices of $(j_1,\dots,j_N)$ may lead to this 
particular $((\ell_c), (\eta^{(c)}))$, we proceed as follows. There are at most 
\be{count1} 
\big( 2M \varepsilon_1 \big)^C\,
\exp\big[M\big(H_{\tau|K}(Q) + 2 \varepsilon_1\big) \big]^C 
\leq \exp\big[N \big(H_{\tau|K}(Q) + \delta_2\big)\big]  
\ee
possible choices for the word lengths inside these good subsentences. Indeed, 
by the first line of (\ref{sentprop}), at most $2M \varepsilon_1$ different 
elements of $\mathscr{B}$ can start at any given position $\ell_c$ and, by 
(\ref{eqnsizeIj}), each of them can be cut in at most $\exp\big[M(H_{\tau|K}(Q) 
+ 2 \varepsilon_1)\big]$ different ways to obtain an element of $\mathscr{A}$. 
In (\ref{count1}), $\delta_2 = \delta_2(\varepsilon_1,\delta_1,M)$ can be made 
arbitrarily small by choosing $M$ large and $\varepsilon_1, \delta_1$ small. 
Furthermore, there are at most 
\be{count2} 
{N-C(M-1) \choose C} \leq \exp[\delta_3 N] 
\ee
possible choices of the $m_c$'s, where $\delta_3=\delta_3(\delta_1,M)$ can be made 
arbitrarily small by choosing $M$ large and $\delta_1$ small.

Next, we estimate the value of $\prod_{i=1}^N \rho(j_i-j_{i-1})$ for any 
$(j_1,\dots,j_N)$ leading to the given $((\ell_c), (\eta^{(c)}))$. In view of the 
fifth line of (\ref{sentprop}), we have 
\be{count1ext}
\begin{aligned}
&\prod_{i=1}^N \,\rho(j_i-j_{i-1})^{\textstyle \1_{\{\mbox{\footnotesize the $i$-th word falls inside 
the $C$ good subsentences}\}}}\\  
&\qquad \leq \exp\big[ C M \big(\E_Q[\log\rho(\tau_1)] + \varepsilon_1\big)\big]\\[0.2cm] 
&\qquad \leq \exp\big[ N \big(\E_Q[\log\rho(\tau_1)] + \delta_4\big)\big],  
\end{aligned}
\ee
where $\delta_4 = \delta_4(\varepsilon_1,\delta_1,M)$ can be made arbitrarily small 
by choosing $M$ large and $\varepsilon_1,\delta_1$ small. The filling subsentences 
have to exactly fill up the gaps between the good subsentences and so, for a given 
choice of $(\ell_c)$, $(\eta^{(c)})$ and $(m_c)$, the contribution to 
$\prod_{i=1}^N \rho(j_i-j_{i-1})$ from the filling subsentences is $\prod_{c=1}^C 
\rho^{* m_c}(\ell_c - \ell_{c-1} - |\eta^{(c-1)}|)$ (the term for $c=1$ is to be 
interpreted as $\rho^{* m_1}(\ell_1)$, and $\rho^{* 0}$ as $\delta_0$). By 
Lemma~\ref{lemmarhoconv}, using (\ref{eq:algtail}), 
\be{count4}
\begin{aligned}
&\prod_{c=1}^C \rho^{* m_c}\big(\ell_c - \ell_{c-1} - |\eta^{(c-1)}|\big)\\  
&\qquad \leq (C_\rho \vee 1)^C \left(\prod_{c=1}^C m_c^{\alpha+1}\right) 
\prod_{c=1}^C \big( (\ell_c - \ell_{c-1} - |\eta^{(c-1)}|) \vee 1 \big)^{-\alpha}\\
&\qquad \leq  (C_\rho \vee 1)^C \Big( \frac{N-CM}{C}\Big)^{(\alpha+1)C} 
\prod_{c=1}^C \big( (\ell_c - \ell_{c-1} - |\eta^{(c-1)}|) \vee 1 \big)^{-\alpha}\\
&\qquad \leq \exp[ N \delta_5 ] \prod_{c=1}^C 
\big( (\ell_c - \ell_{c-1} - |\eta^{(c-1)}|) \vee 1 \big)^{-\alpha}, 
\end{aligned}
\ee  
where $\delta_5=\delta(\delta_1,M)$ can be made arbitrarily small by choosing $M$ 
large and $\delta_1$ small. For the second inequality, we have used the fact that 
the product $\prod_{c=1}^C m_c^{\alpha+1}$ is maximal when all factors are equal. 

Combining (\ref{RNform}--\ref{count4}), we obtain 
\be{count5}
\begin{aligned}
\Pr\big(R_N \in \mathcal{O} \mid X \big) 
 & \leq \exp\Big[ N \Big( H_{\tau|K}(Q) + \E_Q[\log \rho(\tau_1)] 
+ \delta_2 + \delta_3 + \delta_4 + \delta_5 \Big) \Big] \\
& \qquad \times \sum_{\mbox{\footnotesize $(\ell_c)$, $(\eta^{(c)})$ good}} 
\,\,\prod_{c=1}^C \big( (\ell_c - \ell_{c-1} - |\eta^{(c-1)}|) \vee 1 \big)^{-\alpha}.
\end{aligned} 
\ee
Combining (\ref{count5}) with Lemma~\ref{targetlemma} below, and recalling the 
identity in (\ref{eqnsre1}), we obtain the result in Proposition~\ref{propub} 
for $\rho$ satisfying (\ref{eq:algtail}), with $\mathcal{O}$ defined in (\ref{Opnedef}) 
and $\varepsilon = \delta_2+\delta_3+\delta_4+\delta_5+\delta_6$. Note that 
$\varepsilon$ can be made arbitrarily small by choosing $\varepsilon_1, \delta_1$ small 
and $M$ large.

\subsection{Step 4: Cost of finding good sentences}
\label{S3.4} 

\begin{lemma} 
\label{targetlemma} 
For $\varepsilon_1,\delta_1>0$ and $M \in \N$, 
\be{target} 
\begin{aligned}
&\limsup_{N\to\infty} \frac{1}{N} \log 
\left[\sum_{\mbox{\footnotesize $(\ell_c)$, $(\eta^{(c)})$ {\rm good}}}\,\, 
\prod_{c=1}^C \big( (\ell_c - \ell_{c-1} - |\eta^{(c-1)}|) \vee 1 \big)^{-\alpha}\right]\\ 
&\qquad \leq - \alpha \, m_Q H(\Psi_Q \mid \nu^{\otimes \N}) + \delta_6 
\quad \mbox{a.s.}, 
\end{aligned}
\ee 
where $\delta_6=\delta(\varepsilon_1,\delta_1,M)$ can be made arbitrarily small by 
choosing $M$ large and $\varepsilon_1$, $\delta_1$ small. 
\end{lemma} 

\begin{proof}
Note that, by the fourth line of (\ref{sentprop}), for any $\eta \in \mathscr{B}$
(recall (\ref{Bdef})) and $k \in \N$,
\be{Prstartest}
\Pr\big( \eta \; \mbox{starts at position} \; k \: \mbox{in} \; X\big) 
 \le \exp\big[ M \big(m_Q \E_{\Psi_Q}[\log \nu(X_1)] + \varepsilon_1\big)\big].  
\ee
Combining this with (\ref{eqnsizeofB}), we get 
\be{eq:findB} 
\begin{aligned}
\Pr\big( \mbox{some} & \; \mbox{element of $\mathscr{B}$ starts at position} 
         \; k \: \mbox{in} \; X\big) \\
&  \le \exp\big[ M \big(m_Q \E_{\Psi_Q}[\log \nu(X_1)] + \varepsilon_1\big)\big] 
       \times \exp\big[M\big(m_Q H(\Psi_Q) + \varepsilon_1\big)\big] \\
& = \exp\big[ -M \big(m_Q H(\Psi_Q \mid\nu^{\otimes\N}) - 2 \varepsilon_1\big) \big],
\end{aligned}
\ee
where we use (\ref{HPsirel}).

Next, we \emph{coarse-grain} the sequence $X$ into blocks of length 
\be{defL}
L := \lfloor M(m_Q -\varepsilon_1) \rfloor,
\ee 
and compare the coarse-grained sequence with a {\em low-density Bernoulli sequence}. 
To this end, define a $\{0,1\}$-valued sequence $(A_l)_{l\in\N}$ inductively 
as follows. Put $A_0:=0$, and, for $l \in \N$ given that $A_{0}, A_{1}, \dots, A_{l-1}$ 
have been assigned values, 
define $A_l$ by distinguishing the following two cases:
\begin{itemize}
\item[(1)]
If $A_{l-1}=0$, then 
\be{defAs21}
A_l := \left\{
\begin{array}{ll} 
1,  & \parbox{30em}{if in $X$ there is a word $\eta \in \mathscr{B}$ starting 
in $((l-1)L, lL]$,}\\[2ex]
0, & \mbox{otherwise}. 
\end{array}\right.
\ee
\item[(2)]
If $A_{l-1}=1$, then  
\be{defAs22}
A_l := \left\{
\begin{array}{ll} 
1, & \parbox{28em}{if in $X$ there are words $\eta, \eta' \in \mathscr{B}$ 
starting in $((l-2)L,(l-1)L]$, respectively, $((l-1)L,lL]$ and occurring 
disjointly,}\\[5ex]
0, \;\; & \mbox{otherwise}. 
\end{array}
\right.
\ee
\end{itemize}

Put 
\be{eq:defp} 
p:= L\,\exp\big[ -M \big(m_Q H(\Psi_Q \mid\nu^{\otimes\N}) - 2 \varepsilon_1\big) \big]. 
\ee
Then we claim
\be{eq:Adominated}
\Pr( A_1=a_1,\dots,A_n=a_n) \le p^{a_1+\cdots+a_n}, \quad 
n \in \N, \: a_1,\dots,a_n \in \{0,1\}.
\ee
In order to verify (\ref{eq:Adominated}), fix $a_1,\dots,a_n \in \{0,1\}$ with 
$a_1+\cdots+a_n=m$. By construction, for the event in the left-hand side of 
(\ref{eq:Adominated}) to occur there must be $m$ non-overlapping elements of 
$\mathscr{B}$ at certain positions in $X$. By (\ref{eq:findB}), the occurrence of 
any $m$ fixed starting positions has probability at most
\be{expkM}
\exp\big[ -m M \big(m_Q H(\Psi_Q \mid\nu^{\otimes\N}) - 2 \varepsilon_1\big) \big], 
\ee 
while the choice of the $a_l$'s dictates that there are at most $L^m$ possibilities 
for the starting points of the $m$ words. 

By (\ref{eq:Adominated}), we can couple the sequence $(A_l)_{l\in\N}$ with an i.i.d.\ 
Bernoulli($p$)-sequence $(\omega_l)_{l\in\N}$ such that 
\be{eq:Acoupledomega} 
A_l \le \omega_l \qquad \forall\,l \in \N \quad \mbox{a.s.}
\ee
(Note that (\ref{eq:Adominated}) guarantees the existence of such a coupling 
for any fixed $n$. In order to extend this existence to the infinite sequence, 
observe that the set of functions depending on finitely many coordinates is dense 
in the set of continuous increasing functions on $\{0,1\}^\N$, and use the results 
in Strassen \cite{Str65}.)

Each admissible choice of $\ell_1, \dots, \ell_C$ in (\ref{target}) leads to a 
$C$-tuple $i_1 < \dots < i_C$ such that $A_{i_1} = \cdots = A_{i_C} = 1$ (since 
it cuts out non-overlapping words, which is compatible with (\ref{defAs21}--\ref{defAs22})), 
and for any such $(i_1,\dots,i_C)$ there are at most $L^C$ different admissible
choices of the $\ell_c$'s. Thus, we have
\be{prodestells}
\sum_{\mbox{\footnotesize $(\ell_c)$, $(\eta^{(c)})$ good}}\,\, 
\prod_{c=1}^C 
\big( (\ell_c - \ell_{c-1} - |\eta^{(c-1)}|) \vee 1 \big)^{-\alpha} 
\leq L^C L^{-\alpha} \sum_{0 < i_1 < \cdots < i_C < \infty 
\atop A_{i_1} = \dots = A_{i_C}=1} 
\prod_{c=1}^C (i_{c}-i_{c-1})^{-\alpha}. 
\ee

Using 
(\ref{defG}) and recalling the definition of $\phi(\alpha,p)$
in (\ref{LDlimsup}), we have
\be{eqlimsupgood} 
\limsup_{N\to\infty} \frac{1}{N} \log\,[\,\mbox{\rm r.h.s. (\ref{prodestells})}\,]  
\leq \frac{1-\delta_1}{M}\Big( \log\big( M m_Q\big) - \phi(\alpha, p) \Big)
\qquad (\omega, A)-a.s.  
\ee
From (\ref{eq:defp}) we know that $\log(1/p) \sim M (m_Q H(\Psi_Q \mid\nu^{\otimes\N}) 
- 2 \varepsilon_1)$ as $M\to\infty$ and so, by Lemma~\ref{lowtoylemma}, we have
\be{rhsest1} 
\mbox{\rm r.h.s. (\ref{eqlimsupgood})} \leq 
 -(1-\varepsilon_2) \alpha \big(m_Q H(\Psi_Q \mid\nu^{\otimes\N}) - 2\varepsilon_1\big)
\ee
for any $\varepsilon_2 \in (0, 1)$, provided $M$ is large enough. This completes the proof of 
Lemma~\ref{targetlemma}, and hence of Proposition~\ref{propub} for $Q \in 
\mathcal{P}^{\mathrm{erg,fin}}(\widetilde{E}^\N)$.
\end{proof}

\subsection{Step 5: Removing the assumption of ergodicity}
\label{S3.5} 

Sections~\ref{S3.1}--\ref{S3.4} contain the main ideas behind the proof of 
Proposition~\ref{propub}. In the present section we extend the bound from
$\mathcal{P}^{\mathrm{erg,fin}}(\widetilde{E}^\N)$ to $\mathcal{P}^{\mathrm{inv,fin}}
(\widetilde{E}^\N)$. This requires setting up a variant of the argument in 
Sections~\ref{S3.1}--\ref{S3.4} in which the ergodic components of $Q$ are
``approximated with a common length scale on the letter level''. This turns 
out to be technically involved and to fall apart into 6 substeps. 

Let $Q \in \mathcal{P}^{\mathrm{inv,fin}}(\widetilde{E}^\N)$ have 
a non-trivial ergodic decomposition 
\be{}
Q = \int_{\mathcal{P}^{\mathrm{erg}}(\widetilde{E}^\N)} Q' \, W_Q(dQ'), 
\ee 
where $W_Q$ is a probability measure on $\mathcal{P}^{\mathrm{erg}}(\widetilde{E}^\N)$
(Georgii~\cite{Ge88}, Proposition 7.22). We may assume w.l.o.g.\ that 
$H(Q \mid q_{\rho,\nu}^{\otimes\N})< \infty$, otherwise we can simply employ the 
annealed bound. Thus, $W_Q$ is in fact supported on $\mathcal{P}^{\mathrm{erg,fin}}
(\widetilde{E}^\N) \cap \{Q'\colon\, H(Q' \mid q_{\rho,\nu}^{\otimes\N})<\infty\}$. 

Fix $\varepsilon > 0$. In the following steps, we will construct an open neighbourhood 
$\mathcal{O}(Q) \subset \mathcal{P}^{\mathrm{inv}}(\widetilde{E}^\N)$ of $Q$ 
satisfying (\ref{mainub}) (for technical reasons with $\varepsilon$ replaced by some 
$\varepsilon'=\varepsilon'(\varepsilon)$ that becomes arbitrarily small as $\varepsilon 
\downarrow 0$). 

\subsubsection{Preliminaries}

Observing that
\be{}
m_Q = \int_{\mathcal{P}^{\mathrm{erg}}(\widetilde{E}^\N)} m_{Q'}\,W_Q(dQ')<\infty, 
\qquad
H(Q|q_{\rho,\nu}^{\otimes\N}) = \int_{\mathcal{P}^{\mathrm{erg}}(\widetilde{E}^\N)} 
H(Q'|q_{\rho,\nu}^{\otimes\N})\,W_Q(dQ') < \infty,
\ee 
we can find $K_0,K_1,m^*>0$ and a compact set
\be{}
\mathscr{C} \subset \mathcal{P}^{\mathrm{inv}}(\widetilde{E}^\N) 
\cap \mathrm{supp}(W_Q) \cap \{Q\colon\, H(\cdot|q_{\rho,\nu}^{\otimes\N}) 
\leq K_0\}
\ee 
such that 
\begin{eqnarray} 
&& \sup \{H(\Psi_P \mid \nu^{\otimes \N}) 
\colon\, P \in \mathscr{C} \} \le K_1, \\
&&\sup \{m_P \colon\, P \in \mathscr{C} \} \le m^*,\\
&&\label{eq:unif_integrable}
\mbox{the family} \; \{ \mathscr{L}_P(\tau_1)\colon\,P \in \mathscr{C} \} \; 
\mbox{is uniformly integrable}, \\
\label{eq:WQClarge}
&& W_Q(\mathscr{C}) \ge 1-\varepsilon/2, \\ 
\label{eq:intCHQlarge}
&& \int_{\mathscr{C}} H(Q'|q_{\rho,\nu}^{\otimes\N}) \, W_Q(dQ') 
\ge H(Q|q_{\rho,\nu}^{\otimes\N})-\varepsilon/2, \\
\label{eq:intCHPsiQlarge}
&& \int_{\mathscr{C}} m_{Q'} H(\Psi_{Q'}|\nu^{\otimes\N})\, W_Q(dQ') 
\ge m_Q H(\Psi_{Q}|\nu^{\otimes\N})-\varepsilon/2.
\end{eqnarray}
In order to check (\ref{eq:unif_integrable}), observe that $\E_Q[\tau_1]<\infty$ implies 
that there is a sequence $(c_n)$ with $\lim_{n\to\infty} c_n=\infty$ such that 
\be{}
\E_Q\big[\tau_1 \1_{\{\tau_1 \ge c_n\}}\big] \le 
\frac{6}{\pi^2 n^3}\,\frac{\varepsilon}{6}, 
\quad n \in \N. 
\ee
Put 
\be{}
\widehat{A}_n := \{ Q' \in \mathcal{P}^{\mathrm{inv}}(\widetilde{E}^\N)\colon\, 
\E_{Q'}\big[\tau_1 \1_{\{\tau_1 \ge c_n\}}\big] > 1/n\}
\ee 
and $A:=\cap_{n\in\N} (\widehat{A}_n)^c$. Each $\widehat{A}_n$ is open, 
hence $A$ is closed, and by the Markov inequality we have 
\be{}
W_Q\Big(\big\{ Q'\colon\,\E_{Q'}\big[\tau_1 \1_{\{\tau_1 \ge c_n\}}\big] > 1/n 
   \big\}\Big)
\leq n \E_Q\big[\tau_1 \1_{\{\tau_1 \ge c_n\}}\big] 
\leq \frac{6}{\pi^2 n^2}\,\frac{\varepsilon}{6}.
\ee
Thus, 
\be{}
W_Q(A^c) =  W_Q\big( {\textstyle \cup_{n\in\N} \widehat{A}_n  } \big) 
\le \frac{\varepsilon}{6} \sum_{n\in\N} \frac{6}{\pi^2 n^2} 
= \frac{\varepsilon}{6}.
\ee
This implies that the mapping 
\be{eq:mQHPsiQ.cont} 
Q' \mapsto m_{Q'} H(\Psi_{Q'}|\nu^{\otimes \N}) 
\quad \mbox{is lower semicontinuous on $\mathscr{C}$}. 
\ee
Indeed, if ${\rm w}-\lim_{n\to\infty} Q'_n = Q''$ and $(Q'_n) \subset \mathscr{C}$,
then $\lim_{n\to\infty} \E_{Q'_n}[\tau_1]=\lim_{n\to\infty} m_{Q'_n}=m_{Q''}=\E_{Q''}[\tau_1]$ 
and ${\rm w}-\lim_{n\to\infty} \Psi_{Q'_n}=\Psi_{Q''}$ by uniform integrability  (see Birkner~
\cite{Bi08}, Remark~7). 
\smallskip 

Furthermore, we can find $N_0, L_0 \in \N$ with $L_0 \leq N_0$ and a finite set $\widetilde{W} 
\subset \widetilde{E}^{N_0}$ such that the following holds. Let 
\be{} 
W:= \Big\{ \pi_{L_0} (\theta^i \kappa(\zeta))\colon\, 
\zeta=(\zeta^{(1)},\dots,\zeta^{(N_0)}) \in \widetilde{W}, 0 \leq i < |\zeta^{(1)}| \Big\}
\ee
be the set of words of length $L_0$ obtained by concatenating sentences from $\widetilde{W}$, 
possibly shifting the ``origin'' inside the first word and restricting to the first $L_0$ letters.
Then, denoting by $\mathscr{D}$ the set of all 
$P\in\mathcal{P}^{\mathrm{inv,fin}}(\widetilde{E}^\N)\cap\mathscr{C}$
that satisfy  
\be{eq:Wtildedef1}
\sum_{\zeta \in \widetilde{W}} P(\zeta) \geq 1 - \frac{\varepsilon}{3 c_{\lceil 3/\varepsilon \rceil}}, 
\quad \forall\, \xi \in W : \Psi_P(\xi) \leq 
\frac{1+\varepsilon/2}{m_P} \E_P \Big[ \1_{\widetilde{W}}(\pi_{N_0} Y)
  \sum_{i=0}^{\tau_1-1}\1_{\{\xi\}}(\pi_{L_0}\theta^i\kappa(Y))\Big]  
\ee 
\begin{align} 
\label{eq:HP.bd} 
H(P \mid q_{\rho,\nu}^{\otimes \N}) + \varepsilon/4 \geq
\frac{1}{N_0} \sum_{\zeta \in \widetilde{W}} P(\zeta) 
\log\frac{P(\zeta)}{q_{\rho,\nu}^{\otimes N_0}(\zeta)} \geq & 
\: H(P \mid q_{\rho,\nu}^{\otimes \N}) 
- \varepsilon/4,\\ 
\label{eq:HPsiP.bd}
m_P H(\Psi_P \mid \nu^{\otimes\N}) + \varepsilon/4 \geq
\frac{m_P}{L_0} \sum_{w \in W} \Psi_P(w) \log\frac{\Psi_P(w)}{\nu^{\otimes L_0}(w)}
\geq & \: m_P H(\Psi_P \mid \nu^{\otimes\N}) - \varepsilon/4,
\end{align}
we can choose $N_0$, $L_0$ and $\widetilde{W}$ so large that 
the following inequalities hold:
\begin{eqnarray} 
W_Q(\mathscr{D}) & \geq & 1 - 3\varepsilon/4,\\
\int_{\mathscr{D}} H(P \mid q_{\rho,\nu}^{\otimes\N}) \, W_Q(dP) 
& \geq & H(Q \mid q_{\rho,\nu}^{\otimes\N})- 3\varepsilon/4, \\
\int_{\mathscr{D}} m_{P} H(\Psi_{P}|\nu^{\otimes\N})\, W_Q(dP) 
& \geq & m_Q H(\Psi_{Q} \mid \nu^{\otimes\N})-3\varepsilon/4.
\end{eqnarray}
We may choose the set $\widetilde{W}$ in such a way that 
\be{}
\delta_{\widetilde{W}} := 
\min\{ q_{\rho,\nu}^{\otimes N_0}(\zeta)\colon \zeta \in \widetilde{W} \} 
\cdot
\frac{\min\{ \nu^{\otimes L_0}(\xi)\colon\,\xi \in W \}}
{\max\{ |\zeta^{(1)}|\colon\,\zeta \in \widetilde{W} \}}
\cdot \frac1{|\widetilde{W}|}
> 0. 
\ee

\subsubsection{Approximating with a given length scale on the letter level}
\label{ss.sect:appr.given.length.scale}

\newcommand{\lepsilon}{\delta}

For $P \in\mathcal{P}^{\mathrm{inv,fin}}(\widetilde{E}^\N)$, we put 
\be{eq:deltaPW} 
\delta_{P,\widetilde{W}}:= 
\delta_{\widetilde{W}} \cdot
\left( 
\min\big\{ P(\zeta)\colon\,\zeta\in \widetilde{W}, P(\zeta)>0 \big\} \wedge 
\min\big\{ \Psi_P(\xi)\colon\, \xi \in W, \Psi_P(\xi) > 0 \big\} 
\right).
\ee
For $\lepsilon > 0$ and $L \in \N$, we say that $P \in\mathcal{P}^{\mathrm{inv,fin}}
(\widetilde{E}^\N)$ can be $(\lepsilon, L)$-{\em approximated} if there exists a finite 
subset $\mathscr{A}_P \subset \widetilde{E}^{\lceil L/m_P\rceil}$ of ``$P$-typical'' 
sentences, each consisting of $\approx L/m_P$ words 
(we assume that $L > N_0 m_P$), such that
\be{eq:Ptyp1} 
P_{|_{\mathscr{F}_{\lceil L/m_P\rceil}}}(\mathscr{A}_P) 
> 1- \frac12 \delta \cdot \delta_{P,\widetilde{W}}
\ee
and, for all $z=(y^{(1)},\dots,y^{(\lceil L/m_P\rceil)}) \in \mathscr{A}_P$, 
\be{eq:Ptyp2}
\begin{aligned} 
P(z) & \in \Big[\exp\big[-\lceil L/m_P\rceil (H(Q)+\lepsilon)\big], 
\exp\big[-\lceil L/m_P\rceil (H(Q)-\lepsilon)\big]\Big], \\
|\kappa(z)| & \in [L(1-\lepsilon),L(1+\lepsilon)], \\
P\big(K^{(\lceil L/m_P \rceil)}=z\big) & \in 
\Big[\exp\big[-L (H(\Psi_Q)+\lepsilon)\big], 
\exp[-L (H(\Psi_Q)-\lepsilon)\big]\Big], \\
\sum_{k=1}^{|\kappa(z)|} \log \nu(\kappa(z)_k) 
& \in [L(1-\lepsilon),L(1+\lepsilon)]\,\,\E_{\Psi_P}\big[\log\nu(X_1)\big],\\
\sum_{i=1}^{\lceil L/m_P\rceil} \log \rho(|y^{(i)}|) 
& \in [(L/m_P)(1-\lepsilon),(L/m_P)(1+\lepsilon)]\,\,\E_P\big[\log\rho(\tau_1)\big], \\
|\{ z' \in \mathscr{A}_P\colon\,\kappa(z)=\kappa(z') \}| & \leq 
\exp\big[(L/m_P) (H_{\tau|K}(P) + \lepsilon)\big]. 
\end{aligned}
\ee 
By the third and the 
fourth
line of (\ref{eq:Ptyp2}) we have, using (\ref{HPsirel}), 
\be{}
\Pr\big(X\;\mbox{starts with some element of}\;\kappa(\mathscr{A}_P)\big) 
\leq \exp\Big[-L(1-2\lepsilon) H\big(\Psi_Q \mid \nu^{\otimes\N}\big)\Big]. 
\ee
For $P$ that can be $(\lepsilon, L)$-approximated, define an open neighbourhood of 
$P$ via
\be{def:UPepsL}
\mathcal{U}_{(\lepsilon, L)}(P) 
:= \left\{P' \in \mathcal{P}^{\mathrm{inv}}(\widetilde{E}^\N) :
\frac{P'(z)}{P(z)} \in (1-\lepsilon \cdot \delta_{P,\widetilde{W}},1+\lepsilon \cdot 
\delta_{P,\widetilde{W}}) \;\;
\forall \, z \in \mathscr{A}_P \right\}, 
\ee
where $\mathscr{A}_P=\mathscr{A}_P(\lepsilon,L)$ is the set from (\ref{eq:Ptyp1}--\ref{eq:Ptyp2}). 
By the results of Section~\ref{S3.1} and the above, for given $P\in\mathcal{P}^{\mathrm{erg,fin}}
(\widetilde{E}^\N) \cap \mathscr{C}$ and $\lepsilon_0>0$ there exist $\lepsilon' \in (0,
\lepsilon_0)$ and $L'$ such that 
\be{Pcanap}
\forall \, L'' \ge L' \colon\,\mbox{ $P$ can be $(\lepsilon', L'')$-approximated.} 
\ee

Assume that a given $P \in \mathscr{D}$ can be $(\lepsilon, L)$-approximated for some $L$ such 
that $\lceil L /m_P \rceil \ge N_0$. We claim that then for any $P' \in \mathscr{D} \cap 
\mathcal{U}_{(\lepsilon, L)}(P)$, 
\begin{align}
\label{eq:Pzeta.approxprimeAp} 
P'(\widetilde{E}^{\lceil L /m_P \rceil} \setminus \mathscr{A}_P) 
\le & \hspace{0.5em} 2 \lepsilon \cdot \delta_{P,\widetilde{W}},
\\
\label{eq:Pzeta.approx}
\forall\, \zeta \in \widetilde{W}\colon\,
P'(\zeta) \leq &  
\begin{cases} (1+3\lepsilon) P(\zeta) & \mbox{if $P(\zeta)>0$}, \\
2 \lepsilon \cdot \delta_{P,\widetilde{W}} 
\: (\leq 2\lepsilon \cdot 
\min\{ q_{\rho,\nu}^{\otimes N_0}(\zeta') \colon\, \zeta' \in \widetilde{W} \}) & \mbox{otherwise}, 
\end{cases} \\ 
\label{eq:PsiPxi.approx}
\forall\, \xi \in W\colon\,
m_{P'} \Psi_{P'}(\xi) \leq &
\begin{cases} (1+\varepsilon/2) 
(1+3\lepsilon) m_P \Psi_P(\xi) & \mbox{if $\Psi_P(\xi)>0$}, \\
(1+\varepsilon/2) 2\lepsilon 
\min\{ \nu^{\otimes L_0}(\xi') \colon\, \xi' \in W \} & \mbox{otherwise}, 
\end{cases} \\ 
\label{eq:mPprime.approx}
m_{P'} \geq & (1-3\lepsilon)(m_P - \varepsilon) \; (\geq (1-3\lepsilon-\varepsilon) m_P).  
\end{align} 
(\ref{eq:Pzeta.approxprimeAp}) follows from (\ref{eq:Ptyp1}) and (\ref{def:UPepsL}). 
To verify (\ref{eq:Pzeta.approx}), note that, for $\zeta \in \widetilde{W}$, 
\be{}
\begin{aligned}
P'(\zeta) & \le 
\sum_{z \in \mathscr{A}_P \colon\, \pi_{N_0}(z)=\zeta} \hspace{-2em} P'(z) \; + 
\sum_{z \in \widetilde{E}^{\lceil L/m_P\rceil} \setminus \mathscr{A}_P \colon\, \pi_{N_0}(z)=\zeta} 
\hspace{-3em} P'(z) \\ 
& \le (1+\lepsilon) \hspace{-1em} 
\sum_{z \in \mathscr{A}_P : \pi_{N_0}(z)=\zeta} \hspace{-2em} P(z) \; + 
P'\big(\widetilde{E}^{\lceil L/m_P\rceil} \setminus \mathscr{A}_P \big) 
\end{aligned} 
\ee
and use (\ref{eq:Pzeta.approxprimeAp}) on the last term in the second line, 
observing that $\delta_{P,\widetilde{W}} \leq P(\zeta)$ whenever $\zeta \in \widetilde{W}$ 
and $P(\zeta)>0$.  
To verify (\ref{eq:PsiPxi.approx}), observe that, for $\xi \in W$ (recall the definition of $\Psi_{P'}$ 
from (\ref{PsiQdef})), using (\ref{eq:Wtildedef1}), 
\be{}
\begin{aligned}
(1+\varepsilon/2)^{-1} m_{P'} \Psi_{P'}(\xi) & \leq 
\sum_{\zeta \in \widetilde{W}} P'(\zeta) \sum_{i=0}^{|\zeta^{(1)}|-1} 
\1_{\{\xi\}}\big( \pi_{L_0} ( \theta^{i} \kappa(\zeta) ) \big) \\
& \leq  (1+\lepsilon) m_P \Psi_P(\xi) + 
\sum_{\zeta \in \widetilde{W} \colon\, P(\zeta)=0} |\zeta^{(1)}| P'(\zeta) 
\end{aligned}
\ee
and that the sum in the second line above is bounded by 
$|\widetilde{W}| \cdot \max\{ |\zeta^{(1)}| : \zeta \in \widetilde{W} \} 
\cdot 2\lepsilon \cdot \delta_{P,\widetilde{W}}$, 
which is not more than $2\lepsilon m_P \Psi_P(\xi)$ if $\Psi_P(\xi)>0$ and not 
more than $2\delta \min\{ \nu^{\otimes L_0}(\xi') \colon\, \xi' \in W \}$ otherwise. 
Lastly, to verify (\ref{eq:mPprime.approx}), note that 
\be{}
P'(\zeta) \ge (1-3\lepsilon) P(\zeta) \;\;\; \forall \, \zeta \in \widetilde{W}
\ee
(which can be proved in the same way as (\ref{eq:Pzeta.approx})), so that
\be{eq:mPprime.approx2}
m_{P'} = \sum_{y\in\widetilde{E}} |y| P'(y) \ge 
\sum_{\zeta \in \widetilde{W}} |\zeta^{(1)}| P'(\zeta) 
\ge (1-3\lepsilon) 
\sum_{\zeta \in \widetilde{W}} |\zeta^{(1)}| P(\zeta). 
\ee
Furthermore,  
\be{uppu}
m_P \le \sum_{\zeta \in \widetilde{W}} |\zeta^{(1)}| P(\zeta) 
+ c_{\lceil 3/\varepsilon\rceil} P\big(\widetilde{E}^{N_0} \setminus \widetilde{W} \big) 
+ \sum_{y \in \widetilde{E}\colon\, |y|> c_{\lceil 3/\varepsilon\rceil}} |y| P(y). 
\ee
Observing that the second and the third term on the right-hand side are each at most $\varepsilon/3$, 
we find that (\ref{eq:mPprime.approx2}--\ref{uppu}) imply (\ref{eq:mPprime.approx}). 
\smallskip

Finally, observe that (\ref{eq:Pzeta.approx}--\ref{eq:mPprime.approx}) imply that 
there exists $\delta_0 \,(=\delta_0(\varepsilon)) \, > 0$ with the following property: 
For any $P, P' \in \mathscr{D}$ such that $P$ can be $(\lepsilon, L)$-approximated 
for some $L$ with $\lceil L /m_P \rceil \ge N_0$ and $\delta \leq \delta_0$ 
and $P' \in \mathcal{U}_{(\lepsilon, L)}(P)$, we have 
\begin{eqnarray} 
\label{eq:HPprime.bd1}
H(P' \mid q_{\rho,\nu}^{\otimes \N}) & \le & 
(1+\varepsilon) \Big( H(P \mid q_{\rho,\nu}^{\otimes \N}) + \varepsilon \Big) 
\quad \mbox{and}
\\
\label{eq:HPsiPprime.bd1}
m_{P'} H(\Psi_{P'} \mid \nu^{\otimes \N}) & \le & 
(1+\varepsilon) \Big( m_{P} H(\Psi_{P} \mid \nu^{\otimes \N}) + \varepsilon \Big).
\end{eqnarray}
Here, (\ref{eq:HPprime.bd1}) follows from the observation 
\begin{equation}
\begin{aligned} 
&H(P' \mid q_{\rho,\nu}^{\otimes \N}) - \frac{\varepsilon}{4}\\ 
&\leq 
\frac{1}{N_0} \sum_{\zeta \in \widetilde{W}} P'(\zeta) 
\log\frac{P'(\zeta)}{q_{\rho,\nu}^{\otimes N_0}(\zeta)} \\
&\leq \frac{1+3\lepsilon}{N_0} \sum_{\zeta \in \widetilde{W}} P(\zeta) 
\log\frac{(1+3\lepsilon)P(\zeta)}{q_{\rho,\nu}^{\otimes N_0}(\zeta)} 
+ \frac{1}{N_0} \sum_{\zeta \in \widetilde{W} : P(\zeta)=0} \hspace{-1.5em} P'(\zeta) 
\log\frac{\min\{ q_{\rho,\nu}^{\otimes N_0}(\zeta') \colon \zeta' \in \widetilde{W} \}}
{q_{\rho,\nu}^{\otimes N_0}(\zeta)} \\
&\leq  (1+3\lepsilon) \Big( H(P \mid q_{\rho,\nu}^{\otimes \N}) + \frac{\varepsilon}{4} \Big)
+ \frac{1+3\lepsilon}{N_0} \log(1+3\lepsilon). 
\end{aligned}
\end{equation}
Similarly, observing that 
\begin{equation} 
\begin{aligned} 
m_{P'} & \sum_{\xi \in W} \Psi_{P'}(\xi) \log 
\frac{m_{P'}\Psi_{P'}(\xi)}{m_{P'} \nu^{\otimes L_0}(\xi)} \\
& \leq 
\Big(1+\frac{\varepsilon}{2}\Big)(1+3\lepsilon) m_P \sum_{\xi \in W} \Psi_P(\xi) 
\log\frac{(1+\varepsilon/2)(1+3\lepsilon)m_P \Psi_P(\xi)}{(1-3\lepsilon-\varepsilon)
m_P \nu^{\otimes L_0}(\xi)} \\
& \hspace{2em} {} 
+ m_{P'} \hspace{-1.5em} \sum_{\xi \in W \colon\, \Psi_P(\xi)=0} \hspace{-1.5em} \Psi_{P'}(\xi)
\log\frac{(1+\varepsilon/2)2\lepsilon \min\{ \nu^{\otimes L_0}(\xi') \colon \xi' \in W \}}{\nu^{\otimes L_0}(\xi)}
\\ 
& \leq \Big(1+\frac{\varepsilon}{2}\Big)(1+3\lepsilon) L_0 m_P H(\Psi_P\mid\nu^{\otimes \N}) 
+ \varepsilon/2+ m^* \log\frac{(1+3\lepsilon) (1+\varepsilon/2)}{1-3\lepsilon-\varepsilon} 
\end{aligned}
\end{equation}
we obtain (\ref{eq:HPsiPprime.bd1}) in view of (\ref{eq:HPsiP.bd}).

\subsubsection{Approximating the ergodic decomposition}
 
In the previous subsection, we have approximated a given 
$P \in \mathcal{P}^{\mathrm{erg,fin}}(\widetilde{E}^\N)$, 
i.e., we have constructed a certain neighbourhood of $P$ w.r.t.\ the weak topology, which 
requires only conditions on the frequencies of sentences whose concatenations are $\approx L$ 
letters long. While the required $L$ will in general vary with $P$, we now want to construct a 
compact $\mathscr{C}'\subset\mathscr{C}$ such that $W_Q(\mathscr{C}')$ is still close to $1$ 
and all $P \in \mathscr{C}'$ can be approximated on the same scale $L$ (on the letter level). 
To this end, let 
\be{}
\mathscr{D}_{\varepsilon', L'} := 
\big\{ P\in \mathscr{D} \colon\, 
\mbox{$P$ can be $(\varepsilon',L')$-approximated} \big\}.
\ee
By (\ref{Pcanap}), we have
\be{}
\bigcup_{ {\varepsilon' \in (0,\varepsilon/2)} \atop {L' \in \N}} \mathscr{D}_{\varepsilon', L'} 
= \mathcal{P}^{\mathrm{erg,fin}}(\widetilde{E}^\N) \cap \mathscr{C}, 
\ee
so, in view of (\ref{eq:WQClarge}--\ref{eq:intCHPsiQlarge}), we can 
choose 
\be{}
0 < \varepsilon_1 
< \frac{\varepsilon}{2 m^* (1 \vee K_1)} \wedge \frac{\delta_0}{2}
\ee
and $L \in \N$ such that 
\begin{eqnarray}
W_Q(\mathscr{D}_{\varepsilon_1, L}) & \ge & 1-\varepsilon, \\
\label{eq:intHQD1big}
\int_{\mathscr{D}_{\varepsilon_1, L}} H(Q'\mid q_{\rho,\nu}^{\otimes\N}) \, W_Q(dQ') 
& \ge & H(Q \mid q_{\rho,\nu}^{\otimes\N})-\varepsilon, \\
\label{eq:intHPsiD1big}
\int_{\mathscr{D}_{\varepsilon_1, L}} m_{Q'} H(\Psi_{Q'}\mid \nu^{\otimes\N})\, W_Q(dQ') 
& \ge & m_Q H(\Psi_{Q}\mid \nu^{\otimes\N})-\varepsilon.
\end{eqnarray}
For $P \in \mathscr{D}_{\varepsilon_1, L}$, let   
\be{}
\mathcal{U}'(P) := 
\left\{P' \in \mathcal{P}^{\mathrm{inv}}(\widetilde{E}^\N)\colon\,
\frac{P'(z)}{P(z)} \in \left(1-\frac{\varepsilon_1}{2} \delta_{P,\widetilde{W}},
1+\frac{\varepsilon_1}{2} \delta_{P,\widetilde{W}} \right) \;\;
\forall \, z \in \mathscr{A}_P \right\}, 
\ee
where $\mathscr{A}_P$ is the set from (\ref{eq:Ptyp1}--\ref{eq:Ptyp2}) that appears in the 
definition of $\mathcal{U}_{(\varepsilon_1, L)}(P)$ and 
$\delta_{P,\widetilde{W}}$ is defined in (\ref{eq:deltaPW}). Note that $\mathcal{U}'(P) \subset 
\mathcal{U}_{(\varepsilon_1, L)}(P)$. Indeed, $\inf_{P \in \mathscr{D}_{\varepsilon_1, L}}
\mathrm{dist}(\mathcal{U}'(P),\mathcal{U}_{(\varepsilon_1, L)}(P)^c)>0$ if we metrize the 
weak topology. Consequently, 
\be{}
\mathscr{C}':= 
\mathscr{C} \cap \overline{\cup_{P \in \mathscr{D}_{\varepsilon_1, L}} \mathcal{U}'(P)} 
\; \big( \supset \mathscr{D}_{\varepsilon_1, L} \big)
\ee
is compact and satisfies $W_Q(\mathscr{C}') \ge 1-\varepsilon$, and 
\be{}
\mathscr{C}' \subset \bigcup_{P \in \mathscr{D}_{\varepsilon_1, L}} 
\mathcal{U}_{(\varepsilon_1, L)}(P)
\ee
is an open cover. By compactness there exist $R\in\N$ and (pairwise different) 
$Q_1,\dots,Q_R \in \mathcal{P}^{\mathrm{erg,fin}}(\widetilde{E}^\N) \cap \mathscr{C}$
such that 
\be{}
\mathcal{U}_{(\varepsilon_1, L)}(Q_1) \cup \cdots \cup 
\mathcal{U}_{(\varepsilon_1, L)}(Q_R) \supset \mathscr{C}', 
\ee
where $\mathcal{U}_{(\varepsilon_1, L)}(Q_r)$ is of the type 
(\ref{def:UPepsL}) with a set $\mathscr{A}_r \subset \widetilde{E}^{M_r}$ 
satisfying (\ref{eq:Ptyp1}--\ref{eq:Ptyp2}) with $P$ replaced by $Q_r$, 
and $M_r=\lceil L/m_{Q_r} \rceil$. 

For $z \in \cup_{n\in\N}\widetilde{E}^n$ consider the 
probability measure on $[0,1]$ given  by $\mu_{Q,z}(B) := 
W_Q(\{Q' \in \mathcal{P}^{\mathrm{erg,fin}}(\widetilde{E}^\N) \colon\, Q'(z) \in B \})$, 
$B \subset [0,1]$ measurable. Observing that 
\be{}
\bigcup_{r=1}^R \bigcup_{z \in \mathscr{A}_r} 
\big\{ u \in [0,1]\colon\,\mbox{$u$ is an atom of $\mu_{Q,z}$}\big\}
\ee
is at most countable, we can find $\varepsilon_2 \in 
[\varepsilon_1, \varepsilon_1+\varepsilon_1^2)$ (note that still 
$\varepsilon_2 < 2\varepsilon_1$) and $\widetilde{\delta}>0$ 
such that 
\be{}
\begin{aligned} 
&W_Q\left(
\left\{ Q' \in \mathcal{P}^{\mathrm{erg,fin}}(\widetilde{E}^\N)\colon\, 
\begin{array}{cc} 
{Q'(z)}/{Q_r(z)} \in 
[1-(\varepsilon_2+\widetilde{\delta}) \delta_{Q_r,\widetilde{W}},
1-(\varepsilon_2-\widetilde{\delta}) \delta_{Q_r,\widetilde{W}}] \; \mbox{or} \\
{Q'(z)}/{Q_r(z)} \in 
[1+(\varepsilon_2-\widetilde{\delta}) \delta_{Q_r,\widetilde{W}},
1+(\varepsilon_2+\widetilde{\delta}) \delta_{Q_r,\widetilde{W}}] \\ 
\mbox{for some $r \in \{1,\dots,R\}$ and $z \in \mathscr{A}_r$}
\end{array}
\right\}\right) \\
&\qquad\qquad \le \frac{\varepsilon}{1 \vee K_0 \vee m^* K_1}.
\end{aligned}
\ee
Define ``disjointified'' versions of the $\mathcal{U}_{(\varepsilon, L)}(Q_r)$
as follows. For $r=1,\dots,R$, put iteratively 
\be{}
\widetilde{\mathcal{U}}_r := 
\left\{ Q' \in \mathcal{P}^{\mathrm{inv}}(\widetilde{E}^\N)\colon\, 
\begin{array}{l} 
Q'(z) \in Q_r(z) (1-\varepsilon_2 \delta_{Q_r,\widetilde{W}}, 
1+\varepsilon_2 \delta_{Q_r,\widetilde{W}}) \;\;
\mbox{for all $z \in \mathscr{A}_r$} \\
\mbox{and for each $r'<r$ there 
is $z' \in \mathscr{A}_{r'}$ such 
that} \\ \; 
Q'(z') \not\in Q_{r'}(z') [1-(\varepsilon_2+\widetilde{\delta}) \delta_{Q_{r'},\widetilde{W}},
1+(\varepsilon_2+\widetilde{\delta}) \delta_{Q_{r'},\widetilde{W}}]
\end{array}
\right\}. 
\ee
It may happen that some of the $\widetilde{\mathcal{U}}_r$ are empty or satisfy 
$W_Q(\widetilde{\mathcal{U}}_r)=0$. We then (silently) remove these and re-number 
the remaining ones. Note that each $\widetilde{\mathcal{U}}_r$ is an open subset 
of $\mathcal{P}^{\mathrm{inv}}(\widetilde{E}^\N)$ and 
\be{eq:WQUtildeslarge}
W_Q\big( \cup_{r=1}^R \widetilde{\mathcal{U}}_r \big) 
= \sum_{r=1}^R W_Q(\widetilde{\mathcal{U}}_r) \ge 1-2\varepsilon, 
\ee
since $W_Q\big(\mathscr{C}' \setminus \cup_{r=1}^R \widetilde{\mathcal{U}}_r\big) 
\leq \varepsilon$. 
\medskip

\label{missing}
For $r=1,\dots,R$, we have, using (\ref{eq:HPprime.bd1}--\ref{eq:HPsiPprime.bd1}) 
and the choice of $\varepsilon_2 \, (\le 2 \varepsilon_1 \leq \delta_0)$, 
\begin{eqnarray} 
W_Q(\widetilde{\mathcal{U}}_r \cap \mathscr{D}) 
\Big( H(Q_r \mid q_{\rho,\nu}^{\otimes\N}) + \varepsilon \Big) 
& \geq & \frac1{1+\varepsilon} \int_{\widetilde{\mathcal{U}}_r \cap \mathscr{D}} 
H(Q' \mid q_{\rho,\nu}^{\otimes\N}) \, W_Q(dQ') , \\
W_Q(\widetilde{\mathcal{U}}_r \cap \mathscr{D}) 
\Big( m_{Q_r} H(\Psi_{Q_r} \mid {\nu}^{\otimes\N}) + \varepsilon \Big) 
& \geq & \frac1{1+\varepsilon}
\int_{\widetilde{\mathcal{U}}_r \cap \mathscr{D}} m_{Q'} 
H(\Psi_{Q'} \mid {\nu}^{\otimes\N}) \, W_Q(dQ'),
\end{eqnarray}
so that altogether, using (\ref{eq:intHQD1big}--\ref{eq:intHPsiD1big}), 
\be{eq:Ifinappr1}
\begin{aligned}
\sum_{r=1}^R W_Q(\widetilde{\mathcal{U}}_r) & 
\Big\{H(Q_r \mid q_{\rho,\nu}^{\otimes\N}) + (\alpha-1)m_{Q_r}H(\Psi_{Q_r}\mid \nu^{\otimes \N})\Big\}\\ 
&\geq 
\frac1{1+\varepsilon}
\Big( H(Q \mid q_{\rho,\nu}^{\otimes\N}) + (\alpha-1) m_Q H(\Psi_{Q} \mid \nu^{\otimes \N}) 
\Big) 
- 2\alpha\varepsilon. 
\end{aligned}
\ee

\subsubsection{More layers: long sentences with the right pattern frequencies}

For $z \in \cup_{n\in\N} \widetilde{E}^n$ and $\xi = (\xi^{(1)},\dots,\xi^{(\widetilde{M})}) 
\in \widetilde{E}^M$ (with $M>|z|)$, let 
\be{} 
\mathrm{freq}_{z}(\xi) = 
\frac{1}{M} \big| \{ 1 \le i \le M-|z|+1\,\colon\, 
(\xi^{(i)},\dots,\xi^{(i+|z|-1)}) = z \} \big| 
\ee
be the empirical frequency of $z$ in $\xi$. Note that, for any $P \in \mathcal{P}^{\mathrm{erg,fin}}
(\widetilde{E}^\N)$, $z \in \cup_{n\in\N}\widetilde{E}^n$ and $\varepsilon'>0$, we have 
\be{eq:Pfreq.right}
\lim_{M\to\infty} 
P\Big(\big\{\xi \in \widetilde{E}^M\colon\, 
\mathrm{freq}_{z}(\xi) \in P(z) (1-\varepsilon',1+\varepsilon') \big\}\Big) = 1
\ee
and 
\be{eq:Pmwl.right}
\lim_{M\to\infty} 
P\Big(\big\{\xi \in \widetilde{E}^M\colon\, 
|\kappa(\xi)| \in M(m_P-\varepsilon',m_P+\varepsilon')\big\}\Big) = 1.
\ee
For $\widetilde{M} \in \N$ and $r \in \{1,\dots,R\}$, put 
\be{}
V_{r,\widetilde{M}} := 
\left\{ \xi \in \widetilde{E}^{\widetilde{M}}\colon\,
\begin{array}{l} 
|\kappa(\xi)| \in \widetilde{M}(m_{Q_r}-\varepsilon_2,m_{Q_r}+\varepsilon_2), \\
\mathrm{freq}_z(\xi) \in Q_r(z) 
(1-\varepsilon_2\delta_{Q_r,\widetilde{W}},1+\varepsilon_2\delta_{Q_r,\widetilde{W}}) \;
\mbox{for all $z \in \mathscr{A}_r$},  \\
\mbox{and for each $r'<r$ there is a $z' \in \mathscr{A}_{r'}$ such that} \; \\
\mathrm{freq}_{z'}(\xi) 
\not\in Q_{r'}(z') [1-(\varepsilon_2+\widetilde{\delta})\delta_{Q_{r'},\widetilde{W}},
1+(\varepsilon_2+\widetilde{\delta})\delta_{Q_{r'},\widetilde{W}}]
\end{array}
\right\}. 
\ee 
Note that when $|E|<\infty$, also $|V_{r,\widetilde{M}}|<\infty$. Furthermore, 
$V_{r,\widetilde{M}} \cap V_{r',\widetilde{M}} = \emptyset$ for $r \neq r'$. For 
$\xi \in V_{r,\widetilde{M}}$, we have 
\be{}
\Big| \Big\{ 1 \le i \le \widetilde{M}-M_r+1\colon\, 
\big(\xi^{(i)},\xi^{(i+1)},\dots,\xi^{(i+M_r-1)}\big) 
\in \mathscr{A}_{r} \Big\} \Big| 
\geq \widetilde{M}(1-2\varepsilon_2), 
\ee
in particular, there are at least $K_r := \lfloor \widetilde{M} (1-3\varepsilon_2)/M_r \rfloor$ 
elements $z_1,\dots, z_{K_r} \in \mathscr{A}_{r}$ (not necessarily distinct) appearing in this 
order as disjoint subwords of $\xi$. The $z_k$'s can for example be constructed in a ``greedy'' 
way, parsing $\xi$ from left to right as in Section~\ref{S3.2} (see, in particular, (\ref{itis})). 
This implies, in particular, that 
\be{eq:weightgoodloops}
\begin{aligned} 
\prod_{i=1}^{\widetilde{M}} \rho(|\xi^{(i)}|) 
&\leq \prod_{k=1}^{K_r}\,\,\prod_{w \; \mbox{\footnotesize in} \; z_k} \rho(|w|) 
\leq \Big( \exp\big[ (1-\varepsilon_2) 
\widetilde{M}_r \E_{Q_r}[\log \rho(\tau_1)] \big] \Big)^{K_r} \\
&\leq \exp\Big[ (1-4\varepsilon_2) 
\widetilde{M} \E_{Q_r}[\log \rho(\tau_1)] \Big] 
\leq \exp\Big[ \widetilde{M} \E_{Q_r}[\log \rho(\tau_1)] + 
\widetilde{M} \varepsilon c'_\rho \Big]
\end{aligned}
\ee 
if $\widetilde{M}$ is large enough, where 
$c'_\rho:=\sup_{k \in \mathrm{supp}(\rho)} \{-\log(\rho(k))/k\} \, (< \infty)$ and we use that 
$\varepsilon_2 m^* \le \varepsilon$ by definition. 
Furthermore, for each $r \in \{1,\dots, R\}$ and $\eta \in 
V_{r,\widetilde{M}}$, we have 
\be{eq:recutVr}
\big| \big\{ \zeta \in V_{r,\widetilde{M}} \colon\, 
\kappa(\zeta)=\kappa(\eta) \big\} \big| 
\leq \exp\Big[ \widetilde{M} (H_{\tau|K}(Q_r)+\delta_1) \Big],  
\ee
where $\delta_1$ can be made arbitrarily small by choosing $\varepsilon$ small. (Note that the 
quantity on the left-hand side is the number of ways in which $\kappa(\eta)$ can be ``re-cut'' 
to obtain another element of $V_{r,\widetilde{M}}$.) In order to check (\ref{eq:recutVr}), we 
note that any $\zeta \in V_{r,\widetilde{M}}$ must contain at least $K_r$ disjoint subsentences 
from $\mathscr{A}_r$, and each $z \in \mathscr{A}_r \subset \widetilde{E}^{M_r}$ satisfies
$|\kappa(z)| \geq L$. Hence there are at most 
\be{}
{\widetilde{M}(m_{Q_r}+\varepsilon_2) - K_r (L -1) 
  \choose K_r } \le 2^{4\varepsilon_2 \widetilde{M} m_{{Q}_r}} 
\le 2^{4\varepsilon_2 m^* \widetilde{M}}
\ee
choices for the positions in the letter sequence $\kappa(\eta)$ where the concatenations of 
the disjoint subsentences from $\mathscr{A}_r$ can begin, and there are at most
\be{}
{ \widetilde{M} - K_r (M_r-1) \choose K_r } 
\le 2^{3\varepsilon_2 \widetilde{M}}
\ee
choices for the positions in the word sequence $\zeta$ where the subsentences from $\mathscr{A}_r$ 
can begin. By construction (recall the last line of (\ref{eq:Ptyp2})), each $z \in \mathscr{A}_r$ 
can be ``re-cut'' in not more than $\exp[(L/m_{Q_r})(H_{\tau|K}(Q_r) + \varepsilon_2)]$ many ways. 
Combining these observations with the fact that 
\be{}
\Big( \exp\big[ (L/m_{Q_r}) (H_{\tau|K}(Q_r) + \varepsilon_2) \big] \Big)^{K_r}
\leq \exp\Big[ \frac{\widetilde{M}}{M_r} M_r (H_{\tau|K}(Q_r) + \varepsilon_2) \Big], 
\ee
we get (\ref{eq:recutVr}) with $\delta_1 := \varepsilon_2 + 3\varepsilon_2\log 2 
+ 4\varepsilon_2 m^* \log 2$. 

We see from (\ref{eq:Pfreq.right}--\ref{eq:Pmwl.right}) and the definitions of $\widetilde{\mathcal{U}}_r$
and $V_{r,\widetilde{M}}$ that, for any $\varepsilon'>0$
\be{} 
\bigcup_{\widetilde{M}\in\N} 
\Big\{ P \in \widetilde{\mathcal{U}}_r\colon\, P(V_{r,\widetilde{M}}) 
> 1 - \varepsilon' \Big\} = \widetilde{\mathcal{U}}_r.
\ee
Put $\varepsilon_3 := \varepsilon_2 \min_{r=1,\dots,R} W_Q(\widetilde{\mathcal{U}}_r) 
\, (\leq \varepsilon_2)$. 
We can choose $\widetilde{M}$ so large that 
\be{}
W_Q\Big(\Big\{ P \in \widetilde{\mathcal{U}}_r\colon\, P(V_{r,\widetilde{M}}) 
> 1 - \frac{\varepsilon_3}{4} \Big\}\Big) 
> W_Q(\widetilde{\mathcal{U}}_r)\left(1-\frac{\varepsilon_2}{2}\right), 
\quad r=1,\dots,R.
\ee

For $M' > \widetilde{M}$ and $r=1,\dots,R$, put 
\be{eq:defWrm1}
W_{r,M'} := 
\big\{ \zeta \in \widetilde{E}^{M'}\colon\,
\mathrm{freq}_{V_{r,\widetilde{M}}}(\zeta) > 1-\varepsilon_3/2 \big\}. 
\ee
Note that for $r \neq r'$ (because $V_{r,\widetilde{M}} \cap V_{r',\widetilde{M}} = \emptyset$) 
there cannot be much overlap between $\zeta \in W_{r,M'}$ and $\eta \in W_{r',M'}$: 
\be{eq:overlapcontrol1}
\max\{ k\colon\, \mbox{$k$-suffix of $\zeta$}=\mbox{$k$-prefix of $\eta$} \}
\le \varepsilon_3 M' 
\ee
(here, the $k$-prefix of $\eta \in \widetilde{E}^n$, $k<n$, consists of the first $k$ words, 
the $k$-suffix of the last $k$ words). To see this, note that any subsequence of length $k$ 
of $\zeta$ must contain at least $(k-\varepsilon_3 M'/2)_+$ positions where a sentence from 
$V_{r,\widetilde{M}}$ starts, and any subsequence of length $k$ of $\eta$ must contain at 
least $(k-\varepsilon_3 M'/2)_+$ positions where a sentence from $V_{r',\widetilde{M}}$ starts, 
so any $k$ appearing in (\ref{eq:overlapcontrol1}) must satisfy $2(k-\varepsilon_3 M'/2)_+ \le k$, 
which enforces $k \le \varepsilon_3 M'$. 

Observe that 
(\ref{eq:defWrm1}) implies that we may choose $M'$ so large that for $r=1,\dots,R$, 
\be{eq:WrMsubwords}
\mbox{each $\zeta \in W_{r,M'}$ contains at least} \; 
(1-\varepsilon_3)\frac{M'}{\widetilde{M}} \;
\mbox{disjoint subsentences from $V_{r,\widetilde{M}}$}.
\ee

For $P \in \mathcal{P}^{\mathrm{erg,fin}}(\widetilde{E}^\N)$ with 
$P(V_{r,\widetilde{M}})>1-\varepsilon_3/3$ we have 
\be{}
\lim_{M'\to\infty} P(W_{r,M'}) = 1, 
\ee
and hence 
\be{}
\bigcup_{M'>\widetilde{M}} 
\Big\{ P \in \widetilde{\mathcal{U}}_r \colon\, P(W_{r,M'}) > 1-\varepsilon_2 \Big\} 
\supset \Big\{ P \in \widetilde{\mathcal{U}}_r\colon\, P(V_{r,\widetilde{M}}) 
> 1 - \varepsilon_3/3 \Big\},
\ee
and so we can choose $M'$ so large that 
\be{eq:Mprimechoice}
W_Q\Big(\big\{ P \in \widetilde{\mathcal{U}}_r\colon\, P(W_{r,M'}) 
> 1 - \varepsilon_2 \big\}\Big) 
> W_Q(\widetilde{\mathcal{U}}_r)(1-\varepsilon_2),
\quad r=1,\dots,R.
\ee

Now define 
\be{}
\mathcal{O}(Q) := 
\Big\{ Q' \in \mathcal{P}^{\mathrm{inv}}(\widetilde{E}^\N)\colon\,
Q'(W_{r,M'}) > W_Q(\widetilde{\mathcal{U}}_r)(1-2\varepsilon_2), \; 
r=1,\dots,R \Big\}. 
\ee
Note that $\mathcal{O}(Q)$ is open in the weak topology on $\mathcal{P}^{\mathrm{inv}}
(\widetilde{E}^\N)$, since it is defined in terms of requirements on certain finite marginals 
of $Q'$, and that for $r=1,\dots,R$, 
\be{} 
Q(W_{r,M'}) = 
\int_{\mathcal{P}^{\mathrm{erg}}(\widetilde{E}^\N)} Q'(W_{r,M'}) 
\, W_Q(dQ') \ge 
\int_{\widetilde{\mathcal{U}}_r} Q'(W_{r,M'}) \, W_Q(dQ') 
\ge \big(1-\varepsilon_2\big)^2 W_Q(\widetilde{\mathcal{U}}_r) 
\ee
by (\ref{eq:Mprimechoice}), so that in fact $Q \in \mathcal{O}(Q)$.

\subsubsection{Estimating the large deviation probability: good loops and filling loops}

Consider a choice of ``cut-points'' $j_1 < \cdots < j_N$ as appearing in the sum in (\ref{RNform}). 
Note that, by the definition of $\mathcal{O}(Q)$ (recall (\ref{xidef}--\ref{RNaltdef})),
\be{}
R^N_{j_1,\dots,j_N}(X) \in \mathcal{O}(Q)  
\ee
enforces  
\be{eq:enoughWMp}
\big| \big\{ 1 \le i \le N - M'\colon\,
(X|_{(j_{i-1},j_i]},\dots,X|_{(j_{i+M'-1},j_{i+M'}]}) \in 
W_{r,M'} \big\} \big| 
\ge N W_Q(\widetilde{\mathcal{U}}_r) (1-3 \varepsilon_2), 
\quad r=1,\dots,R, 
\ee
when $N$ is large enough. This fact, together with (\ref{eq:overlapcontrol1}), enables us to 
pick at least 
\be{}
J:=\sum_{r=1}^R \lceil (1-4\varepsilon_2)N/M' \rceil W_Q(\widetilde{\mathcal{U}}_r)
\ee 
subsentences $\zeta_1,\dots, \zeta_J$ occurring as disjoint subsentences in this order on 
$\xi_N$ such that 
\be{eq:WrMpfraction}
\big|\big\{ 1 \le j \leq J \colon\, \zeta_j \in W_{r,M'} \big\}\big| > 
(1-4\varepsilon_2) W_Q(\widetilde{\mathcal{U}}_r) \, \frac{N}{M'}, 
\quad r=1,\dots,R, 
\ee
where we note that $J \ge (1-4\varepsilon_2)(1-2\varepsilon)(N/M') 
\, (\geq (1-8\varepsilon)(N/M'))$ by (\ref{eq:WQUtildeslarge}). Indeed, we can 
for example construct these $\zeta_j$'s iteratively in a ``greedy'' way, parsing through $\xi_N$ 
from left to right and always picking the next possible subsentence from one of the $R$ 
types whose count does not yet exceed $(1-4\varepsilon_2) W_Q(\widetilde{\mathcal{U}}_r)\,(N/M')$,
as follows. Let $k_{s,r}$ be total number of subsentences of type $r$ we have chosen after the 
$s$-th step ($k_{0,1}=\cdots=k_{0,R}=0$). If in the $s$-th step we have picked $\zeta_s=(\xi_N^{(p)},
\dots,\xi_N^{(p+M'-1)})$ at position $p$, then let 
\be{}
p' := \min\big\{ i \ge p+M' \colon\, \mbox{at position $i$ in $\xi_N$ 
starts a sentence from $W_{u,M'}$ for some $u\in U_s$}\big\}, 
\ee
where $U_s:=\{r\colon\, k_{r,s} < (1-4\varepsilon_2) W_Q(\widetilde{\mathcal{U}}_r)\,(N/M')\}$, 
pick the next subsentence $\zeta_{s+1}$ starting at position $p'$ (say, of type $u$) and 
increase the corresponding $k_{s+1,u}$. Repeat this until $k_{s,r} \ge (1-4\varepsilon_2) 
W_Q(\widetilde{\mathcal{U}}_r) \, (N/M')$ for $r=1,\dots,R$. 

In order to verify that this algorithm does not get stuck, let $\mathrm{rem}(s,r)$ be the 
``remaining'' number of positions (to the right of the position where the word was picked 
in the $s$-th step) where a subsentence from $W_{r,M'}$ begins on $\xi_N$. By (\ref{eq:enoughWMp}), 
we have 
\be{}
\mathrm{rem}(0,r) \ge N W_Q(\widetilde{\mathcal{U}}_r) (1-3 \varepsilon_2).
\ee
If in the $s$-th step a subsentence of type $r$ is picked, then we have $\mathrm{rem}(s+1,r) \geq 
\mathrm{rem}(s,r)-M'$, and for $r'\neq r$ we have $\mathrm{rem}(s+1,r') \geq \mathrm{rem}(s,r')
-\varepsilon_3 M'$ by (\ref{eq:overlapcontrol1}). Thus, 
\be{}
\begin{aligned}
\mathrm{rem}(s,r) 
& \ge \mathrm{rem}(0,r) - k_{s,r} M' - (s-k_{s,r})\varepsilon_3 M' \\
& = \mathrm{rem}(0,r) -k_{s,r}(1-\varepsilon_3) M' - s \varepsilon_3 M', 
\end{aligned}
\ee
which is $>0$ as long as $k_{s,r} < (1-4\varepsilon_2) W_Q(\widetilde{\mathcal{U}}_r) \, (N/M')$
and $s<J$. 

\medskip\noindent
{\bf A.\ Combinatorial consequences.}
By (\ref{eq:WrMsubwords}) and (\ref{eq:WrMpfraction}), $R^N_{j_1,\dots,j_N}(X) \in \mathcal{O}(Q)$ 
implies that $\xi_N$ contains at least 
\be{}
C:=\sum_{r=1}^R \big\lceil (1-4\varepsilon_2) W_Q(\widetilde{\mathcal{U}}_r) \frac{N}{M'} 
\big\rceil \big\lceil (1-\varepsilon_2) \frac{M'}{\widetilde{M}} \big\rceil  
\quad \Big( \ge (1-5\varepsilon_2)(1-2\varepsilon) \frac{N}{\widetilde{M}} \: \Big)
\ee
disjoint subsentences $\eta_1,\dots,\eta_C$ (appearing in this order in $\xi_N$) such that at least 
\be{eq:respectfreq}
\frac{N}{\widetilde{M}} (1-6\varepsilon_2) W_Q(\widetilde{\mathcal{U}}_r) \; 
\mbox{of the $\eta_c$'s are from $V_{r,\widetilde{M}}$}, \;\; r=1,\dots,R.
\ee
Let $k_1,\dots,k_C$ ($k_{c+1}\ge k_c+\widetilde{M}$, $1\leq c < C$) be the indices where the 
disjoint subsentences $\eta_c$ start in $\xi_N$, i.e.,
\be{}
\eta_c=\Big( \xi_N^{(k_c)}, \xi_N^{(k_c+1)},\dots, 
\xi_N^{(k_c+\widetilde{M}-1)} \Big) 
\in V_{r_c,\widetilde{M}}, \quad i=c,\dots,C, 
\ee
and the $r_c$'s must respect the frequencies dictated by the $W_Q(\widetilde{\mathcal{U}}_r)$'s 
as in (\ref{eq:respectfreq}). Thus, each choice $(j_1,\dots,j_N)$ yielding a non-zero 
summand in (\ref{RNform}) leads to a triple 
\be{}
(\ell_1,\dots,\ell_C), (r_1,\dots,r_C), 
(\overline{\eta}_1, \dots, \overline{\eta}_C) 
\ee
such that $\overline{\eta}_c \in \kappa(V_{r_c,\widetilde{M}})$, $\ell_{c+1} \ge \ell_c 
+ |\overline{\eta}_c|$, the $r_c$'s respect the frequencies as in (\ref{eq:respectfreq}), 
and 
\be{}
\mbox{the word $\overline{\eta}_c$ starts at position $\ell_c$ in $X$ 
for $c=1,\dots,C$.} 
\ee
As in Section~\ref{S3.3}, we call such triples {\em good}, the loops inside the subsentences 
$\eta_i$ {\em good}\/ loops, the others {\em filling}\/ loops. 
\smallskip 

Fix a good triple for the moment. In order to count how many choices of $j_1 < \cdots < j_N$ 
can lead to this particular triple and to estimate their contribution, observe the following: 
\begin{enumerate}
\item 
There are at most 
\be{eq:choices.ki}
{N - C(\widetilde{M}-1) \choose C } \le \exp(\delta_1' N) 
\ee
choices for the $k_1 < \dots < k_C$, where $\delta_1'$ can be made 
arbitrarily small by choosing $\varepsilon$ small and $\widetilde{M}$ large. 
\item 
Once the $k_c$'s are fixed, by (\ref{eq:recutVr}) and (\ref{eq:respectfreq}) 
there are at most 
\be{}
\begin{aligned}
\prod_{r=1}^R &
\Big(\exp\big[ \widetilde{M} (H_{\tau|K}(Q_r)+\delta_1) 
\big]\Big)^{\frac{N}{\widetilde{M}} W_Q(\widetilde{\mathcal{U}}_r)} \\
& = \exp\Big[ N
\sum_{r=1}^R W_Q(\widetilde{\mathcal{U}}_r) (H_{\tau|K}(Q_r)+\delta_1)\Big] 
\end{aligned}
\ee 
choices for the good loops and, by (\ref{eq:weightgoodloops}), for each choice of the good loops 
the product of the $\rho(j_k-j_{k-1})$'s inside the good loops is at most 
\be{}
\begin{aligned} 
\prod_{r=1}^R &
\Big(\exp\big[ 
\widetilde{M} \E_{Q_r}[\log \rho(\tau_1)] \big] + \widetilde{M} c'_\rho \varepsilon
\Big)^{\frac{N}{\widetilde{M}} W_Q(\widetilde{\mathcal{U}}_r)} \\ 
& \leq \exp\Big[ 
N c'_\rho \varepsilon + 
N\sum_{r=1}^R W_Q(\widetilde{\mathcal{U}}_r) \E_{Q_r}[\log \rho(\tau_1)]
\Big].
\end{aligned}
\ee 
\item 
For each choice of the $k_c$'s, the contribution of the filling loops to the weight is 
\begin{eqnarray} 
\lefteqn{\rho^{*(k_1-1)}(\ell_1-1) 
\prod_{c=1}^{C-1} \rho^{*(k_{c+1}-k_c-\widetilde{M})}
(\ell_{c+1}-\ell_c-|\overline{\eta}_c|) } \nonumber \\
& \leq & 
(C_\rho \vee 1)^C k_1^{\alpha+1} 
\prod_{c=1}^{C-1} (k_{c+1}-k_c-\widetilde{M})^{\alpha+1} 
\prod_{c=1}^C \big( (\ell_{c}-\ell_{c-1}-|\overline{\eta}_{c-1}|)
\vee 1 \big)^{-\alpha} \nonumber \\
& \leq & 
(C_\rho \vee 1)^C 
\Big( \frac{N- C \widetilde{M}}{C} \Big)^{(\alpha+1)C} 
\prod_{c=1}^C \big( (\ell_{c}-\ell_{c-1}-|\overline{\eta}_{c-1}|)
\vee 1 \big)^{-\alpha} \nonumber \\ 
\label{eq:weight.fill.loops}
& \leq & e^{\delta_2' N} 
\prod_{c=1}^C \big( (\ell_{c}-\ell_{c-1}-|\overline{\eta}_{c-1}|)
\vee 1 \big)^{-\alpha}, 
\end{eqnarray}
where $\delta_2'$ can be made arbitrarily small by choosing $\varepsilon$ small and 
$\widetilde{M}$ large (and we interpret $\ell_0=0$, $|\overline{\eta}_{0}|=0$). Here, 
we have used Lemma~\ref{lemmarhoconv} in the first inequality, as well as the fact that 
the product $\prod_{c=1}^{C-1} (k_{c+1}-k_c-\widetilde{M})$ is maximal when all factors 
are equal in the second inequality. 
\end{enumerate}

Combining (\ref{eq:choices.ki}--\ref{eq:weight.fill.loops}), we see that 
\be{eq:genub1}
\begin{aligned}
\Pr\big(R_N & \in \mathcal{O}(Q) \big| X \big) \\ 
& \leq 
e^{(\delta_1'+\delta_2'+\delta_1+\varepsilon c'_\rho)N} 
\exp\Big[ N 
\sum_{r=1}^R W_Q(\widetilde{\mathcal{U}}_r) 
\big( H_{\tau|K}(Q_r) + \E_{Q_r}[\log \rho(\tau_1)] \big)\Big] \\
& \hspace{7em} 
{} \times \sum_{(\ell_i), (r_i), (\overline{\eta}_i) \atop \mbox{\tiny good}}
\prod_{i=1}^C \big( (\ell_{i}-\ell_{i-1}-|\overline{\eta}_{i-1}|)\vee 1 \big)^{-\alpha}.
\end{aligned}
\ee 

We claim that $X$-a.s. 
\be{eq:genub2} 
\begin{aligned}
\limsup_{N\to\infty} \frac{1}{N} \log 
\sum_{(\ell_i), (r_i), (\overline{\eta}_i) 
\atop \mbox{\tiny good}} 
&\prod_{i=1}^C \big( (\ell_{i}-\ell_{i-1}-|\overline{\eta}_{i-1}|)
\vee 1 \big)^{-\alpha} \\
& \leq \delta_2   
-\alpha 
\sum_{r=1}^R W_Q(\widetilde{\mathcal{U}}_r) 
m_{Q_r} H\big(\Psi_{Q_r} \mid \nu^{\otimes \N}\big), 
\end{aligned}
\ee
where $\delta_2$ can be made arbitrarily small by choosing $\varepsilon$ small and $L$ large. 
A proof of this is given below. Observe next that (\ref{eq:genub1}--\ref{eq:genub2}) (recall 
also (\ref{eqnsre1})) yield that $X$-a.s.\ 
(with $\delta:=\delta_1'+\delta_2'+\delta_1+\delta_2+\varepsilon c'_\rho$)
\be{eq:genub3} 
\begin{aligned} 
&\limsup_{N\to\infty} \frac{1}{N} \log 
\Pr\big(R_N \in \mathcal{O}(Q) \big| X \big)\\
&\leq \delta - 
\sum_{r=1}^R W_Q(\widetilde{\mathcal{U}}_r) 
\Big( H\big(Q_r\mid q_{\rho,\nu}^{\otimes \N}\big) 
+ (\alpha-1) m_{Q_r} H\big(\Psi_{Q_r} \mid \nu^{\otimes \N}\big)\Big) \\
&\leq \delta +  
2\alpha\varepsilon 
-\frac1{1+\varepsilon}
\int_{\mathcal{P}^{\mathrm{erg}}(\widetilde{E}^N)} 
H\big(Q'\mid q_{\rho,\nu}^{\otimes \N}\big) 
+ (\alpha-1) m_{Q'} H\big(\Psi_{Q'} \mid \nu^{\otimes \N}\big) W_Q(dQ') \\
&= - 
\frac1{1+\varepsilon}I^{\rm fin}(Q) + \delta + 2\alpha\varepsilon 
\end{aligned}
\ee
(use (\ref{eq:Ifinappr1}) for the second inequality, and see (\ref{tildeIQdecomp}) for the 
last equality), which completes the proof. 

\medskip\noindent
{\bf B.\ Coarse-graining $X$ with $R$ colours.} 
It remains to verify (\ref{eq:genub2}), for which we employ a coarse-graining scheme similar 
to the one used in Section~\ref{S3.4} (with block lengths $\lceil (1-\varepsilon_2)L \rceil$,
etc.) To ease notation, we silently replace $L$ by $(1-\varepsilon_2)L$ in the following. 
Split $X$ into blocks of $L$ consecutive letters, define a $\{0,1\}$-valued array $A_{i,r}$, 
$i\in\N$, $r\in\{1,\dots,R\}$ as in Section~\ref{S3.4} inductively: For each $r$, put $A_{0,r}:=0$ 
and, given that $A_{0,r}, A_{1,r}, \dots, A_{l-1,r}$ have been assigned values, define $A_l$ as 
follows:
\begin{itemize}
\item[(1)]
If $A_{l-1,r}=0$, then 
\be{defAs21mult}
A_{l,r} := \left\{
\begin{array}{ll} 
1,  & \parbox{30em}{if in $X$ there is a word from $\kappa(\mathscr{A}_r)$ 
starting in $((l-1)L, lL]$,}\\[0.5ex]
0, & \mbox{otherwise}. 
\end{array}\right.
\ee
\item[(2)]
If $A_{l-1,r}=1$, then  
\be{defAs22mult}
A_l := \left\{
\begin{array}{ll} 
1, & \parbox[t]{30em}{if in $X$ there are two words from $\kappa(\mathscr{A}_r)$
starting in $((l-2)L,(l-1)L]$, respectively, $((l-1)L,lL]$ 
and occurring disjointly,}\\[4ex]
0, \;\; & \mbox{otherwise}. 
\end{array}
\right.
\ee
\end{itemize}
Put 
\be{}
p_r := L \exp\Big( -(1-2\varepsilon_2) L H(\Psi_{Q_r} \mid \nu^{\otimes\N})
             \Big).
\ee
Arguing as in Section~\ref{S3.4}, we can couple the $(A_{i,r})_{i\in\N,1\le r \leq R}$ 
with an array $\overline{\omega}=(\omega_{i,r})_{i\in\N,1\le r \le R}$ such that $A_{i,r} 
\leq \omega_{i,r}$ and the sequence $\big((\omega_{i,1}, \dots,\omega_{i,R})\big)_{i\in\N}$ 
is i.i.d.\ with $\Pr(\omega_{i,r}=1)=p_r$. In particular, for each $r$, $(\omega_{i,r})_{i\in\N}$ 
is a Bernoulli($p_r$)-sequence. There may (and certainly will be if $\Psi_{Q_r}$ and $\Psi_{Q_{r'}}$ 
are similar) an arbitrary dependence between the $\omega_{i,1},\dots,\omega_{i,R}$ for fixed $i$, 
but this will be harmless in the low-density limit we are interested in.

For $r \in \{1,\dots,R\}$, put $d_r:=W_Q(\widetilde{\mathcal{U}}_r)(1-6\varepsilon_2)$, $D_r := 
\lceil (1-\varepsilon_2) \widetilde{M} m_{Q_r}/L \rceil$. If $\eta_c \in V_{r_c,\widetilde{M}}$, 
then 
\be{}
|\kappa(\eta_c)| \in \widetilde{M} m_{Q_{r_c}} (1-\varepsilon_2,1+\varepsilon_2),
\ee 
so $\kappa(\eta_c)$ covers at least $D_{r_c}$ consecutive $L$-blocks of the coarse-graining. 
Furthermore, as $\eta_c$ in turn contains at least $D_{r_c}(1-3\varepsilon_2)$ disjoint subsentences 
from $\mathscr{A}_{r_c}$, we see that at least $D_{r_c}(1-3\varepsilon_2)$ of these blocks must have 
$A_{k,r_c}=1$. Thus, for fixed $X$, we read off from each good triple $(\ell_c), (r_c), 
(\overline{\eta}_c)$ numbers $m_1 < \cdots < m_C$ such that 
\be{eq:condubcoarsegrained}
\begin{aligned}
& m_{c+1} \ge m_c + D_{r_c},\, c=1,\dots,C-1,  \\
& \big|\big\{ m_c \le k < m_c+D_{r_c}\colon\, A_{k,r_c}=1 \}\big| 
  \ge D_{r_c}(1-3\varepsilon_2),\, c=1,\dots,C,\\
& \big|\big\{ 1 \le c \le C \colon\, r_c=r \}\big| \ge d_r C, \, r=1,\dots,R. 
\end{aligned}
\ee
where $m_c$ is the 
index of the $L$-block that contains $\ell_c$. Furthermore, note that for 
a given ``coarse-graining'' $(m_c)$ and $(r_c)$ satisfying (\ref{eq:condubcoarsegrained}), there 
are at most 
\be{eq:cg.overcountcontr}
L^C \Big(2\varepsilon_2 \widetilde{M} \, \max_{r=1,\dots,R} m_{Q_r} 
\Big)^C \le \exp(\delta_3 N) 
\ee
choices for $\ell_c$ and $\overline{\eta}_c$ that lead to a good triple $(\ell_c), (r_c), 
(\overline{\eta}_c)$ with this particular coarse-graining. Indeed, for each $c=1,\dots,C$ 
there are at most $L$ choices for $\ell_c$ and, since each $\eta \in V_{r_c,\widetilde{M}}$ 
satisfies 
\be{}
|\kappa(\eta)| \in \widetilde{M} m_{Q_{r_c}} (1-\varepsilon_2,1+\varepsilon_2),
\ee 
there are at most $2\varepsilon_2\widetilde{M} m_{Q_{r_c}}$ choices for $\overline{\eta}_c$ (note that 
once $\ell_c$ is fixed as a ``starting point'' for a word on $X$, choosing $\overline{\eta}_c$ 
in fact amounts to choosing an ``endpoint''). Note that $\delta_3$ can be made arbitrarily small 
by choosing $\varepsilon$ small and $\widetilde{M}$ large. Finally, (\ref{eq:cg.overcountcontr}) 
and Lemma~\ref{corelemma.plusplus} yield (\ref{eq:genub2}). Indeed, since 
\begin{eqnarray}
\limsup_{N\to\infty} \frac{C}{N} & \leq & \frac{1}{\widetilde{M}}, \\ 
\nonumber 
\sum_{r=1}^R d_r D_r \log p_r & \leq & 
- (1-8\varepsilon_2) \sum_{r=1}^R W_Q(\widetilde{U}_r) \frac{\widetilde{M} m_{Q_r}}{L}
\Big( L H(\Psi_{Q_r} \mid \nu^{\otimes\N}) - \log L \Big) \\
& \leq & 
- \widetilde{M} \sum_{r=1}^R W_Q(\widetilde{U}_r) m_{Q_r} H(\Psi_{Q_r} \mid \nu^{\otimes\N}) 
+ \Big( 8\varepsilon_2 m^* K_1 + \frac{\log L}{L} \Big) \widetilde{M} , 
\end{eqnarray}
by choosing $\varepsilon$ small (note that $\varepsilon_2 m^* K_1 \leq \varepsilon$), 
$L$ and $\widetilde{M}$ large, and $\gamma$ sufficiently 
close to $1/\alpha$, the right-hand side of (\ref{eq:clpp.statement}) is smaller than the 
right-hand side of (\ref{eq:genub2}).

\subsubsection{A multicolour version of the core lemma}

The following is an extension of Lemma~\ref{lowtoylemma}. Let $R\in\N$, $\overline{\omega}_i
=(\omega_{i,1},\dots,\omega_{i,R})\in\{0,1\}^R$, and assume that $(\overline{\omega}_i)_{i\in\N}$ 
is i.i.d.\ with
\be{eq:omega.marginal.old} 
\Pr(\omega_{i,r}=1)=p_r, \quad i\in \N, \: r=1,\dots,R. 
\ee 
Note that there may be an arbitrary dependence between the $\omega_{i,r}$'s for fixed $i$. This 
will be harmless in the limit we are interested in below.

\begin{lemma}
\label{corelemma.plusplus} 
Let $\alpha\in (1,\infty)$, $\varepsilon >0$, $(d_1,\dots,d_R) \in [0,1]^R$ with
$\sum_{r=1}^R d_r \le 1$, $D_1,\dots,D_R \in \N$,
$C \in \N$, put 
\be{} 
S_C(\overline{\omega}) := {\sum}^*_{m_1,\dots,m_C \atop r_1,\dots,r_C} 
\prod_{i=1}^C \big(m_i-m_{i-1}-D_{r_{i-1}}\big)^{-\alpha}, 
\ee
where the sum ${\sum}^*$ extends over all pairs of $C$-tuples $m_0:=0 < m_1 < \cdots < m_C$ 
from $\N^C$ and $(r_1,\dots,r_C) \in \{1,\dots,R\}^C$ satisfying the constraints 
\be{eq:clpp.constr}
\begin{aligned} 
&m_{i+1} \geq m_i + D_{r_i}, \\
& |\{ 1 \le i \le C \colon\, r_i=r\}| \geq d_r C, \,\, r=1,\dots,R, \\
& |\{ m_i \le k < m_i+D_{r_i} \colon\, \omega_{k,r_i}=1 \}| 
\ge D_{r_i}(1-\varepsilon), \,\, i=1,\dots,C. 
\end{aligned} 
\ee
Then $\overline{\omega}$-a.s.
\be{eq:clpp.statement}
\begin{aligned}
&\limsup_{C\to\infty} \frac{1}{C} \log S_C(\overline{\omega}) \\
& \hspace{1em}
\leq \inf_{\gamma \in (1/a,1)} 
\Big\{ \frac{1}{\gamma}\Big( \log \zeta(a\gamma) + h(\underline{d})  
+ d_0 \log R  + \big(\log 2 \big)\sum_{r=1}^R d_r D_r 
+ (1-\varepsilon) \sum_{r=1}^R d_r D_r \log p_r \Big) \Big\}, 
\end{aligned}
\ee
where $h(\underline{d}):=-\sum_{r=0}^R d_r \log d_r$ (with $d_0:=1-d_1-\cdots-d_R$) is the 
entropy of $\underline{d}$. 
\end{lemma}

\begin{proof}
The proof is a variation on the proof of Lemma~\ref{lowtoylemma}. We again estimate fractional 
moments. For $\gamma \in (1/\alpha,1)$, we have 
\be{eq:clpp.fracmom1}
\begin{aligned}
&\E\big[ (S_C)^\gamma \big]\\ 
&\leq  {\sum}^{'}_{r_1,\dots,r_C} 
\sum_{ {m_1,\dots,m_C} \atop {m_{i+1} \geq m_i + D_{r_i}\,\forall\,i} } 
\Pr\Big( \cap_{i=1}^C \big\{ |\{ k \in [m_i,m_i+D_{r_i}-1] \colon\, 
\omega_{k,r_i}=1 \}| \geq (1-\varepsilon) D_{r_i} \big\} \Big) \\[-3ex]
& \hspace{20em} \times \prod_{i=1}^C (m_i-m_{i-1}-D_{r_{i-1}})^{-\alpha \gamma}, 
\end{aligned}
\ee
where the sum ${\sum}^{'}$ extends over all $(r_1,\dots,r_C)$ satisfying the constraint in 
the second line of (\ref{eq:clpp.constr}). Noting that 
\begin{eqnarray*}
\Pr\Big( |\{ k \in [m_i,m_i+D_{r_i}-1] \colon\, \omega_{k,r_i}=1 \}| 
\geq (1-\varepsilon) D_{r_i} \Big) 
& = & \sum_{m=(1-\varepsilon) D_{r_i}}^{D_{r_i}} 
{D_{r_i} \choose k} p_r^m (1-p_r)^{D_{r_i}-m} \\
&\leq & p_r^{(1-\varepsilon) D_{r_i}} 2^{D_{r_i}}
\end{eqnarray*}
and 
\begin{eqnarray*}
\lefteqn{ \big| \big\{ (r_1,\dots,r_C) \in \{1,\dots,R\}^C \colon\,
\mbox{at least $d_r C$ of the $r_i=r$}, \: r=1,\dots,R \big\} \big| } \\
& \le & R^{d_0 C} { C \choose d_0C \; d_1C \; \dots \; d_RC } 
= \exp\Big[ C\big( d_0 \log R + h(\underline{d}) + o(1) \big) \Big], 
\end{eqnarray*}
we see from (\ref{eq:clpp.fracmom1}) that 
\be{}
\begin{aligned}
\E\big[ (S_C)^\gamma \big] \leq 
&\exp\Big[ C\big( d_0 \log R + h(\underline{d}) + o(1) \big) \Big] 
\times \prod_{r=1}^R \big( 2 p_r^{(1-\varepsilon)} \big)^{d_r C D_r} \\
& \hspace{6em} {} 
\times 
\sum_{ {m_1,\dots,m_C} \atop {m_{i+1} \geq m_i + D_{r_i}\,\forall\,i} } 
\prod_{i=1}^C (m_i-m_{i-1}-D_{r_{i-1}})^{-\alpha\gamma} \\
&= \exp C\Big[ d_0 \log R + h(\underline{d}) + \log \zeta(a\gamma) 
+ {\textstyle \sum_{r=1}^R d_r D_r \log 2} 
+ {\textstyle (1-\varepsilon) \sum_{r=1}^R d_r D_r \log p_r} \Big], 
\end{aligned}
\ee
which yields (\ref{eq:clpp.statement}) as in the proof of 
Lemma~\ref{lowtoylemma}. 
\end{proof}

\subsection{Step 6: Weakening the tail assumption} 
\label{S3.6}

We finally show how to go from (\ref{eq:algtail}) to (\ref{rhocond}). Suppose that $\rho$ 
satisfies (\ref{rhocond}) with a certain $\alpha \in (1,\infty)$. Then, for any $\alpha' 
\in (1,\alpha)$, there is a $C_\rho(\alpha')$ such that (\ref{eq:algtail}) holds for
this $\alpha'$. Hence, as shown in Sections \ref{S3.1}--\ref{S3.4}, for any $\varepsilon>0$ 
we can find a neighbourhood $\mathcal{O}(Q)\subset\mathcal{P}^{\mathrm{inv,fin}}
(\widetilde{E}^\N)$ of $Q$ such that 
\be{Prubneigh}
\limsup_{N\to\infty} \frac{1}{N} 
\log \Pr\big(R_N \in \mathcal{O}(Q) \mid X\big) 
\leq  - H(Q \mid q_{\rho,\nu}^{\otimes\N}) -(\alpha'-1)\,m_{Q}\,H(\Psi_{Q} \mid \nu^{\otimes \N}) 
 + \frac{\varepsilon}{2} \qquad X-a.s.
\ee
The right-hand side is $\leq - I^{\mathrm{fin}}(Q) + \varepsilon$ for 
$\alpha'$ sufficiently close to $\alpha$, so that we again get (\ref{mainub}). 
\end{proof}


\section{Lower bound} 
\label{S4}

The following lower bound will be used in Section~\ref{S5} to derive 
the lower bound in the definition of the LDP.

\begin{proposition}
\label{proplb} 
For any $Q \in \mathcal{P}^{\mathrm{inv,fin}}(\widetilde{E}^\N)$ and any open 
neighbourhood $\mathcal{U}(Q) \subset \mathcal{P}^{\mathrm{inv}}(\widetilde{E}^\N)$ 
of $Q$, 
\be{mainlb} 
\liminf_{N\to\infty} \frac{1}{N} \log 
\Pr\big( R_N \in \mathcal{U}(Q) \mid X \big) 
\geq -I^{\mathrm{fin}}(Q) \quad X-a.s.
\ee
\end{proposition}

\begin{proof} 
Suppose first that $Q \in \mathcal{P}^{\mathrm{erg,fin}}(\widetilde{E}^\N)$. Then,
informally, our strategy runs as follows. In $X$, look for the first string of length 
$\approx N m_Q$ that looks typical for $\Psi_Q$. Make the first jump long enough so 
as to land at the start of this string. Make the remaining $N-1$ jumps typical for $Q$. 
The probability of this strategy on the exponential scale is the conditional specific 
relative entropy of word lengths under $Q$ w.r.t.\ $\rho^{\otimes\N}$ given the 
concatenation, i.e., $\approx \exp[N(H_{\tau|K}(Q)+\E_Q[\log\rho(\tau_1)])]$, times
the probability of the first long jump.  In order to find a suitable string, we have 
to skip ahead in $X$ a distance $\approx \exp[N m_Q H(\Psi_Q \mid \nu^{\otimes\N})]$. 
By (\ref{rhocond}), the probability of the first jump is therefore $\approx 
\exp[-N\alpha\,m_Q H(\Psi_Q\mid\nu^{\otimes\N})]$. In view of (\ref{eqnratefctexplicit}) 
and (\ref{eqnsre1}), this yields the claim. In the actual proof, it turns out to be 
technically simpler to employ a slightly different strategy, which has the same 
asymptotic cost, where we look not only for one contiguous piece of 
``$\Psi_Q$-typical'' letters but for a sequence of $\lceil N/M \rceil$ 
pieces, each of length $\approx M m_Q$. Then we let $N\to\infty$, followed by $M\to\infty$. 
\smallskip

More formally, we choose for $\mathcal{O}(Q)$ an open neighborhood $\mathcal{O}^\prime
\subset\mathcal{O}$ of the type introduced in Section \ref{S3.2}, and we estimate 
$\Pr(R_N \in \mathcal{O}^\prime \mid X)$ from below by using (\ref{RNaltdef}--\ref{xiprop}).

Assume first that $Q$ is ergodic. We can then assume that the neighbourhood $\mathcal{U}$
is given by 
\be{low1}
\mathcal{U} = \big\{ Q' \in \mathcal{P}^{\mathrm{inv}}(\widetilde{E}^\N)\colon\, 
(\pi_{L_u} Q')(\zeta_u) \in (a_u,b_u), \; u=1,\dots,U \big\}
\ee
for some $U \in \N$, $L_1, \dots, L_U \in \N$, $0\le a_u < b_u \le 1$ and $\zeta_u \in 
\widetilde{E}^{L_u}$, $u=1,\dots,U$. As in Section~\ref{S3.1}, by ergodicity of $Q$ we 
can find for each $\varepsilon > 0$ a sufficiently large $M \in \N$ and a set $\mathscr{A} 
=\{z_1,\dots,z_A\} \subset \widetilde{E}^M$ of ``$Q$-typical sentences'' satisfying 
(\ref{sentprop}--\ref{probsent}) (with $\varepsilon_1=\delta_1=\varepsilon$, say), 
and additionally 
\be{low2}
\frac{1}{M} \big| \{ 0 \le j \le M-L_i\colon\, 
\pi_{L_u}(\widetilde{\theta}^j z_a) = \zeta_u \}\big|
\in (a_u, b_u), \quad a=1,\dots,A, \, u=1,\dots,U.
\ee
Let $\mathscr{B}:=\kappa(\mathscr{A})$. Then from (\ref{sentprop}--\ref{probsent}) 
we have that, for each $b \in \mathscr{B}$, 
\be{low3}
|I_b| = | \{ z \in \mathscr{A}\colon\,\kappa(z)=b \} | 
\geq \exp\big[ M(H_{\tau|K}(Q)-2\varepsilon) \big],
\ee
and 
\be{eq:probfindB}
\Pr(\mbox{$X$ begins with some element of $\mathscr{B}$}) 
\geq \exp\big[ - M m_Q(H(\Psi_Q \mid \nu^{\otimes \N})+2\varepsilon) \big]. 
\ee
Let 
\be{low4}
\begin{aligned}
\sigma^{(M)}_1 
&:= \min\{i\colon\,\mbox{$\theta^i X$ begins with some element of $\mathscr{B}$}\},\\
\sigma^{(M)}_l &:=\min\{i>\sigma^{(M)}_{l-1}+M (m_Q+\varepsilon)\colon\, 
\mbox{$\theta^i X$ begins with some element of $\mathscr{B}$}\},
\quad l=2,3,\dots.
\end{aligned}
\ee 
Restricting the sum in (\ref{RNform}) over $0<j_1<\dots<j_N<\infty$ such that $j_1=\sigma^{(M)}_1$, 
$j_2-j_1, \dots, j_{M}-j_{M-1}$ are the word lengths corresponding to the $z_a$'s compatible with 
$\pi_{M m_Q} (\theta^{\tau_1} X)$, $j_{M+1}=\sigma^{(M)}_2$, etc., we see that 
\be{low5}
\frac{1}{N} \log \Pr(R_N \in \mathcal{U} \mid X) \geq 
 H_{\tau|K}(Q) + \E_Q[\log\rho(\tau_1)] - 3 \varepsilon 
- \alpha \frac{1}{N} \sum_{l=1}^{\lfloor N/M\rfloor} 
\log\big(\sigma^{(M)}_l-\sigma^{(M)}_{l-1}\big)
\ee
for $N$ sufficiently large. Hence $X$-a.s. 
\be{low6}
\begin{aligned}
\liminf_{N\to\infty} 
\frac{1}{N} \log \Pr(R_N \in \mathcal{U} \mid X) & \ge 
H_{\tau|K}(Q) + \E_Q[\log\rho(\tau_1)] - 3 \varepsilon 
- \alpha \frac{1}{M} \E[\log \sigma^{(M)}_1] \\
& \ge H_{\tau|K}(Q) + \E_Q[\log\rho(\tau_1)] - \alpha m_Q(H(\Psi_Q \mid \nu^{\otimes \N}) 
- 6 \varepsilon \\
& = -I^\mathrm{fin}(Q) - 6 \varepsilon, 
\end{aligned}
\ee 
where we have used (\ref{eq:probfindB}) in the second inequality. Now let $\varepsilon
\downarrow 0$. 

It remains to remove the restriction of ergodicity of $Q$, analogously to the proof of 
Birkner \cite{Bi08}, Proposition~2. To that end, assume that $Q \in \mathcal{P}^\mathrm{inv,fin}
(\widetilde{E}^\N)$ admits a non-trivial ergodic decomposition. Then, for each $\varepsilon > 0$, 
we can find $Q_1, \dots, Q_R \in \mathcal{P}^\mathrm{erg,fin}(\widetilde{E}^\N)$, $\lambda_1, 
\dots, \lambda_R \in (0,1)$, $\sum_{r=1}^R \lambda_r=1$ such that $\lambda_1 Q_1 + \cdots 
+ \lambda_R Q_R \in \mathcal{U}$ and 
\be{low7}
\sum_{i=1}^R \lambda_r I^\mathrm{fin}(Q_r) \le I^\mathrm{fin}(Q) + \varepsilon
\ee
(for details see Birkner \cite{Bi08}, p.~723; employ the fact that both terms in $I^\mathrm{fin}$ 
are affine). For each $r=1,\dots,R$, pick a small neighbourhood $\mathcal{U}_r$ of $Q_r$ such 
that 
\be{low8}
Q'_r \in \mathcal{U}_r, \, r=1,\dots,R \quad \Longrightarrow \quad 
\sum_{i=1}^R \lambda_r Q'_r \in \mathcal{U}. 
\ee
Using the above strategy for $Q_1$ for $\lambda_1 N$ loops, then the strategy for $Q_2$ 
for $\lambda_2 N$ loops, etc., we see that 
\be{Low9}
\liminf_{N\to\infty} \frac{1}{N} \Pr(R_N \in \mathcal{U} \mid X) \ge 
- \sum_{i=1}^R \lambda_r I^\mathrm{fin}(Q_r) - 6 \varepsilon 
\ge - I^\mathrm{fin}(Q) - 7 \varepsilon. 
\ee
\end{proof}


\section{Proof of Theorem~\ref{mainthm}}
\label{S5}

\begin{proof}
The proof comes in 3 steps. We first prove that, for each word length 
truncation level $\tr \in \N$, the family $\Pr([R_N]_\tr \in \cdot \mid X)$, 
$N\in\N$, $X$-a.s.\ satisfies an LDP on 
\be{ag29}
\mathcal{P}_\tr^{\mathrm{inv}}(\widetilde{E}^\N) 
= \big\{ Q \in \mathcal{P}^{\mathrm{inv}}(\widetilde{E}^\N)\colon\, 
Q(|Y^{(1)}| \le \tr) = 1 \big\}
\ee
(recall (\ref{trunword}--\ref{trunprob})) with a deterministic rate
function $I^\mathrm{fin}([Q]_\tr)$ (this is essentially the content
of Propositions~\ref{proplb} and \ref{propub}). Note that $[Q]_\tr=Q$
for $Q \in \mathcal{P}_\tr^{\mathrm{inv}}(\widetilde{E}^\N)$, and that
$\mathcal{P}_\tr^{\mathrm{inv}}(\widetilde{E}^\N)$ is a closed subset
of $\mathcal{P}^{\mathrm{inv}}(\widetilde{E}^\N)$, in particular, a 
Polish space under the relative topology (which is again the weak 
topology). After we have given the proof for fixed $\tr$, we let 
$\tr\to\infty$ and use a projective limit argument to complete the 
proof of Theorem~\ref{mainthm}.

\medskip\noindent
{\bf 1.} Fix a truncation level $\tr \in \N$. Propositions \ref{proplb} 
and \ref{propub} combine to yield the LDP on $\mathcal{P}_\tr^{\mathrm{inv}}
(\widetilde{E}^\N)$ in the following standard manner. Note that any $Q \in 
\mathcal{P}_\tr^{\mathrm{inv}}(\widetilde{E}^\N)$ satisfies $m_Q < \infty$. 

\medskip\noindent
{\bf 1a.} Let $\mathcal{O} \subset \mathcal{P}_\tr^{\mathrm{inv}}(\widetilde{E}^\N)$
be \emph{open}. Then, for any $Q \in \mathcal{O}$, there is an open neighbourhood
$\mathcal{O}(Q) \subset \mathcal{P}_\tr^{\mathrm{inv}}(\widetilde{E}^\N)$ of 
$Q$ such that $\mathcal{O}(Q) \subset \mathcal{O}$. The latter inclusion, 
together with Proposition~\ref{proplb}, yields 
\be{LDPlb1}
\liminf_{N\to\infty} \frac{1}{N} \log \Pr\big( [R_N]_\tr \in \mathcal{O} \mid X \big)\\
\geq - I^{\mathrm{fin}}(Q) \qquad X-a.s.
\ee
Optimising over $Q\in\mathcal{O}$, we get 
\be{LDPlb2} 
\liminf_{N\to\infty} \frac{1}{N} \log \Pr\big( [R_N]_\tr \in \mathcal{O} \mid X \big)\\
\geq - \inf_{Q \in \mathcal{O}} I^{\mathrm{fin}}(Q) \qquad X-a.s.
\ee
Here, note that, since $\mathcal{P}_\tr^{\mathrm{inv}}(\widetilde{E}^\N)$ is Polish, 
it suffices to optimise over a countable set generating the weak topology, allowing
us to transfer the $X$-a.s.\ limit from points to sets (see, e.g., Comets~\cite{Co89}, 
Section~III).

\medskip\noindent
{\bf 1b.} Let $\mathcal{K} \subset \mathcal{P}_\tr^{\mathrm{inv}}(\widetilde{E}^\N)$ be 
\emph{compact}. Then there exist $M \in \N$, $Q_1,\dots,Q_M \in \mathcal{K}$ 
and open neighbourhoods $\mathcal{O}(Q_1),\dots,\mathcal{O}(Q_M) \subset 
\mathcal{P}_\tr^{\mathrm{inv}}(\widetilde{E}^\N)$ such that $\mathcal{K} \subset 
\cup_{m=1}^M \mathcal{O}(Q_m)$. The latter inclusion, together with 
Proposition~\ref{propub}, yields
\be{LDPub1}
\limsup_{N\to\infty} \frac{1}{N} \log \Pr\big( [R_N]_\tr \in \mathcal{K} \mid X \big)
\leq - \inf_{1 \leq m \leq M} I^{\mathrm{fin}}(Q_m) + \varepsilon \qquad X-a.s.
\qquad \forall\,\varepsilon>0.
\ee
Extending the infimum to $Q \in \mathcal{K}$ and letting $\varepsilon \da 0$ afterwards,
we obtain
\be{LDPub2}
\limsup_{N\to\infty} \frac{1}{N} \log \Pr\big( [R_N]_\tr \in \mathcal{K} \mid X \big)
\leq - \inf_{Q \in \mathcal{K}} I^{\mathrm{fin}}(Q) \qquad X-a.s.
\ee
{\bf 1c.} Let $\mathcal{C} \subset \mathcal{P}_\tr^{\mathrm{inv}}(\widetilde{E}^\N)$ be 
\emph{closed}. Because $Q \mapsto H(Q \mid q_{\rho,\nu}^{\otimes\N})$ has compact level sets, for 
any $M<\infty$ the set $\mathcal{K}_M = \mathcal{C} \cap \{Q \in 
\mathcal{P}_\tr^{\mathrm{inv}}(\widetilde{E}^\N) \colon\,H(Q \mid q_{\rho,\nu}^{\otimes\N}) \leq M\}$ 
is compact. Hence, doing \emph{annealing} on $X$ and using (\ref{LDPub2}), 
we get 
\be{LDPub3}
\limsup_{N\to\infty} \frac{1}{N} \log \Pr\big( [R_N]_\tr \in \mathcal{C} \mid X \big)
\leq \max\left\{-M, - \inf_{Q \in \mathcal{K}_M} I^{\mathrm{fin}}(Q)\right\} 
\qquad X-a.s.
\ee
Extending the infimum to $Q \in \mathcal{C}$ and letting $M \to \infty$ afterwards,
we arrive at
\be{LDPub4}
\limsup_{N\to\infty} \frac{1}{N} \log \Pr\big( [R_N]_\tr \in \mathcal{C} \mid X \big)
\leq - \inf_{Q \in \mathcal{C}} I^{\mathrm{fin}}(Q) \qquad X-a.s.
\ee
Equations (\ref{LDPlb2}) and (\ref{LDPub4}) complete the proof of the 
conditional LDP for $[R_N]_\tr$.

\medskip\noindent
{\bf 2.} It remains to remove the truncation of word lengths. We know from Step~1 
that, for every $\tr \in \N$, the family $\Pr([R_N]_\tr \in \cdot \mid X)$, $N \in \N$, 
satisfies the LDP on $\mathcal{P}^{\mathrm{inv}}([\widetilde{E}]_\tr^\N)$ with rate 
function $I^{\mathrm{fin}}$. Consequently, by the Dawson-G\"artner projective limit 
theorem (see Dembo and Zeitouni~\cite{DeZe98}, Theorem~4.6.1), the family $\Pr(R_N 
\in \cdot \mid X)$, $N \in \N$, satisfies the LDP on $\mathcal{P}^{\mathrm{inv}}
(\widetilde{E}^\N)$ with rate function
\be{Iratesup}
I^\mathrm{que}(Q) = \sup_{\tr \in \N} I^{\mathrm{fin}}([Q]_\tr), 
\qquad Q \in \mathcal{P}^{\mathrm{inv}}(\widetilde{E}^\N).
\ee
The $\sup$ may be replaced by a $\limsup$ because the truncation may start at any level. 
For $Q \in \mathcal{P}^{\mathrm{inv,fin}}(\widetilde{E}^\N)$, we have $\lim_{\tr\to\infty} 
I^\mathrm{fin}([Q]_\tr) = I^\mathrm{fin}(Q)$ by Lemma~\ref{lemma:appendix}, and so we get 
the claim if we can show that $\limsup$ can be replaced by a limit, which is done in Step 3. 
Note that $I^\mathrm{que}$ \emph{inherits} from $I^{\mathrm{fin}}$ the properties qualifying 
it to be a rate function: this is part of the projective limit theorem. For $I^{\mathrm{fin}}$ 
these properties are proved in Section~\ref{S6}.

\medskip\noindent
{\bf 3.} Since $I^\mathrm{que}$ is lower semi-continuous, it is equal to its lower 
semi-continuous regularisation
\be{lscreg}
\widetilde{I}^\mathrm{que}(Q) := 
\sup_{\mathcal{O}(Q)} \inf_{Q' \in \mathcal{O}(Q)} I^\mathrm{que}(Q'), 
\ee
where the supremum runs over the open neighborhoods of $Q$. For each $\tr \in \N$, 
$[Q]_\tr \in \mathcal{P}^{\mathrm{inv,fin}}(\widetilde{E}^\N)$, while 
${\rm w}-\lim_{\tr\to\infty}[Q]_\tr = Q$. So, in particular,
\be{ag52b}
I^\mathrm{que}(Q)=\widetilde{I}^\mathrm{que}(Q) 
\leq \sup_n \inf_{\tr \ge n} I^\mathrm{fin}([Q]_\tr) 
= \liminf_{\tr\to\infty} I^\mathrm{fin}([Q]_\tr),  
\ee
implying that in fact 
\be{Ieqs}
I^\mathrm{que}(Q) = \lim_{\tr\to\infty} I^\mathrm{fin}([Q]_\tr), \quad 
Q \in \mathcal{P}^{\mathrm{inv}}(\widetilde{E}^\N). 
\ee
\epr

Lemma~\ref{trunclimcont} in Appendix~\ref{SA}, together with (\ref{Ieqs}), shows that 
$I^\mathrm{que}(Q) = I^\mathrm{fin}(Q)$ for $Q\in\mathcal{P}^{\mathrm{inv,fin}}(\widetilde{E}^\N)$, 
as claimed in the first line of (\ref{eqgndefinitionI}). 


\section{Proof of Theorem~\ref{ratefct}}
\label{S6}

\begin{proof}
The proof comes in 5 steps.

\medskip\noindent
{\bf 1.} Every $Q \in \mathcal{P}^{\mathrm{inv}}(\widetilde{E}^\N)$ can be decomposed as
\be{ergdecomp}
Q = \int_{\mathcal{P}^{\mathrm{erg}}(\widetilde{E}^\N)} Q'\,W_Q(dQ')
\ee
for some unique probability measure $W_Q$ on $\mathcal{P}^{\mathrm{erg}}
(\widetilde{E}^\N)$ (Georgii~\cite{Ge88}, Proposition 7.22). If $Q \in 
\mathcal{P}^{\mathrm{inv,fin}}(\widetilde{E}^\N)$, then $W_Q$ is concentrated 
on $\mathcal{P}^{\mathrm{erg,fin}}(\widetilde{E}^\N)$ and so, by 
(\ref{eq:defmQ}--\ref{PsiQdef}),
\be{mQPsiQdecomp} 
m_Q = \int_{\mathcal{P}^{\mathrm{erg,fin}}(\widetilde{E}^\N)} m_{Q'}\,W_Q(dQ'), 
\quad
\Psi_Q = \int_{\mathcal{P}^{\mathrm{erg,fin}}(\widetilde{E}^\N)} 
\frac{m_{Q'}}{m_Q}\,\Psi_{Q'}\,W_Q(dQ').
\ee 
Since $Q \mapsto H(Q \mid q_{\rho,\nu}^{\otimes\N})$ and $\Psi \mapsto H(\Psi \mid 
\nu^{\otimes \N})$ are affine (see e.g.\ Deuschel and Stroock \cite{DeStr89}, 
Example~4.4.41), it follows from (\ref{eqnratefctexplicit}) and 
(\ref{ergdecomp}--\ref{mQPsiQdecomp}) that 
\be{tildeIQdecomp}
I^{\mathrm{fin}}(Q) = \int_{\mathcal{P}^{\mathrm{erg,fin}}(\widetilde{E}^\N)} 
I^{\mathrm{fin}}(Q')\,W_Q(dQ').
\ee
Since $Q \mapsto W_Q$ is affine, (\ref{tildeIQdecomp}) shows that $I^{\mathrm{fin}}$ 
is affine on $\mathcal{P}^{\mathrm{inv,fin}}(\widetilde{E}^\N)$. 

\medskip\noindent
{\bf 2.} Let $(Q_n)_{n \in \N} \subset \mathcal{P}^{\mathrm{inv,fin}}(\widetilde{E}^\N)$ 
be such that ${\rm w}-\lim_{n\to\infty} Q_n = Q \in \mathcal{P}^{\mathrm{inv,fin}}
(\widetilde{E}^\N)$. By Proposition~\ref{propub}, for any $\varepsilon>0$ we can find an 
open neighbourhood $\mathcal{O}(Q)\subset\mathcal{P}^{\mathrm{inv}}(\widetilde{E}^\N)$ of 
$Q$ such that 
\be{eq:couldwebesostupid} 
\limsup_{N\to\infty} \frac{1}{N} 
\log \Pr\big(R_N \in \mathcal{O}(Q) \mid X\big) 
\leq -I^{\mathrm{fin}}(Q) + \varepsilon \quad X-a.s.
\ee
On the other hand, for $n$ large enough so that $Q_n \in \mathcal{O}(Q)$, we have from
Proposition~\ref{proplb} that
\be{eq:yeswecould}
\liminf_{N\to\infty} \frac{1}{N} \log \Pr\big(R_N \in \mathcal{O}(Q) \mid X\big)
\geq -I^{\mathrm{fin}}(Q_n) \quad X-a.s.
\ee
Combining (\ref{eq:couldwebesostupid}--\ref{eq:yeswecould}), we get that, for any 
$\varepsilon>0$,  
\be{}
\liminf_{n\to\infty} I^{\mathrm{fin}}(Q_n) \ge I^{\mathrm{fin}}(Q) - \varepsilon. 
\ee  
Now let $\varepsilon\downarrow 0$, to conclude that $I^{\mathrm{fin}}$ is lower semicontinuous 
on $\mathcal{P}^{\mathrm{inv,fin}}(\widetilde{E}^\N)$ (recall also (\ref{Ieqs})).

\medskip\noindent
{\bf 3.} From (\ref{eqnratefctexplicit}) we have
\be{lbinv}
I^{\mathrm{fin}}(Q) \geq H(Q \mid q_{\rho,\nu}^{\otimes\N}) \qquad \forall\, 
Q \in \mathcal{P}^{\mathrm{inv,fin}}(\widetilde{E}^\N)
\ee
Since $\{Q\in\mathcal{P}^{\mathrm{inv}}(\widetilde{E}^\N)\colon\, H(Q \mid 
q_{\rho,\nu}^{\otimes\N}) \leq C\}$ is compact for all $C<\infty$ (see, e.g., Dembo and 
Zeitouni~\cite{DeZe98}, Corollary~6.5.15), it follows that $I^{\mathrm{fin}}$ has compact 
level sets on $\mathcal{P}^{\mathrm{inv,fin}}(\widetilde{E}^\N)$. 

\medskip\noindent
{\bf 4.} As mentioned at the end of Section~\ref{S5}, $I^\mathrm{que}$ inherits from 
$I^{\mathrm{fin}}$ that it is lower semicontinuous and has compact level sets.  In particular, 
$I^{\mathrm{que}}$ is the lower semicontinuous extension of $I^{\mathrm{fin}}$ from
$\mathcal{P}^{\mathrm{inv,fin}}(\widetilde{E}^\N)$ to $\mathcal{P}^{\mathrm{inv}}(\widetilde{E}^\N)$. 
Moreover, since $I^{\mathrm{fin}}$ is affine on $\mathcal{P}^{\mathrm{inv,fin}}(\widetilde{E}^\N)$ 
and $I^\mathrm{que}$ arises as the truncation limit of $I^\mathrm{fin}$ (recall (\ref{ag52b})), it 
follows that $I^\mathrm{que}$ is affine on $\mathcal{P}^{\mathrm{inv}}(\widetilde{E}^\N)$.

\medskip\noindent
{\bf 5.} It is immediate from (\ref{eqgndefinitionI}--\ref{eqnratefctexplicit}) that 
$q_{\rho,\nu}^{\otimes\N}$ is the unique zero of $I^\mathrm{que}$.

\end{proof}


\section{Proof of Theorem \ref{mainthmboundarycases}}
\label{S7}

\begin{proof}
The extension is an easy generalisation of the proof given in Sections~\ref{S3}--\ref{S4}.

\medskip\noindent
(a) Assume that $\rho$ satisfies (\ref{rhocond}) with $\alpha=1$. Since the LDP upper bound 
holds by the \emph{annealed} LDP (compare (\ref{Iann}) and (\ref{eqnratefctexplicit})), it 
suffices to prove the LDP lower bound. To achieve this, we first show that for any $Q \in 
\mathcal{P}^{\mathrm{inv,fin}}(\widetilde{E}^\N)$ and $\varepsilon>0$ there exists an open 
neighbourhood $\mathcal{O}(Q)\subset\mathcal{P}^{\mathrm{inv}}(\widetilde{E}^\N)$ of $Q$ 
such that 
\be{eq:lbalpha1}
\liminf_{N\to\infty} \frac{1}{N} \log \Pr\big( R_N \in \mathcal{O}(Q) \mid X \big) 
\geq - I^\mathrm{ann}(Q) - \varepsilon 
\qquad \mbox{$X$--a.s.}
\ee 
After that, the extension from $\mathcal{P}^{\mathrm{inv,fin}}(\widetilde{E}^\N)$ to
$\mathcal{P}^{\mathrm{inv}}(\widetilde{E}^\N)$ follows the argument in Section~\ref{S5}.

In order to verify (\ref{eq:lbalpha1}), observe that, by our assumption on $\rho(\cdot)$, 
for any $\alpha'>1$ there exists a $C_{\alpha'}>0$ such that  
\be{}
\frac{\rho(n)}{n^{\alpha'}} \geq C_{\alpha'} \qquad \forall\, n \in \mathrm{supp}(\rho).
\ee
Picking $\alpha'$ so close to $1$ that $(\alpha'-1) m_Q H(\Psi_Q | \nu^{\otimes \N})<
\varepsilon/2$, we can trace through the proof of Proposition~\ref{proplb} in Section~\ref{S4} 
to construct an open neighbourhood $\mathcal{O}(Q)\subset\mathcal{P}^{\mathrm{inv}}
(\widetilde{E}^\N)$ of $Q$ satisfying
\be{} 
\begin{aligned}
&\liminf_{N\to\infty} \frac{1}{N} \log \Pr\big( R_N \in \mathcal{O}(Q) \mid X \big)\\ 
&\qquad \geq - H(Q\mid q_{\rho,\nu}^{\otimes \N}) 
- (\alpha'-1) m_Q H(\Psi_Q \mid \nu^{\otimes \N}) - \varepsilon/2 
\geq -I^{\rm ann}(Q) - \varepsilon \quad X-a.s.,
\end{aligned}
\ee
which is (\ref{eq:lbalpha1}).

\medskip\noindent 
 (b) We only give a sketch of the argument. Assume
$\alpha=\infty$ in (\ref{rhocond}). For $Q \in
\mathcal{P}^{\mathrm{inv,fin}} (\widetilde{E}^\N)$, the lower bound
(which is non-zero only when $Q \in \mathscr{R}_\nu$) follows from
Birkner~\cite{Bi08}, Proposition 2, or can alternatively be obtained
from the argument in Section~\ref{S4}. Now consider a $Q \in
\mathcal{P}^{\mathrm{inv}}(\widetilde{E}^\N)$ with $m_Q=\infty$, 
$H(Q \mid q_{\rho,\nu}^{\otimes\N}) < \infty$ and 
$\lim_{\tr\to\infty} m_{[Q]_\tr} H( \Psi_{[Q]_\tr} \mid \nu^{\otimes \N}) = 0$, 
let $\mathcal{O}(Q) \subset \mathcal{P}^{\mathrm{inv}}(\widetilde{E}^\N)$ be 
an open neighbourhood of $Q$. 
For simplicity, we assume $\mathrm{supp}(\rho)=\N$. 
Fix $\varepsilon > 0$. We can find a sequence $\delta_N \downarrow 0$ such that 
\be{eq:thm14b1}
\max \Big\{ -\frac{1}{N} \log \rho(n) \colon\, n \le \lceil N \delta_N \rceil \Big\} \le \varepsilon. 
\ee
Furthermore, 
\be{} 
\frac{1}{N} h\left(Q_{|_{\mathscr{F}_N}}
\mid q_{\rho,\nu}^{\otimes N} \right) \geq H(Q \mid q_{\rho,\nu}^{\otimes\N}) - \varepsilon
\ee
for $N \ge N_0=N_0(\varepsilon, Q)$, and we can find $\tr_0 \in \N$ such that 
\be{}
\frac{1}{N_0} h\left( ([Q]_\tr)_{|_{\mathscr{F}_{N_0}}}
\mid q_{\rho,\nu}^{\otimes N_0} \right) \geq 
\frac{1}{N_0} h\left( Q_{|_{\mathscr{F}_{N_0}}}
\mid q_{\rho,\nu}^{\otimes N_0} \right) - \varepsilon  
\ee 
for $\tr \geq \tr_0$. Hence 
\be{eq:thm14b2}
H([Q]_\tr \mid q_{\rho,\nu}^{\otimes\N}) \geq H(Q \mid q_{\rho,\nu}^{\otimes\N}) - 2\varepsilon 
\quad \mbox{for $\tr\ge\tr_0$}. 
\ee
We may also assume that $[Q]_\tr \in \mathcal{O}(Q)$ for $\tr\geq\tr_0$. For a given $N \geq N_0$, 
pick $\tr(N) \geq \tr_0$ so large that $m_{[Q]_{\tr(N)}} H( \Psi_{[Q]_{\tr(N)}} \mid \nu^{\otimes \N}) 
\leq \delta_N/2$. Using the strategy described at the beginning of Section~\ref{S4}, we can 
construct a neighbourhood $\mathcal{O}_N \subset \mathcal{O}(Q)$ of $[Q]_{\tr(N)}$ such that 
the conditional probability $\Pr(R_N \in \mathcal{O}_N | X)$ is bounded below by 
\be{eq:thm14b3}
\exp\big[-N(H([Q]_\tr \mid q_{\rho,\nu}^{\otimes\N}) - \varepsilon) \big] 
\times \mbox{the cost of the first jump}, 
\ee
where the first jump takes us to a region of size $\approx N m_{[Q]_{\tr(N)}}$ on which  
the medium looks ``$\Psi_{[Q]_{\tr(N)}}$-typical''. Since, in a typical medium, the 
size of the first jump will be 
\be{eq:thm14b4}
\approx \exp\big[ N m_{[Q]_{\tr(N)}} H( \Psi_{[Q]_{\tr(N)}} \mid \nu^{\otimes \N}) \big] 
\leq \exp[ N \delta_N ], 
\ee
we obtain from (\ref{eq:thm14b1}) and (\ref{eq:thm14b2}--\ref{eq:thm14b4}) that 
\be{}
\Pr(R_N \in \mathcal{O}(Q) | X ) \ge 
\exp\big[ - N ( H(Q \mid q_{\rho,\nu}^{\otimes\N}) + 4 \varepsilon) \big] 
\ee
for $N$ large enough. 
\smallskip 

For the upper bound we can argue as follows: 
For $Q \in \mathcal{P}^{\mathrm{inv}}(\widetilde{E}^\N)$ put 
\be{}
r(Q) := \limsup_{\tr\to\infty} m_{[Q]_{\tr(N)}} H( \Psi_{[Q]_{\tr(N)}} \mid \nu^{\otimes \N}).  
\ee
Since $\rho$ satisfies the bound (\ref{eq:algtail}) for any 
$\alpha > 1$, we obtain from the upper bound in Theorem~\ref{mainthm} that 
the rate function at $Q$ is at least  
\be{} 
\limsup_{\tr\to\infty} I^{\mathrm{fin}}([Q]_\tr) = 
H(Q \mid q_{\rho,\nu}^{\otimes\N}) + (\alpha-1) r(Q), 
\ee
hence equals $\infty$ if $r(Q)>0$. On the other hand, if $r(Q)=0$, then this is 
simply the annealed bound.
\end{proof}


\section{Proof of Corollary\ \ref{LDPcount}}
\label{Scount}

\bpr
Let $E$ be a Polish space with metric $d_E$ (equipped with its 
Borel-$\sigma$-algebra $\mathscr{B}_E$). We can choose a sequence of 
nested finite partitions $\mathscr{A}_c=\{A_{c,1},\dots,A_{c,n_c}\}$, 
$c\in\N$, of $E$ with the property that 
\be{eq:finepart}
\forall \, x \in E : \: \lim_{c\to\infty} \mathrm{diam}\big( \langle x \rangle_c \big) = 0,
\ee
where the \emph{coarse-graining map} $\langle \cdot \rangle_c$ maps an 
element of $E$ to the element of $\mathscr{A}_c$ it is contained in. 
Each $\mathscr{A}_{c}=\langle E \rangle_{c}$ is a finite set, which we equip 
with the discrete metric $d_c$. Extend $\langle \cdot \rangle_c$ to 
$\langle E \rangle_{c'}$ for each $c'>c$ via $\langle A_{c',i'} \rangle_c 
= A_{c,i}$ if $A_{c',i'} \subset A_{c,i}$. Then the collection $\mathscr{A}_c$, 
$\langle \cdot \rangle_c$, $c \in \N$, forms a projective family, 
and the projective limit 
\be{}
F = \big\{ (\xi_1,\xi_2,\dots) : \xi_c \in \mathscr{A}_c, 
            \langle \xi_{c'} \rangle_c = \xi_c, 1 \le c < c' \big\}
\ee
is again a Polish space with the metric 
\be{}
d_{F}\big((\xi_1,\xi_2,\dots), (\eta_1,\eta_2,\dots)\big) 
:= \sum_{c=1}^\infty 2^{-c} d_c\big( \xi_c, \eta_c \big). 
\ee
We equip $F$ with its Borel-$\sigma$-algebra $\mathscr{B}_F$. 
We can identify $E$
with a subset of $F$ via $\iota : x \mapsto \big( \langle x \rangle_c
\big)_{c\in\N}$, since $\iota$ is injective by (\ref{eq:finepart}). 
Note that $\iota(E)$ is a measurable subset of $F$ 
(in general $\iota(E) \neq F$; it is easy to see that $\iota(E)$ 
is a closed subset of $F$ when $E$ is compact; for non-compact $E$ use the 
one-point compactification of $E$).

Note that the topology generated by
$d_{F}$ on $\iota(E)$ is finer than the original topology
generated by $d_E\,$: By
(\ref{eq:finepart}), for each $x \in E$ and $\varepsilon>0$, there is
an $\varepsilon'>0$ such that the $d_{F}$-ball of radius
$\varepsilon'$ around $x$ is contained in the $d$-ball of radius
$\varepsilon$ around $x$. We will make use of the fact that 
\be{eq:compatsigmaalgebras}
\mbox{the trace of $\mathscr{B}_F$ on $\iota(E)$ agrees with 
the image of $\mathscr{B}_E$ under $\iota$.} 
\ee
To check this, note that 
for any $x \in E$, the function 
\be{}
F \ni \xi \mapsto 
\begin{cases} d_E\big(\iota^{-1}(\xi), x\big), & \xi \in \iota(E), \\
\infty, & \mbox{otherwise}, 
\end{cases} 
\ee
can be pointwise approximated by functions that are constant on 
$\iota(A_{c,i})$, $i=1,\dots,n_c$, and is therefore 
$\mathscr{B}_F$-measurable.  

We extend $\langle \cdot \rangle_c$ in the obvious way to $E^N$
and $\widetilde{E}^N$, $N \in \N \cup \{\infty\}$
(via coordinate-wise coarse-graining), and then to $\mathcal{P}(E^N)$,
$\mathcal{P}(\widetilde{E}^N)$, $N\in\N$, and finally to 
$\mathcal{P}^{\mathrm{inv}}(E^\N)$ and
$\mathcal{P}^{\mathrm{inv}}(\widetilde{E}^\N)$ (by taking image
measures). Note that $\langle\cdot\rangle_{c}$ and $[\cdot]_\tr$
commute, and
\be{}
m_Q = m_{\langle Q \rangle_{c}}, \;\; 
\langle\Psi_Q\rangle_{c} =\Psi_{\langle Q \rangle_c}, \quad 
Q \in \mathcal{P}^{\mathrm{inv}}(\widetilde{E}^\N). 
\ee

By Theorem~\ref{mainthm}, for each $c\in\N$
the family 
\be{ag51}
\Pr(\langle R_N \rangle_\mathrm{c}\in\cdot\mid X), \quad N \in \N,
\ee
$X$-a.s.\ satisfies the LDP with deterministic rate function 
\be{}
I^\mathrm{que}_c(Q) = 
\begin{cases} 
I^\mathrm{fin}_c(Q) := 
H\big(Q \mid \langle q_{\rho,\nu}^{\otimes\N}\rangle_c \big) 
+ (\alpha-1) m_Q H\big(\Psi_Q \mid \langle \nu^{\otimes \N}\rangle_\mathrm{c} \big), 
& Q \in \mathcal{P}^{\mathrm{inv,fin}}(\langle\widetilde{E}\rangle_c^\N), \\[1ex]
\lim_{\tr\to\infty} I^\mathrm{fin}_c([Q]_\tr), & \mbox{if } m_Q = \infty. 
\end{cases}
\ee
Hence, by the Dawson-G\"artner \emph{projective limit theorem} 
(see Dembo and Zeitouni~\cite{DeZe98}, 
Theorem~4.6.1), the family $\Pr(R_N\in\cdot\mid X)$, $N\in\N$, $X$-a.s.\ satisfies the LDP on 
$\mathcal{P}^{\mathrm{inv}}(\widetilde{F}^\N)$ with rate function
\be{eq:ratefcttrunc}
I_F^\mathrm{que}(Q) =
\sup_{\mathrm{c} \in \N} I^\mathrm{que}_{c}(\langle Q\rangle_c), 
\qquad Q \in \mathcal{P}^{\mathrm{inv}}(\widetilde{F}^\N).
\ee

The following lemma follows from Deuschel and Stroock \cite{DeStr89}, 
Lemma~4.4.15. 
\begin{lemma}
\label{lemma:relentrlimit}
Let $G$ be a Polish space, let $\mathscr{A}_c=\{A_{c,1},\dots,A_{c,n_c}\}$, 
$c=1,2,\dots$ be a sequence of nested finite partitions of $G$ such that 
$\lim_{c\to\infty} \mathrm{diam}\big( \langle x \rangle_c \big) = 0$ for 
all $x \in G$ (with a coarse-graining map defined as above). 
Then we have, for 
$\mu, \nu \in \mathcal{P}(G)$, 
\be{eq:relentrlim}
h( \langle \mu \rangle_c \mid \langle \nu \rangle_c) 
\nearrow h(\mu \mid \nu) \quad \mbox{as $c \to \infty$}.
\ee
\end{lemma}

Let 
\begin{eqnarray}
\label{eq:defPEF}
\mathcal{P}_E^{\mathrm{inv}}(F^\N) & := & 
\left\{ \Phi \in \mathcal{P}^{\mathrm{inv}}(F^\N) : \pi_1\Phi(\iota(E))=1 \right\}, \\
\label{eq:defPEtildeF}
\mathcal{P}_E^{\mathrm{inv}}(\widetilde{F}^\N) & := & 
\left\{ Q \in \mathcal{P}^{\mathrm{inv}}(\widetilde{F}^\N) : \pi_1Q(\iota(\widetilde{E}))=1 
\right\}. 
\end{eqnarray}
Note that (\ref{eq:compatsigmaalgebras}) allows to view each $\Phi \in 
\mathcal{P}^{\mathrm{inv}}(E^\N)$ as an element of 
$\mathcal{P}_E^{\mathrm{inv}}(F^\N)$ and each 
$Q \in \mathcal{P}^{\mathrm{inv}}(\widetilde{E}^\N)$ as an element of 
$\mathcal{P}_E^{\mathrm{inv}}(\widetilde{F}^\N)$ via the 
identification of $E$ and $\iota(E) \subset F$. 
In particular, we can view $\nu^{\otimes\N}$ as an element of 
$\mathcal{P}_E^{\mathrm{inv}}({F}^\N)$ and 
$q_{\rho,\nu}^{\otimes\N}$ as an element of $\mathcal{P}_E^{\mathrm{inv}}(\widetilde{F}^\N)$. 
We will make use of the fact that, since each real-valued 
$d_E$-continuous function on $\iota(E)$ is automatically 
$d_F$-continuous, the weak topology on 
$\mathcal{P}_E^{\mathrm{inv}}(\widetilde{F}^\N)$ is finer than the 
weak topology on $\mathcal{P}^{\mathrm{inv}}(\widetilde{E}^\N)$. 
\smallskip

Fix $Q \in \mathcal{P}^{\mathrm{inv,fin}}(\widetilde{F}^\N)$. 
Note that the functions 
\be{eq:functincrbc}
(N,c) \mapsto \frac1N 
h\left( \langle \pi_N Q \rangle_c \mid \langle q_{\rho,\nu}^{\otimes N} \rangle_c \right)
\quad \mbox{and} \quad 
(L,c) \mapsto \frac1L
h\left( \langle \pi_L \Psi_Q \rangle_c \mid \langle \nu^{\otimes L} \rangle_c \right)
\ee
are non-decreasing in both coordinates. 
Then deduce from (\ref{eq:ratefcttrunc}) and (\ref{eqnratefctexplicit}) that 
\begin{eqnarray} \label{eq:IfinPolish1}
I^\mathrm{que}_F(Q) & = & 
\sup_{\mathrm{c} \in \N} 
\left\{ H(\langle Q \rangle_c \mid \langle q_{\rho,\nu}^{\otimes\N}\rangle_c )
+ (\alpha-1) \, m_{\langle Q \rangle_c} 
\,H(\Psi_{\langle Q\rangle_c} \mid \langle \nu^{\otimes \N}\rangle_c)
\right\} \nonumber \\
& = & 
\sup_{\mathrm{c} \in \N} 
\left\{ \sup_{N \in \N} \frac1N 
h\left( \langle \pi_N Q \rangle_c \mid \langle q_{\rho,\nu} \rangle_c^{\otimes N} \right) 
+ (\alpha-1) \, m_{Q} \sup_{L \in \N} \frac1L 
h\left( \langle \pi_L \Psi_Q \rangle_c \mid \langle \nu^{\otimes L} \rangle_c \right)
\right\} \nonumber \\
& = & 
\sup_{N \in \N} \frac1N \sup_{\mathrm{c} \in \N} 
h\left( \langle \pi_N Q \rangle_c \mid \langle q_{\rho,\nu}^{\otimes N} \rangle_c \right) 
+ (\alpha-1) \, m_{Q} \sup_{L \in \N} \frac1L \sup_{\mathrm{c} \in \N} 
h\left( \langle \pi_L \Psi_Q \rangle_c \mid \langle \nu^{\otimes L} \rangle_c \right) 
\nonumber \\
& = & 
\sup_{N \in \N} \frac1N 
h\left( \pi_N Q \mid q_{\rho,\nu}^{\otimes N} \right) 
+ (\alpha-1) \, m_{Q} \sup_{L \in \N} \frac1L 
h\left( \pi_L \Psi_Q \mid \nu^{\otimes L} \right) \nonumber \\
& = & H(Q \mid q_{\rho,\nu}^{\otimes\N}) + (\alpha-1) \, m_{Q} 
\,H(\Psi_{Q} \mid \nu^{\otimes \N}),
\end{eqnarray}
where we have used Lemma~\ref{lemma:relentrlimit} in the fourth line. 
Note that in the third line interchanging the suprema and
splitting out the supremum over the sum is justified
because of (\ref{eq:functincrbc}). 

For $Q \in \mathcal{P}^{\mathrm{inv}}(\widetilde{F}^\N)$ with 
$m_Q=\infty$ we see from 
(\ref{eq:ratefcttrunc}), (\ref{eqgndefinitionI}) and (\ref{eq:IfinPolish1}) that 
\begin{eqnarray}
I^\mathrm{que}_F(Q) & = & \sup_{\mathrm{c} \in \N} I^\mathrm{que}_{c}(\langle Q\rangle_c) 
= \sup_{\mathrm{c} \in \N} \sup_{\tr \in \N} 
\left\{ H([\langle Q \rangle_c]_\tr \mid \langle q_{\rho,\nu}^{\otimes\N}\rangle_c) 
+ (\alpha-1) \, m_{[Q]_\tr} 
\,H(\langle \Psi_{[Q]_\tr} \rangle_c \mid \langle \nu^{\otimes \N}\rangle_c)
\right\} \nonumber \\
& = & 
\sup_{\tr \in \N} \left\{ 
H([Q]_\tr \mid q_{\rho,\nu}^{\otimes\N}) + (\alpha-1) \, m_{[Q]_\tr} 
\, H(\Psi_{[Q]_\tr} \mid \nu^{\otimes \N})
\right\}. 
\end{eqnarray}
Note that $\sup_{\tr\in\N}$ can be replaced by $\tr\to\infty$ by 
arguments analogous to Step~3 in the proof of Theorem~\ref{mainthm}.

Finally, we transfer the LDP from $\mathcal{P}^{\mathrm{inv}}(\widetilde{F}^\N)$ 
to $\mathcal{P}^{\mathrm{inv}}(\widetilde{E}^\N)$. To this end, we 
first verify that the rate function is concentrated on 
$\mathcal{P}_E^{\mathrm{inv}}(\widetilde{F}^\N)$. 
Put 
\be{}
F'':=\{ y \in \widetilde{F} : 
\mbox{$y$ contains at least one letter from $F \setminus \iota(E)$}\}.
\ee
Then $q_{\rho,\nu}(F'')=0$. For 
$Q \in \mathcal{P}^{\mathrm{inv}}(\widetilde{F}^\N) \setminus 
\mathcal{P}_E^{\mathrm{inv}}(\widetilde{F}^\N)$ we have $\pi_1 Q(F'')>0$, 
and hence 
\be{}
I^{\mathrm{que}}_F(Q) \ge H(Q \mid q_{\rho,\nu}^{\otimes\N}) 
\ge h(\pi_1 Q \mid q_{\rho,\nu}) = \infty.
\ee 
Thus, by Dembo and Zeitouni\ \cite{DeZe98}, Lemma~4.1.5, 
the family $\Pr(R_N\in\cdot\mid X)$ satisfies for $\nu^{\otimes\N}$-a.s.\ 
all $X$ an LDP on $\mathcal{P}_E^{\mathrm{inv}}(\widetilde{F}^\N)$ with 
rate $N$ and with rate function given by 
(\ref{eqgndefinitionI}--\ref{eqnratefctexplicit}). 

To conclude the proof, observe that we can identify 
$\mathcal{P}^{\mathrm{inv}}(\widetilde{E}^\N)$ and 
$\mathcal{P}_E^{\mathrm{inv}}(\widetilde{F}^\N)$, and that the weak topology on 
$\mathcal{P}^{\mathrm{inv}}(\widetilde{E}^\N)$, which is `built' on $d_E$, 
is not finer than that which $\mathcal{P}_E^{\mathrm{inv}}(\widetilde{F}^\N)$ 
inherits from $\mathcal{P}^{\mathrm{inv}}(\widetilde{F}^\N)$, which is 
`built' on $d_F$ (recall the discussion following 
(\ref{eq:defPEF}--\ref{eq:defPEtildeF})). 
Consequently, the LDP carries over. 
\epr


\appendix

\section{Appendix: Continuity under truncation limits}
\label{SA}

The following lemma implies (\ref{conlimtr}).

\bl{trunclimcont} 
\label{lemma:appendix}
For all $Q \in \mathcal{P}^{\mathrm{inv,fin}}(\widetilde{E}^\N)$,
\be{trunclimconts}
\begin{aligned}
\lim_{\tr\to\infty} H([Q]_\tr \mid q_{\rho,\nu}^{\otimes\N}) 
&= H(Q \mid q_{\rho,\nu}^{\otimes\N}),\\
\lim_{\tr\to\infty} m_{[Q]_\tr} H(\Psi_{[Q]_\tr} \mid \nu^{\otimes\N}) 
&= m_{Q} H(\Psi_{Q} \mid \nu^{\otimes\N}). 
\end{aligned}
\ee
\end{lemma}

\begin{proof} 
The proof is not quite standard, because $Q$ and $[Q]_\tr$, respectively, $\Psi_Q$
and $\Psi_{[Q]_\tr}$ are not ``$\bar{d}$-close'' when $\tr$ is large, so that we 
cannot use the fact that entropy is ``$\bar{d}$-continuous'' (see Shields~\cite{Shi96}).

Lower semi-continuity yields $\liminf_{\tr\to\infty} {\rm l.h.s.} \geq {\rm r.h.s.}$
for both limits, so we need only prove the reverse inequality. Note that, for all 
$Q \in \mathcal{P}^{\mathrm{inv,fin}}(\widetilde{E}^\N)$, 
\be{Hbdsfin}
H(Q) \le h(Q_{|_{\mathscr{F}_1}}) \leq 
h\big( \LL_Q(\tau_1)\big) + m_{Q} \log |E| < \infty, 
\quad H(\Psi_Q) \leq \log |E| < \infty, 
\quad H(Q \mid q_{\rho,\nu}^{\otimes\N}) < \infty. 
\ee
For $Z$ a random variable, we write $\LL_Q(Z)$ to denote the law of $Z$ under $Q$.

\subsection{Proof of first half of (\ref{trunclimconts})}
\label{appA1}

\begin{proof}
Since $q_{\rho,\nu}^{\otimes\N}$ is a product measure, we have for, any $\tr \in \N$, 
\begin{equation}
\begin{aligned}
H([Q]_\tr \mid q_{\rho,\nu}^{\otimes\N}) 
&= -H([Q]_\tr) - \E_{[Q]_\tr}\left[ \log \rho(\tau_1)\right] 
- \E_{[Q]_\tr}\left[ \sum_{i=1}^{\tau_1} \log \nu\left(Y^{(1)}_i\right)\right]\\
&= -H([Q]_\tr) - \E_{Q}\left[\log \rho(\tau_1 \wedge\tr)\right] 
- \E_{Q}\left[ \sum_{i=1}^{\tau_1 \wedge \tr} \log \nu\left(Y^{(1)}_i\right)\right].
\end{aligned}
\end{equation}
By dominated convergence, using that $m_Q<\infty$ and $\log\rho(n)\leq C\log(n+1)$ for
some $C<\infty$, we see that as $\tr\to\infty$ the last two terms in the second line 
converge to 
\be{ag47}
- \E_{Q}\big[ \log\rho(\tau_1)\big] 
- \E_{Q}\left[ \sum_{i=1}^{\tau_1} \log \nu\left(Y^{(1)}_i\right)\right].
\ee
Thus, it remains to check that 
\begin{equation}
\label{eq:HQtlim}
\lim_{\tr\to\infty} H([Q]_\tr) = H(Q). 
\end{equation}

Obviously, $H([Q]_\tr) \le H(Q)$ for all $\tr\in\N$ (indeed, $h({[Q]_\tr}_{|_{\mathscr{F}_N}}) 
\leq h(Q_{|_{\mathscr{F}_N}})$ for all $N,\tr\in\N$, because $[Q]_\tr$ is the image measure of 
$Q$ under the truncation map). For the asymptotic converse, we argue as follows. A decomposition 
of entropy gives  
\begin{equation}
\label{eq:hNdecomp} 
h(Q_{|_{\mathscr{F}_N}}) 
= h({[Q]_\tr}_{|_{\mathscr{F}_N}}) 
+ \int_{[\widetilde{E}]_\tr^N} 
h\Big( \mathscr{L}_Q\big(\pi_N Y \mid  
\pi_N [Y]_\tr = z\big) \Big)\,(\pi_N [Q]_\tr)(dz), 
\end{equation}
where $\pi_N$ is the projection onto the first $N$ words, and $\mathscr{L}_Q(\pi_N Y \mid 
\pi_N [Y]_\tr = z)$ is the conditional distribution of the first $N$ words given their 
truncations. We have 
\be{ag2} 
h\Big( \mathscr{L}_Q\big(\pi_N Y \mid \pi_N [Y]_\tr = z\big) \Big) 
\leq \sum_{i=1}^N h\Big( \mathscr{L}_Q\big(Y_i \mid \pi_N [Y]_\tr = z\big) \Big) 
\ee
and
\begin{equation}
\begin{aligned} 
\label{ag3}
\int_{[\widetilde{E}]_\tr^N} 
&h\Big( \mathscr{L}_Q\big(Y_i \mid \pi_N [Y]_\tr = z\big) \Big)
\,(\pi_N [Q]_\tr)(dz) \\
&\leq 
\int_{[\widetilde{E}]_\tr^N} 
h\Big( \mathscr{L}_Q\big(Y_i \mid [Y_i]_\tr = z_i\big) \Big) 
\,(\pi_N [Q]_\tr)(dz) \\
&= \int_{[\widetilde{E}]_\tr} 
h\Big( \mathscr{L}_Q\big(Y_1 \mid [Y_1]_\tr = y\big) \Big) 
\,(\pi_1 [Q]_\tr)(dy), \qquad 1 \leq i \leq N,
\end{aligned}
\end{equation}
where the inequality in the second line comes from the fact that conditioning on 
less increases entropy, and the third line uses the shift-invariance. Combining 
(\ref{eq:hNdecomp}--\ref{ag3}) and letting $N\to\infty$, we obtain
\be{ag4} 
H(Q) \leq H([Q]_\tr) + \int_{[\widetilde{E}]_\tr} 
h\Big( \mathscr{L}_Q\big(Y_1 \mid [Y_1]_\tr = y\big) \Big) 
\,(\pi_1 [Q]_\tr)(dy), 
\ee
and so it remains to check that the second term in the right-hand side vanishes 
as $\tr\to\infty$. 

Note that this term equals (write $\varepsilon$ for the empty word and $w \cdot w'$ 
for the concatenation of words $w$ and $w'$) 
\begin{equation}
\label{ag5} 
\begin{aligned}
&-\sum_{w \in \tilde{E} \atop \tau(w)=\tr} [Q]_\tr(w) 
\sum_{w' \in \tilde{E} \cup \{\varepsilon\}} \frac{Q(w\cdot w')}{[Q]_\tr(w)} 
\log \left[\frac{Q(w\cdot w')}{[Q]_\tr(w)}\right]\\
&= - \sum_{w'' \in \tilde{E} \atop \tau(w'') \ge \tr} Q(w'') \log Q(w'') 
+ \sum_{w'' \in \tilde{E} \atop \tau(w'') \ge \tr} Q(w'') 
\log\,[Q]_\tr([w'']_\tr). 
\end{aligned}
\end{equation}
But
\be{ag6} 
0 \geq \sum_{ {w'' \in \tilde{E}} \atop {\tau(w'') \ge \tr} } 
Q(w'') \log\,[Q]_\tr([w'']_\tr)
\geq \sum_{ {w'' \in \tilde{E} } \atop {\tau(w'') \ge \tr} } 
Q(w'') \log Q(w''),
\ee
and so the right-hand side of (\ref{ag5}) vanishes as $\tr\to\infty$.
\end{proof}

\subsection{Proof of second half of (\ref{trunclimconts})}
\label{appA2}

Note that $\lim_{\tr\to\infty} m_{[Q]_\tr}=m_Q$ and ${\rm w}-\lim_{\tr\to\infty} 
\Psi_{[Q]_\tr} = \Psi_{Q}$ by dominated convergence, implying that 
\be{hups}
\liminf_{\tr\to\infty} H(\Psi_{[Q]_\tr} \mid \nu^{\otimes\N}) 
\geq H(\Psi_{Q} \mid \nu^{\otimes\N}).
\ee 
So it remains to check the reverse inequality. Since $\nu^{\otimes\N}$ is product measure,
we have 
\be{ag7} 
H(\Psi_{[Q]_\tr} \mid \nu^{\otimes\N}) 
= - H(\Psi_{[Q]_\tr}) 
- \frac{1}{m_{[Q]_\tr}} \E_Q\left[ \sum_{i=1}^{\tau_1 \wedge \tr} 
\log \nu\left(Y^{(1)}_i\right) \right]. 
\ee
By dominated convergence, as $\tr\to\infty$ the second term converges to 
\be{ag8}
\frac{1}{m_{Q}} \E_Q\left[\sum_{i=1}^{\tau_1} \log \nu\left(Y^{(1)}_i\right) \right] 
= \int_E \Psi_Q(dx)\,\log\nu(x). 
\ee
Thus, it remains to check that 
\begin{equation} 
\label{eq:entropylimit}
\lim_{\tr\to\infty} H(\Psi_{[Q]_\tr}) = H(\Psi_Q). 
\end{equation} 
We will first prove (\ref{eq:entropylimit}) for ergodic $Q$, in which case 
$[Q]_\tr$, $\Psi_Q$, $\Psi_{[Q]_\tr}$ are ergodic (Birkner~\cite{Bi08}, Remark~5). 

For $\Psi \in \mathcal{P}^{\mathrm{erg}}(E^{\N})$ and $\varepsilon \in (0,1)$, 
let 
\be{}
\mathcal{N}_n(\Psi, \varepsilon) = 
\min\big\{ \# A \colon\, A \subset E^n, \Psi(A \times E^\infty) \ge \varepsilon \big\}
\ee
be the $(n,\varepsilon)$ covering number of $\Psi$. For any $\varepsilon \in (0,1)$,
we have 
\be{eq:coveringentropy}
\lim_{n\to\infty} \frac{1}{n} \log \mathcal{N}_n(\Psi, \varepsilon) = H(\Psi) 
\ee
(see Shields \cite{Shi96}, Theorem~I.7.4). The idea behind (\ref{eq:entropylimit}) 
is that there are $\approx \exp[n H(\Psi_Q)]$ ``$\Psi_Q$-typical'' sequences of 
length $n$, and that a ``$\Psi_{[Q]_\tr}$-typical'' sequence arises from a 
``$\Psi_Q$-typical'' sequence by eliminating a fraction $\delta_\tr$ of the 
letters, where $\delta_\tr\to 0$ as $\tr\to\infty$. Hence $\mathcal{N}_n(\Psi_Q,\varepsilon)$ 
cannot be much larger than $\mathcal{N}_n(\Psi_{[Q]_\tr}, \varepsilon)$ (on an exponential 
scale), implying that $H(\Psi_Q)-H(\Psi_{[Q]_\tr})$ must be small. 

To make this argument precise, fix $\varepsilon > 0$ and pick $N_0$ so large that 
\be{eq:kappalarge1}
Q\big( |\kappa(Y^{(1)},\dots,Y^{(N)})| 
\in N m_Q [1-\varepsilon,1+\varepsilon]\big) > 1-\varepsilon 
\qquad \mbox{for $N \ge N_0$}.
\ee
Pick $\tr_0 \in \N$ so large that for $\tr \ge \tr_0$ and $N \ge N_0$,
\be{eq:fewtruncations1}
Q\Big( {\textstyle \sum_{i=1}^N (\tau_i-\tr)_+} < N\varepsilon \Big) 
> 1-\varepsilon/2, \; Q\big(\tau_1 \leq \tr \big) > 1-\varepsilon/2, \;
m_{[Q]_\tr} > (1-\varepsilon) m_Q. 
\ee
For $n \ge \lceil N_0/m_Q\rceil$, we will construct a set $B \subset E^n$ such that 
\be{eq:BPsiQtypical}
\Psi_Q(B \times E^\infty) \ge \tfrac12, \qquad 
|B| \leq \exp\big[ n (H(\Psi_{[Q]_\tr})+\delta) \big], 
\ee
where $\delta$ can be made arbitrarily small by choosing $\varepsilon$ small in 
(\ref{eq:kappalarge1}--\ref{eq:fewtruncations1}). Hence, by the asymptotic cover 
property (\ref{eq:coveringentropy}), we have $H(\Psi_Q) \le (1+\delta) H(\Psi_{[Q]_\tr})$
and 
\be{HPsitrunclowerbound} 
\liminf_{\tr\to\infty} H(\Psi_{[Q]_\tr}) \ge H(\Psi_Q), 
\ee
completing the proof of (\ref{eq:entropylimit}).

We verify (\ref{eq:BPsiQtypical}) as follows. Put $N:= \lceil n m_Q (1+2\varepsilon)\rceil$. 
By (\ref{eq:kappalarge1}--\ref{eq:fewtruncations1}) and the asymptotic cover property 
(\ref{eq:coveringentropy}) for $\Psi_{[Q]_\tr}$, there is a set $A \subset \widetilde{E}^N$ 
such that 
\be{hups2}
\E_Q\big[\tau_1 \1_A(Y^{(1)},\dots,Y^{(N)}) \big] > (1-\varepsilon) m_Q
\ee 
and
\be{hups3}
\begin{aligned}
& |\kappa(y^{(1)},\dots,y^{(N)})| \ge n(1+\varepsilon), \quad 
\tau(y^{(1)}) \le \tr, \quad
\sum_{i=1}^N (\tau(y^{(i)})-\tr)_+ < N\varepsilon,\\
&\forall\,(y^{(1)},\dots,y^{(N)}) \in A,
\end{aligned}
\ee 
while the set  
\be{hups4}
B' := \big\{ \kappa([y^{(1)}]_\tr,\dots,
[y^{(N)}]_\tr)|_{(0,\lceil (1-\varepsilon)n\rceil]}\colon\, (y^{(1)},\dots,y^{(N)}) \in A \big\}
\subset E^{\lceil (1-\varepsilon)n \rceil]} 
\ee 
satisfies 
\be{eq:B'size}
|B'| \le \exp\big[ n (H(\Psi_{[Q]_\tr}) + \varepsilon) \big].
\ee

Put 
\be{hups5} 
B:= \big\{ \kappa(y^{(1)},\dots,y^{(N)})|_{(0,n]}\colon 
(y^{(1)},\dots,y^{(N)}) \in A \big\}\subset E^n. 
\ee
Observe that each $x' \in B'$ corresponds to at most 
\be{ag48}
|E|^{\varepsilon n} {n \choose \varepsilon n} \leq 
\exp\big[-n (\varepsilon \log \varepsilon + (1-\varepsilon) \log (1-\varepsilon))
+ n \varepsilon \log |E| \big] 
\ee
different $x \in B$, so that 
\be{eq:BonB'bound} 
|B| \leq 
|B'| \exp\big[ -n (\varepsilon \log \varepsilon + (1-\varepsilon) \log (1-\varepsilon))
+ n \varepsilon \log |E| \big]. 
\ee
We have
\be{ag9}
\begin{aligned} 
m_Q \Psi_Q(B \times E^\infty) 
& \geq \E_Q\left[ \sum_{k=0}^{\tau_1-1} \1_{B\times E^\infty} 
\big(\theta^k\kappa(Y)\big) \1_A(Y^{(1)},\dots,Y^{(N)})\right] \\
&= \E_Q\left[ \sum_{k=0}^{\tau_1 \wedge \tr -1} \1_{B'\times E^\infty}
\big(\theta^k\kappa([Y]_\tr)\big) \1_A(Y^{(1)},\dots,Y^{(N)})\right] \\
&\geq \E_Q\left[ \sum_{k=0}^{\tau_1 \wedge \tr-1} \1_{B'\times E^\infty}
\big(\theta^k\kappa([Y]_\tr)\big) \right] - \varepsilon m_Q \\
&=  m_{[Q]_\tr} \Psi_{[Q]_\tr}(B' \times E^\infty) - \varepsilon m_Q, 
\end{aligned}
\ee
so that, finally, 
\be{eq:PsiQBlarge} 
\Psi_Q(B \times E^\infty) \geq \frac{m_{[Q]_\tr}}{m_{Q}}\,
\Psi_{[Q]_\tr}(B' \times E^\infty) - \varepsilon \geq \tfrac12. 
\ee
Combining (\ref{eq:B'size}), (\ref{eq:BonB'bound}) and (\ref{eq:PsiQBlarge}),
we obtain (\ref{eq:BPsiQtypical}) with 
\be{ag49}
\delta = -\big(\varepsilon \log \varepsilon + (1-\varepsilon) \log (1-\varepsilon)\big)
+ \varepsilon \big(1 + \log |E| \big). 
\ee
Since $\limsup_{\tr\to\infty} H(\Psi_{[Q]_\tr}) \le H(\Psi_Q)$ by upper semi-continuity 
of $H$ (see e.g.\ Georgii \cite{Ge88}, Proposition.~15.14), this concludes the proof of 
(\ref{eq:entropylimit}) for ergodic $Q$.

For general $Q \in \mathcal{P}^{\mathrm{inv,fin}}(\widetilde{E}^\N)$, we recall the 
ergodic decomposition formulas stated in (\ref{ergdecomp}--\ref{mQPsiQdecomp}).
These yield 
\be{hups8}
\Psi_{[Q]_\tr} = \int_{\mathcal{P}^{\mathrm{erg,fin}}(\widetilde{E}^\N)} 
\frac{m_{[Q']_\tr}}{m_{[Q]_\tr}}\,\Psi_{[Q']_\tr}\,W_Q(dQ') 
\ee
and
\be{eq:HPsiQdecomp}
H(\Psi_{[Q]_\tr}) = 
\int_{\mathcal{P}^{\mathrm{erg,fin}}(\widetilde{E}^\N)} 
\frac{m_{[Q']_\tr}}{m_{[Q]_\tr}}\,H(\Psi_{[Q']_\tr})\,W_Q(dQ'),
\ee
because specific entropy is affine. The integrand inside (\ref{eq:HPsiQdecomp}) 
is non-negative and, by the above, converges to $\frac{m_{Q'}}{m_Q} H(\Psi_{Q'})$ as 
$\tr\to\infty$. Hence, by Fatou's lemma, 
\be{ag50}
\liminf_{\tr\to\infty} H(\Psi_{[Q]_\tr}) \ge 
\int_{\mathcal{P}^{\mathrm{erg,fin}}(\widetilde{E}^\N)} 
\frac{m_{Q'}}{m_Q}\,H(\Psi_{Q'}) \,W_Q(dQ') = H(\Psi_Q), 
\ee
which concludes the proof.
\end{proof}




\begin{thebibliography}{AA} 

\bibitem{Bi03}
M.\ Birkner, 
\emph{Particle Systems with Locally Dependent Branching:
Long-Time Behaviour, Genealogy and Critical Parameters}, 
Dissertation, Johann Wolfgang Goethe-Universit\"{a}t 
Frankfurt am Main, 2003.\\
\url{http://publikationen.ub.uni-frankfurt.de/volltexte/2003/314/}

\bibitem{Bi08} 
M.\ Birkner, 
Conditional large deviations for a sequence of words, 
Stoch.\ Proc.\ Appl.\ 118 (2008) 703--729. 

\bibitem{BGdH08} M.\ Birkner, A.\ Greven, F.\ den Hollander, 
Large deviations for the collision local time of transient 
random walks and intermediate phases in interacting stochastic 
systems, preprint 2008.

\bibitem{Co89} 
F.\ Comets, 
Large deviation estimates for a conditional probability distribution. 
Applications to random interaction Gibbs measures,
Probab.\ Theory Relat.\ Fields 80 (1989) 407--432.

\bibitem{DeZe98} A.\ Dembo and O.\ Zeitouni, 
\emph{Large Deviations Techniques and Applications}, 2nd Edition, 
Springer, 1998.

\bibitem{DeStr89}
J.-D.\ Deuschel and D.W.\ Stroock, 
\emph{Large Deviations}, Academic Press, Boston, 1989. 

\bibitem{Ge88}
H.-O.\ Georgii, 
\emph{Gibbs Measures and Phase Transitions}, 
de Gruyter Studies in Mathematics 9, Walter de Gruyter, Berlin, 1988.

\bibitem{dHo00}
F.\ den Hollander, \emph{Large Deviations}, Fields Institute Monographs 14,
American Mathematical Society, Providence, RI, 2000.

\bibitem{HoLi81}
R.\ Holley and T.\ M.\ Liggett, 
Generalized potlatch and smoothing processes, 
Z.\ Wahrscheinlichkeitstheorie verw.\ Gebiete 55 (1981) 165-195. 

\bibitem{Shi96} 
P.C.\ Shields, 
\emph{The Ergodic Theory of Discrete Sample Paths}, 
American Mathematical Society, Providence, RI, 1996. 

\bibitem{Str65} 
V.\ Strassen, 
The existence of probability measures with given marginals,
Ann.\ Math.\ Statist.\ 36 (1965) 423--439.

\bibitem{To07} 
F.\ Toninelli, 
Disordered pinning models and copolymers: beyond annealed bounds, 
Ann.\ Appl.\ Probab.\ 18 (2008) 1569-1587. 

\bibitem{Wa82} 
P.\ Walters, 
\emph{An Introduction to Ergodic Theory}, 
Graduate Texts in Mathematics 79, Springer, New York, 1982. 

\end{thebibliography}
\end{document}